\documentclass[10pt]{article}

\usepackage[english,strings]{babel}     
\usepackage[utf8]{inputenc}    
\usepackage[T1]{fontenc}       
\usepackage{longtable}         
\usepackage{exscale}           
\usepackage[sort]{cite}        
\usepackage{array}             
\usepackage{fancyhdr}          
\usepackage[a4paper]{geometry} 
\usepackage{xspace}            
\usepackage{tikz-cd}           
\usetikzlibrary{arrows,calc,  
through,backgrounds,matrix,
	decorations.pathmorphing,
	positioning,babel} 

\usepackage{nchairx}           
\usepackage[sectionbib         
          ]{chapterbib}        
\usepackage{appendix}
\usepackage[expansion=false    
           ]{microtype}        
\usepackage[nottoc]{tocbibind} 
\usepackage[backref=page,      
           hypertexnames=false,
           final=true,         
           pdfpagelabels       
           ]{hyperref}         
\usepackage{titlesec}
\usepackage{tocloft}

\titleformat{\subsection}[runin]%
{\normalfont\normalsize\bfseries}{\thesubsection}{1em}{}

\newcommand{\cE}{\mathcal{E}}
\newcommand{\cF}{\mathcal{F}}

\newcommand{\cI}{\mathcal{I}}
\newcommand{\cJ}{\cI}

\newcommand{\cL}{\mathcal{L}}
\newcommand{\cM}{\mathcal{M}}
\newcommand{\cN}{\mathcal{N}}
\newcommand{\cO}{\mathcal{O}}

\newcommand{\fX}{\mathfrak{X}}

\newcommand{\bN}{\mathbb{N}}
\newcommand{\bR}{\mathbb{R}}
\newcommand{\bZ}{\mathbb{Z}}

\newcommand{\nE}{\boldsymbol{E}}

\newcommand{\vf}{\mathfrak{X}}
\newcommand{\can}{\mathrm{can}}

\DeclareMathOperator{\sym}{Sym} 
\DeclareMathOperator{\Tdim}{Totdim}

\newcommand{\chairxauthorbibfont}{\textsc}

\newcommand{\chairxtitlebibfont}{\textit}

\definecolor{forest}{rgb}{0.15,0.45,0.1}
\definecolor{purple}{rgb}{0.5,0,0.5}

\title{Deformations of Lagrangian $NQ$-submanifolds}
\author{
  \textbf{Miquel Cueca}\thanks{ \texttt{miquel.cuecaten@mathematik.uni-goettingen.de} Department of Mathematics,  University of G\"{o}ttingen; Bunsenstraße 3-5, 37073 Göttingen 
		(Germany) },
  \textbf{Jonas Schnitzer}\thanks{ \texttt{jonas.schnitzer@math.uni-freiburg.de} Department of Mathematics, University of Freiburg;
   Ernst-Zermelo-Strasse 1,   79104 Freiburg   (Germany)}
   }

\date{}
\begin{document}
\selectlanguage{english}

\maketitle
\begin{abstract}
In this paper we prove graded versions of the Darboux Theorem and 
Weinstein's Lagrangian tubular neighbourhood Theorem in order to study the deformation theory 
of Lagrangian $NQ$-submanifolds of degree $n$ symplectic $NQ$-manifolds. 
Using Weinstein's Lagrangian tubular neighbourhood Theorem, we attach to every 
Lagrangian $NQ$-submanifold an $L_\infty$-algebra, which controls its deformation theory. 
The main examples are coisotropic submanifolds of Poisson manifolds and (higher) Dirac structures with support
in (higher) Courant algebroids. 
\end{abstract}

\tableofcontents

\section{Introduction}

On the one hand, the study of higher structures is one of the key aspects of geometry
and mathematical physics in the 21st century. In particular, higher versions of symplectic
geometry have emerged naturally in research connected to field theories \cite{AKSZ, mnev:bv, scw:bv},
mechanics \cite{rog:cat, zam:hom}  algebraic geometry \cite{ptvv, pym:shif} or Lie
theory \cite{get:shi, alan:grou} among others. Due to their crucial importance, 
different communities have been developing non-equivalent frameworks to deal with higher 
symplectic structures. In this note we work with $\bN$-graded $Q$-manifolds ($NQ$-manifolds for short) and symplectic 
structures therein, see e.g. \cite{AKSZ, Roytenberg, Severa}.   

On the other hand, deformation theory focuses on the infinitesimal study of moduli spaces \cite{art:ver}.  
The classical philosophy states that every deformation problem is governed by a differential graded 
Lie algebra (dgLa) via solutions of the Maurer-Cartan equation modulo gauge action \cite{nij:lie, sta:lie}. 
However, in many situations the moduli problems form a higher stack, see \cite{hin:for, lur:mod, pri:def}. 
Therefore, one is forced to consider Maurer-Cartan elements modulo gauge on $L_\infty$-algebras, i.e. the 
\emph{up to homotopy} versions of dgLas as we explain in Appendix \ref{app:Linfty}, see e.g. \cite{CattaneoFelder,
Gualtieri, oh:coi} for concrete examples.

In the present paper we combine the above topics and  study the deformation problem of Lagrangian 
$NQ$-submanifolds inside  symplectic $NQ$-manifolds. Our main achievements are: 
\begin{itemize}
    \item We prove a graded version of the Darboux Theorem, see Corollary
    \ref{Cor: Darboux}. 
    \item We prove a graded version of Weinstein's Lagrangian tubular neighbourhood theorem, see Theorem \ref{Thm: dgWeinstein}.
    \item For a Lagrangian $NQ$-submanifold and a tubular neighbourhood, we construct an $L_\infty$-algebra and show that
    two different tubular neighourhoods give isomorphic $L_\infty$-algebras, see Theorem  \ref{Thm: UniqLinfty}.  
    \item We show the formal deformation problem of Lagrangian 
    $NQ$-submanifolds inside symplectic $NQ$-manifolds is controlled by the above $L_\infty$-algebra, 
    see Theorem \ref{Thm:formDefProb}.
    \item We give a geometric characterization for the gauge action in Theorem \ref{Lem:HomPoissEqMC}. 
\end{itemize}
Let us recall in the following the classical \emph{non-graded} results, give a quick overview in our work, 
put it into mathematical context, and finally provide some applications and comment on future directions.   
\subsection*{The framework.} Recall that an \emph{$\bN$-graded manifold} $\cM=(M, C_\cM)$, see 
$\S$\ref{sec:GrMan} for a precise definition,  is a ringed space where the sheaf of functions is 
enlarged with new polynomial variables concentrated in positive degrees that super-commute according 
to their degree. These provide a refinement of the \emph{supermanifolds} commonly used by physicists, 
see \cite{fio:sup, kos:gra, mnev:bv}. We say that $(\cM, Q)$ is an \emph{$NQ$-manifold} if in addition 
$C_\cM$ has a degree $1$ derivation, i.e. a vector field as we introduce in $\S$\ref{sec:vec}, such that
$$Q^2=\frac{1}{2}[Q, Q]=0.$$
Thus, $(C_\cM, Q)$ becomes a cochain complex and we are doing homological algebra, but with a geometric flavour.

The main characters of this work are \emph{symplectic $NQ$-manifolds} $(\cM, \omega, Q)$, see $\S$\ref{sec:sym}, 
and more concretely their \emph{Lagrangian $NQ$-submanifolds} $j:\cL\hookrightarrow(\cM, \omega, Q)$ as introduced 
in Definition \ref{def:lag}. In Table \ref{tab:1} we summarize some examples\footnote{In the second row $(\cM, \omega, Q)$ denotes an arbitrary degree $2$ symplectic $NQ$-manifold. The corresponding Courant algebroid is given by the degree $1$ functions of $\cM$, i.e.  $C_\cM^1=\Gamma E$. As explained in the references, this is a one to one correspondence.} and references that illustrate how 
Lagrangian $NQ$-submanifolds unify seemingly unrelated geometric structures previously considered by mathematicians.  
Although symplectic $NQ$-manifolds and their Lagrangian $NQ$-submanifolds have drawn much attention in recent years 
it is difficult to find a rigorous mathematical treatment in the literature. Therefore, we decided to include Section 
\ref{sec:Pre} to fill this gap. 
\begin{table}[h]
\[
\begin{array}{|c|c|c|c|}
\hline
  \text{Symplectic } Q\text{-manifold}&\text{Classical description} & \text{Lagrangian } Q\text{-subman.}& \text{Ref.}\\
\hline
\begin{array}{c}
     \text{Degree } 1  \\
     (T^*[1]M, \omega_{\can}, Q)
\end{array}
& 
\begin{array}{c}
     \text{Poisson manifold}  \\
     (M, \pi)
\end{array}& 
\begin{array}{c}
     \text{Coisotropic}  \\
     \text{submanifolds}
\end{array}
&\text{\cite{cat:coi, Vor:bia}}
\\
\hline
 \begin{array}{c}
     \text{Degree } 2  \\
     (\cM, \omega, Q)
\end{array}
&
 \begin{array}{c}
     \text{Courant algebroid}  \\
     (E\to M, \langle\cdot,\cdot\rangle, \rho, [\cdot, \cdot])
\end{array}
& 
\begin{array}{c}
     \text{Dirac structures}  \\
     \text{with support}
\end{array}
& \text{\cite{Roytenberg, Severa}}
\\ \hline
\begin{array}{c}
     \text{Degree } 3  \\
     (T^*[3]A[1], \omega_{\can}, Q)
\end{array}
& 
\begin{array}{c}
     H\text{-twisted Lie algebroid}  \\
     (A\to M, \rho, [\cdot,\cdot], \langle\cdot,\cdot\rangle)
\end{array}& 
\begin{array}{c}
     \text{Higher Dirac}  \\
     \text{structures}
\end{array}
& \text{\cite{CuecaCoTan, gru:ht}}
\\
\hline
\begin{array}{c}
     \text{Degree } 6  \\
     \scriptstyle(T^*[6]T[1]M\times \bR[3], \omega_{\can}\times\omega_{\bR[3]}, Q)
\end{array}
& 
\begin{array}{c}
     M_5\text{-theory algebroid}  \\
     (TM\oplus \wedge^2T^*M\oplus \wedge^5 T^*M, \langle\cdot,\cdot\rangle, [\cdot,\cdot])
\end{array}& 
\begin{array}{c}
     \text{No Lagrangian}  \\
     \text{submanifolds}
\end{array}
& \text{\cite{arv:bra}} 
\\
\hline

 \begin{array}{c}
     \text{Degree } n  \\
     (T^*[n]T[1]M, \omega_{\can}, Q)
\end{array}
&
 \begin{array}{c}
     \text{Higher Courant algebroid}  \\
     (TM\oplus \wedge^{n-1}T^*M\to M, \langle\cdot,\cdot\rangle, \rho, [\cdot, \cdot])
\end{array}
& 
\begin{array}{c}
     \text{Higher Dirac}  \\
     \text{structures}
\end{array}
&\text{\cite{CuecaCoTan}}
\\ \hline

\hline
\end{array}\]
\caption{Examples of Lagrangian $NQ$-submanifolds.\label{tab:1}}
\end{table}

\subsection*{The classical deformation problem.} Let $j:L\to (M,\omega)$ be a Lagrangian submanifold 
of the (non-graded) symplectic manifold $(M,\omega)$. The classical deformation problem around $L$ studies the (local) moduli space
of small deformations of $L$ by considering nearby Lagrangian submanifolds modulo Hamiltonian
isotopies. Weinstein's Lagrangian tubular neighbourhood proved in \cite{Weinstein}  shows that any neighbourhood 
 of $L$ in $M$ is symplectomorphic to $(T^*L, \omega_{\can})$. Therefore, it is enough to study the problem 
 of deforming the zero-section of the cotangent bundle. One easily shows that for 
 $$\alpha\in \Omega^1(L),\quad \graph(\alpha)\subseteq T^*L \quad \text{is Lagrangian}\quad \Leftrightarrow \quad \D\alpha=0$$ 
and the Hamiltonian isotopies are given by exact $1$-forms. Thus $H^1_{\mathrm{dR}}(L)$ gives the desired (local) 
moduli space. It is worth to notice that in this case the solution is an abelian Lie algebra, hence linear and commutative.  

\subsection*{Weinstein's Lagrangian tubular neighbourhood.} As in the classical problem, the first step to study deformations of 
Lagrangian $NQ$-submanifolds is to provide a graded version of Weinstein's
Lagrangian tubular neighbourhood in the context of $\bN$-graded symplectic 
manifolds (Theorem \ref{Thm: dgWeinstein}). Our proof follows a classical argument using a tubular 
neighbourhood for graded manifolds, see Theorem \ref{Thm: TubNei}.  It is important to notice that this version of  
Weinstein's Lagrangian tubular neighbourhood does not involve $Q$-structures, so later on we will consider cotangent bundles with canonical 
symplectic structures, but with abritrary homological vector fields. \vspace{3mm}

\subsection*{$L_\infty$-algebra and deformations.} From this point on, our problem diverges from the classical case. 
It was noticed in \cite{Roytenberg} that on an $\mathbb{N}$-graded manifold $\cL$ closed $1$-forms of positive degree are exact, more concretely if
$$\alpha\in \Omega^{1,k}(\cL)\quad\text{with}\  k\neq 0 \ \text{ and } \ \D\alpha=0 \ \text{ then }\ \alpha=\D f\ \text{ for some} \ f\in C^k_{\cL}.$$
Therefore we have that Lagrangian submanifolds near the zero-section in 
$(T^*[n]\cL, \omega_{\can})$, for $n\geq 1$, are in one to one correspondence with $f\in C^n_{\cM}.$ 
Thus Hamiltonian isotopies are not the correct notion of trivial deformation for symplectic 
$\bN$-manifolds.  However, once we add the
$Q$-structure, we will have a clearer picture using Deligne’s deformation theory. Let $Q$ be given by  
 $$Q=\{\theta,\cdot\}_{\can}\in \fX^{1,1}(T^*[n]\cL) \quad \text{with}\quad 0=\{\theta,\theta\}_{\can} \quad \text{for some}\quad\theta\in C^{n+1}_{T^*[n]\cL},$$ 
then we get an $L_\infty$-algebra $(C_\cL^\bullet[n-1], \{l^i\})$ just by realizing that 
 $(C^\bullet_{T^*[n]\cL}[n], C^\bullet_\cL[n], 0_\cL)$ forms a V-algebra, see Definition \ref{def:Valgebra}, 
 and $\theta$ satisfies the hypothesis of Theorem \ref{Them: DerBra}. The 
 $L_\infty$-structure $\{l^i\}$ appear then as derived 
 brackets and they only depend on the $\infty$-jet of $\theta$ in the vector bundle directions at the zero 
 section, see Equation \ref{Eq: LocLinfty}. For a general Lagrangian $NQ$-submanifold $j:\cL\hookrightarrow(\cM,\omega, Q),$
  we show in Theorem \ref{Thm: UniqLinfty} that the isomorphism class of the above $L_\infty$-algebra is independent of 
  the chosen tubular neighbourhood. In particular, this implies the main result in \cite{Gualtieri}.   

In $\S$\ref{sec:defq} we show that the formal deformation problem of $j:\cL\hookrightarrow(\cM,\omega, Q)$ 
as a Lagrangian $NQ$-submanifold is controlled by the Maurer-Cartan elements in the above $L_\infty$-algebra. 
This is done first by transforming the $NQ$-submanifold condition into a flow equation and then making a Taylor 
expansion, see Theorem \ref{Thm:formDefProb}. The aforementioned result unifies and generalizes the deformations 
of Coisotropic submanifolds \cite{CattaneoFelder, oh:coi}, presymplectic structures \cite{zam:pre}(see also
 \cite{zam:reg} for the case of regular Poisson structures) and Dirac structures \cite{Zambon, Gualtieri, wei:cou}. 

For a general Lagrangian $NQ$-submanifold $j:\cL\hookrightarrow(\cM,\omega, Q)$ the $L_\infty$-algebra 
$(C_\cL^\bullet[n], \{l^i\})$ has infinitely many brackets. However, if the bodies of $\cL$ 
and $\cM$ coincide then $(C_\cL^\bullet[n], \{l^i\})$ is in fact an $L_{n+1}$-algebra and
controls the smooth deformation problem
of $\cL$ in $(\cM,\omega, Q)$. In 
this case, $\cL$ becomes a \emph{homotopy Poisson manifold} as considered in \cite{Raj}, see Corollary 
\ref{Prop:Lagsubmanifold}. As explained in Remark \ref{rmk:sp} the above Corollary gives a concrete way 
of producing an \emph{$(n-1)$-shifted Poisson} structure for a given  \emph{$n$-shifted Lagrangian}, 
a principle in derived algebraic geometry \cite{saf:poi}. Another condition to ensure that we can get a hand on 
possible divergence issues is to ask $\theta$ to be an \emph{entire} function as introduced in \cite{zam:ent}.
We do not explore these analytic techniques here; however, we believe that all of them apply in our context.

To  identify the moduli space of deformations correctly,  we also need to incorporate the
gauge action
of the $L_\infty$-algebra $(C_\cL^\bullet[n-1], \{l^i\})$, see Definition \ref{prop: EquMC}. In the case where $\cL$ and $\cM$ have 
the same bodies, ensuring that no convergence issues appear, we are able to give a geometric characterization 
of the elements in the gauge orbit of zero. 
Indeed, we show in Theorem \ref{Lem:HomPoissEqMC} that the gauge action is given by flows of Hamiltonian 
vector fields for $Q$-exact functions. In other words, the homological vector field $Q$ plays the 
role of the de Rham differential in the classical case. Hence for a general Lagrangian 
$NQ$-submanifold $j:\cL\hookrightarrow(\cM,\omega, Q)$ we get that:
\begin{center}
        \begin{tabular}{|c|}
        \hline
             The $L_\infty$-algebra $(C_\cL^\bullet[n-1], \{l^i\})$ controls the formal deformations of $\cL$\\
             as a Lagrangian $NQ$-submanifold modulo flows of $Q$-exact Hamiltonian vector fields.\\
             \hline
        \end{tabular}
\end{center}
Finally, in $\S$\ref{sec:sim} we explore the simultaneous deformations of a Lagrangian $NQ$-submanifold 
together with the ambient homological vector field $Q$, see \cite{Zambon} for related works. 

\subsection*{Applications and further directions.} As shown in Table \ref{tab:1} Lagrangian 
$NQ$-submanifolds appear naturally in physics and Poisson geometry. Therefore their deformation 
theory has already a relevance in this context. To illustrate that, in $\S$\ref{sec:ex} we show 
that many computations and results which have been obtained in seemingly different contexts, can be reinterpreted 
as deformations of a Lagrangian $NQ$-submanifolds 
or give obstructions to their existence. In particular we include some relevant examples of Lagrangian 
$NQ$-submanifolds inside a degree $3$ symplectic $NQ$-manifold that were not treated previously.

Lagrangian $NQ$-submanifolds of $(\cM, \omega, Q)$ also appear in the theory of sigma models, more specifically in \emph{AKSZ} sigma models,
as they provide a reasonable notion of
boundary condition for the AKSZ sigma model with target $(\cM, \omega, Q)$,
see e.g. \cite{AKSZ, hof:top, ike:can, mnev:hol, sev:dua}.  Therefore, it would be interesting to know the behaviour of 
these field theories under deformations of the boundary conditions. 
Some of the results obtained in \cite{hof:top, ike:can} can be reinterpreted as deforming the boundary 
condition to obtain "twisted" models as we explain in Example \ref{ex:tca}. 

An obvious extension of our work would be to consider the deformation problem  of a Lagrangian 
$Q$-submanifolds inside $\bZ$-graded symplectic $Q$-manifolds. Indeed, this is the relevant case
for the BV and BFV formalisms. In particular, in the BV-formalism the gauge fixing condition is 
implemented as a Lagrangian $Q$-submanifold and the gauge independence translates to a deformation invariance of a certain functional, 
see \cite{mnev:bv, scw:bv}. We hope that our work might be a further step towards
a better understanding of the BV-formalism. 

A probably more subtle direction of extending our work is to consider the deformation problem of a shifted Lagrangian 
inside a shifted symplectic algebroid as introduced in \cite{pym:shif}. These are the \emph{up to homotopy} 
versions of Lagrangian $NQ$-submanifolds inside  symplectic $NQ$-manifolds.
One is tempted to say that the solution to this deformation problem must be related to our main results, 
at least up to homotopy. 
A harder problem will be to study their global counterparts. More concretely, it is conjectured in 
\cite{Severa, Sev:let} that $\bN$-graded $Q$-manifolds integrate to Lie $n$-groupoids. In \cite{get:shi}, 
shifted symplectic structures on Lie $n$-groupoids were introduced. With the help of algebraic geometry
\cite{cal:cot, ptvv}, one can also introduce shifted Lagrangian higher subgroupoids and ask for their 
deformation theory. A first step in this direction was done in \cite{car:sym} where deformations of 
symplectic groupoids are studied.

Last but not least, we comment on the quantization problem associated to a Lagrangian $NQ$-submanifold. For $n=1$, 
i.e. a coisotropic submanifold inside a Poisson manifold $j:C\hookrightarrow (M, \pi)$, the $L_\infty$-algebra 
$(C_{N^*[1]C}, \{l^i\})$ we obtain was introduced in \cite{CattaneoFelder} for the sake of quantizing the reduced 
Poisson manifold $(\underline{C}, \pi_{red})$ via the identification $\Cinfty(\underline{C}):=H^0_{l^1}(C_{N^*[1]C})$. 
In fact, they use the observation that the $L_\infty$-structure is actually a $P_\infty$-algebra and use a $\bZ$-graded 
version of the Kontsevich formality to obtain an $A_\infty$-structure. For $n>1$, 
similar ideas in the context of derived algebraic geometry are explored in $\S$4.4 of \cite{pri:lag} using the higher 
versions of formalities proved in \cite{cal:form}. In our case, we actually obtain a shifted version of a $P_\infty$-algebra,
i.e. a shifted $L_\infty$-algebra structure compatible with the multiplication of functions, so it would be 
interesting to see what these actually induce concretely on the cohomology of the $L_\infty$-algebras.

\subsection*{Acknowledgements.}
The authors are grateful to Jonas Lenthe, Luca Vitagliano and Marco Zambon for fruitful discussions,
which helped clarifying many points of this article. M.C. thanks the hospitality of Albert-Ludwigs-Universit\"{a}t 
Freiburg during two short visits as well as the members of the Higher Structures seminar at Georg-August-Universit\"{a}t G\"{o}ttingen.
J.S. is grateful to Georg-August-Universit\"{a}t in G\"{o}ttingen for the support making a short term visit at the department of 
mathematics possible. 

\section{Preliminaries on graded manifolds}\label{sec:Pre}
This section contains the standard material on $\mathbb{N}$-graded manifolds which is
needed later in our study of Lagrangian $NQ$-submanifolds of symplectic $NQ$-manifolds. 
Our presentation closely follows standard references in the world of  $\mathbb{Z}_2$-graded manifolds
(supermanifolds)  or $\mathbb{Z}$-graded manifolds such as \cite{lei:sos, fio:sup, cat:int, del:not, kos:gra, RajThes}. The later are related notions but slightly different.

\subsection{Graded manifolds.}\label{sec:GrMan}

Given a non-negative integer $n\in\bN$, let 
$\mathbf{V}=\oplus_{i=1}^n V_{i}$ be a graded vector space with $\dim V_i=m_i$.
An \emph{$\mathbb{N}$-manifold of degree $n$} (or simply an \emph{$n$-manifold})  
is a ringed space $\cM=(M, C_\cM)$, where $M$ is a smooth manifold and $C_\cM$ is a 
sheaf of graded commutative algebras such that, for all $p\in M$, there exists a 
neighbourhood $U$ and an isomorphism
    \begin{equation}\label{locality}
        C_\cM\at{U}\cong \Cinfty(U)\otimes \sym \mathbf{V}
    \end{equation}
of sheaves of graded commutative algebras, where  $\sym \mathbf{V}$ denotes 
the unital graded symmetric algebra generated by $\mathbf{V}$\footnote{When $n=0$, the vector space $\mathbf{V}=0$ so $\sym \mathbf{V}=\mathbb{R}$.}. We say that $\cM$ has 
\emph{dimension $m_0|m_1|\cdots|m_n$}, where $m_0=\dim M$ and \emph{total dimension} 
$\Tdim(\cM)=\sum_{i=0}^nm_i$ .
The manifold $M$ is known as the \emph{body} of $\cM$ and is frequently denoted by $|\cM|$. As we work in the smooth setting, we denote the global sections of the sheaf $C_\cM$ also by  $C_\cM=C_\cM(M)$, we hope this abuse of notation  will not cause confusion. 
A global section of $C^i_\cM$ is called a \emph{homogeneous function of degree $i$}; 
we write $|f|=i$ for the degree of $f$. Given $p\in M$, we denote by $C_{\cM}\at{p}$ the
 \emph{stalk} of $C_\cM$ at $p$, and given $f\in C_{\cM}\at{U}$ with $p\in U$, we denote by 
 $\mathbf{f}$ its class in $C_{\cM}\at{p}$.

 Let $\cM=(M,C_\cM)$ be an $n$-manifold of dimension $m_0|\cdots|m_n$. A chart $U\subseteq M$ 
 for which  \eqref{locality} holds is called a \emph{chart} of $\cM$.  
 We say that $$\{ x^{\alpha_i}_i\ |\ \ 0\leq i\leq n, \ 1\leq \alpha_i\leq m_i \}$$   
 are \emph{local coordinates} of $\cM\at{U}=(U,C_\cM|_U)$  if 
 $\{x^{\alpha_0}_0\}_{\alpha_0=1}^{m_0}$ are local coordinates on $U$ and 
 $\{x^{\alpha_i}_i\}_{\alpha_i=1}^{m_i}$ form a basis of $V_i$ for all $i =1,\ldots, n$.  
 Hence $C_{\cM}\at{U}$ is generated 
over $\Cinfty(U)$ by $\{x^{\alpha_i}_i\}$ with $|x^{\alpha_i}_i|=i.$

A morphism of $n$-manifolds $\Psi\colon \cM\to \cN$ is a degree preserving morphism of ringed 
spaces, given by a pair $\Psi=(\psi, \psi^\sharp)$, where $\psi\colon M \to N$ is a smooth map 
and $\psi^\sharp\colon C_\cN\to\psi_* C_\cM$ is a degree preserving morphism of sheaves of algebras 
over $N$. We will additionally use the pull-back of functions, for $f\in C_\cN$ then $\Psi^*f\in C_\cM$ 
induced by a pair $\Psi=(\psi, \psi^\sharp)$, as we will introduce in Example \ref{ex:pullback}

An $n$-manifold $\cM=(M, C_\cM)$ gives rise to a tower of graded manifolds
\begin{equation}\label{tower}
M=\cM_0\leftarrow\cM_1\cdots\leftarrow\cM_{n-1}\leftarrow\cM_{n}=\cM,
\end{equation}
where $\cM_i=(M, C_{\cM_i})$ is an $i$-manifold with $C_{\cM_i}$ the subsheaf 
of algebras of $C_\cM$ locally generated by functions of degree $\leq i$, see 
\cite{Roytenberg}.

\begin{remark}[Geometrization of graded manifolds]\label{rem:geo}
In \cite{bur:frob}, it is proven that $n$-manifolds $\cM$ are equivalent to 
\emph{admissible $n$-coalgebra bundles} $({\bf E},\mu)$. Those are graded vector bundles 
${\bf E}=\oplus_{i=-n}^{-1}E_i\to M$ with a particular coalgebra structure 
$\mu\colon {\bf E}\to{\bf E}\otimes {\bf E}$. We will occasionally use this result, 
for more details see \cite{bur:frob}.
\end{remark}

\subsection{Vector fields and tangent vectors.}\label{sec:vec}

Let $\cM=(M,C_\cM)$ be an $n$-manifold. A \emph{vector field of degree $k$} on $\cM$ 
is a degree $k$ derivation $X$ of the graded algebra $C_\cM$ of global sections, i.e.,  
an  $\mathbb R$-linear map $X\colon C_\cM\to C_\cM$ with the property that, for all 
$f, g \in C_\cM$ with $f$ homogeneous, $|X(f)| =|f|+k$ and
    \begin{equation}\label{der}
 X(fg)=X(f)g+(-1)^{|f|k}fX(g).
    \end{equation}
Vector fields give rise to a sheaf of $C_\cM$-modules over $M$. The sheaf of degree 
$k$ vector fields is denoted by $\fX^{1,k}(\cM)$. The graded commutator of vector 
fields, defined for homogeneous vector fields $X, Y$ by
\begin{equation*}
    [X,Y]=XY-(-1)^{|X||Y|}YX,
\end{equation*}
turns $\fX^{1,\bullet}(\cM)$ into a sheaf of graded Lie algebras. If $U\subseteq M$ is a chart of $\cM$ 
and we have local coordinates $\{x_i^{\alpha_i}\ | \ 0\leq i\leq n, 1\leq \alpha_i\leq m_i\}$ then the vector 
fields on $\cM\at{U}$ are generated as $C_\cM\at{U}$-module by
\begin{equation}\label{eq:basisvf}
    \fX^{1,\bullet}(\cM\at{U})=\langle \frac{\partial}{\partial x^{\alpha_i}_i}\ | \ 0\leq i\leq n, 1\leq \alpha_i\leq m_i\rangle,
\end{equation}
where $\frac{\partial}{\partial x^{\alpha_i}_i}$ are defined as derivations acting on coordinates by 
$$\frac{\partial}{\partial x^{\alpha_i}_i}(x^{\beta_j}_j)=\delta_{ij}\delta_{\alpha_i\beta_j}.$$
Note that this already implies that for another coordinate system $\{\hat{x}_j^{\alpha_j}\ | \ 0\leq j\leq n, 1\leq \alpha_j\leq m_i\}$ we have 
\begin{align}\label{Eq: TransVF}
\frac{\partial}{\partial x^{\alpha_i}_i}=\sum_{j=0}^n\sum_{\alpha_j=1}^{m_j} \frac{\partial \hat{x}^{\alpha_j}_j}{\partial x^{\alpha_i}_i}
\frac{\partial }{\partial \hat{x}^{\alpha_j}_j}.
\end{align}

\begin{example}[The Euler vector field]\label{Ex: Eulervf}
Given a graded manifold $\cM=(M,C_\cM)$ there is always a canonical degree $0$ vector field $\cE_u$, 
called the \emph{Euler vector field}, defined by $\cE_u(f) = |f| f$. Let $U\subseteq M$ be a chart of $\cM$ 
with local coordinates $\{x_i^{\alpha_i}\},$ then on $\cM\at{U}$ the Euler vector field has the form
\begin{equation*}
 {\cE_u}\at{U}=\sum_{i=0}^n\sum_{\alpha_i=1}^{m_i} |x^{\alpha_i}_i|\ x^{\alpha_i}_i\frac{\partial}{\partial x^{\alpha_i}_i}.
\end{equation*}
\end{example}

Note that vector fields with a odd degree do not necessarily commute with themselves, but if they do they are of certain interest.   

\begin{definition}[$NQ$-manifolds]\label{ex:Qman}
A degree $1$ vector field on an  $\mathbb{N}$-graded manifold $\cM$, $Q\in\fX^{1,1}(\cM),$ such that $$[Q,Q]=2Q^2=0,$$ is 
called a \emph{homological vector field}. The pair $(\cM , Q)$ is called an \emph{$NQ$-manifold}.
\end{definition}

Note that in the case of Definition \ref{ex:Qman}, the ring of functions 
$(C_\cM, Q)$ is a differential complex, justifying 
the terminology homological vector field.  A $1$-manifold with a 
homological vector field is the same as a Lie algebroid \cite{Vai}; in general, 
$NQ$-manifolds codify Lie $n$-algebroids, see \cite{bon:on}.

\begin{remark}[Multivector fields and the Schouten bracket]\label{rmk:mult}
   For an integer $k\in\mathbb{Z}$ we define the symmetric algebra 
    $$\Big(\bigoplus_{i=0}^\infty \sym^i \fX^{1,\bullet}(\cM)[k], \wedge\Big)
    \quad\text{with}\quad X\wedge Y=(-1)^{1\cdot 1+(|X|-k)(|Y|-k)}Y\wedge X\quad X,Y\in \fX^{1,\bullet}(\cM).$$
    For short, we  write $\sym^i \fX^{1,\bullet}(\cM)[k]=\vf^{i,\bullet}(\cM)$ (omitting the $k$).
    If we extend the Lie bracket using the Leibniz rule, this algebra carries a Poisson bracket of bi-degree 
    $(-1, k)$ that we denote by $[\cdot,\cdot]$ and which is known as the \emph{Schouten bracket}, which is in our case shifted by bidegree $(0,k)$. For $k=0,$ we just refer to $\Big(\bigoplus_i \vf^{i,\bullet}(\cM), \wedge\Big)$
    as the \emph{multivector fields} on $\cM$.
\end{remark}

\begin{remark}[Flows of degree $0$ vector fields]\label{Rem: FlowVec0} Let $\cM=(M, C_\cM)$ be 
an $n$-manifold and $I:=(a,b)$ an interval. A \emph{time-dependent degree $0$ 
vector field}, $X_t\in\mathfrak{X}^{1,0}(\cM)$ for $t\in I,$ is a smooth curve of 
vector fields of degree $0$. Throughout this work we will use flows of such vector fields. 
We define the \emph{flow of} $X_t$ as a smooth curve of  diffeomorphisms of $\cM$,  
$\Phi_t\colon \cM\to \cM,$ with $\Phi_0=\id$ and satisfying the flow equation 
	\begin{align*}
	\frac{\D}{\D t}\Phi_t^* f= \Phi_t^*X_t (f)\qquad \forall f\in C_\cM.
	\end{align*}	 
An easy way of showing the existence of flows is the following. According to Remark \ref{rem:geo},  
the  $n$-manifold $\cM$ is equivalent to an admissible $n$-coalgebra bundle $(\mathbf{E}, \mu)$. 
Under this equivalence, 
$X_t$ induces a degree preserving derivation of the augmented graded vector bundle $\mathbf{E}^*_{\mathbb{R}}\to M$, i.e.  
$X_t\in \Gamma \mathrm{Der}(\nE^*_{\mathbb{R}})$ see Lemma 4.8 and Theorem 4.9 in \cite{bur:frob}, and is a derivation of the coalgebra structure. Therefore, $X_t$ 
has a flow $\Phi_t$ by automorphisms of the graded vector bundle $\mathbf{E}_{\mathbb{R}}^*\to M$. One can check that 
this flow preserves indeed the coalgebra structure and therefore by \cite[Theorem 3.8]{bur:frob} it gives a diffeomorphism of $\cM$.  
In addition, we see that by construction one gets that the flow exists as long as the flow of the 
symbol of the derivation exists 
on the base manifold $M$, see also \cite[Proposition 4.14]{bur:frob}.
\end{remark}

Let $C_{\cM}\at{p}$ be the stalk of $C_\cM$ at $p\in M$. 
A \emph{homogeneous tangent vector of degree $k$ at a point $p\in M$} is 
a linear map $v\colon C_{\cM}\at{p}\to \bR$ satisfying
\begin{equation*}
 v(\mathbf{fg})=v(\mathbf{f})g^0(p)+(-1)^{k|f|}f^0(p)v(\mathbf{g})\quad \forall \mathbf{f, g}\in C_{\cM}\at{p},
\end{equation*}
where $f^0$ and $g^0$ denote the degree zero components of $f$ and $g$, respectively. 
We denote by $T_p\cM$ the space of tangent vectors at a point $p\in M$, observe that 
it is a graded vector space over $\mathbb{R}$ with only non-positive degrees.

For each open subset $U\subseteq M$, any vector field $X\in\vf^{1,\bullet}(\cM\at{U})$ defines a 
tangent vector $X_p$ at each point $p\in U$ by
\begin{equation}\label{tan-vf-gen}
X_p(\mathbf{f})=X(f)^0(p)\quad \text{for}\ \mathbf{f}\in C_\cM\at{p}.
\end{equation}
For $\mathbb{N}$-manifolds, it is clear that only non-positively graded vector fields 
can generate non-zero tangent vectors. Therefore, unlike for 
smooth manifolds, vector fields are not determined by their corresponding tangent 
vectors. Nevertheless, for all $v_p\in T_p\cM$ there exists a vector field $X\in\fX^{1,\bullet}(M)$, 
such that $X_p=v_p$. More concretely, let $U$ be a chart of $\cM$ containing $p\in M$ with coordinates 
$\{ x^{\alpha_i}_i\ |\ 0\leq i\leq n, 1\leq \alpha_i\leq m_i\}$. Then equation \eqref{eq:basisvf} gives a 
basis of vector fields on $U$ and  we conclude that a basis of  
the graded vector space $T_p\cM$ is given by 
$$T_p\cM=\langle \tfrac{\partial}{\partial x^{\alpha_i}_i}\at{p}\ |\ 0\leq i\leq n, 1\leq \alpha_i\leq m_i\rangle.$$
In particular,  $T_p\cM$ identifies with the stalk at $p\in M$ of the sheaf of vector fields $\fX^{1,\bullet}(\cM).$

\begin{remark}\label{tangentcomplex}
Let $(\cM, Q)$ be a $NQ$-manifold. It follows immediately from \eqref{tan-vf-gen} that $Q_p=0$ for all $p\in M$, since $Q(f)$ has at least degree $1$ for all functions $f\in C_\cM$.   Thus the linear part of $Q$ induces a 
differential in $T_p\cM$ and we have that the tangent is naturally a complex as explained in \cite{AKSZ}.     
\end{remark}

The \emph{differential of a map $\Psi\colon \cM\to\cN$ at $p\in M$} is the linear map
\begin{equation*}
\begin{array}{lccl}
T_p\Psi\colon &T_p\cM&\to& T_{\psi(p)}\cN, \\
&v&\longrightarrow&T_p\Psi(v)(\mathbf{g})=v(\mathbf{\psi^* g})
\quad \forall \mathbf{g}\in C_{\cN}\at{\varphi(p)}.
\end{array}
\end{equation*}

We say that a morphism $\Psi\colon \cM\to\cN$ is an 
\emph{immersion/submersion} if for all  $ p\in M$ the map $T_p\Psi$ is injective/surjective.
 Moreover, $\Psi$ is an \emph{embedding} if it is an injective immersion and $\psi$ 
 is an embedding of smooth manifolds. The following result can be proven in the usual way.

\begin{proposition}[Local normal form of immersions cf. \cite{fio:sup}] 
\label{prop: SubMfldCoo}
Let $\Psi\colon \cM\to\cN$ be an immersion and $U\subseteq M$ a chart of $\cM$ around $p$ with local 
coordinates $\{x_i^{\alpha_i}\}$. Then there exist an open $V$ with $\psi(p)\in V$ and  coordinates 
$\{\hat{x}^{\alpha_i}_i ,y^{\beta_i}_i\}$ on $\cN\at{V}$ such that  
$\Psi\at{\psi^{-1}(V)}\colon \cM\at{\psi^{-1}(V)}\to\cN\at{V}$ has the form
\begin{equation*}
\psi(x^{\alpha_0}_0)=(\hat{x}^{\alpha_0}_0,0), \quad 
\psi^* \hat{x}^{\alpha_i}_i= x_i^{\alpha_i}\quad \text{ and }
\quad \psi^* {y}_k^{\beta_k}= 0.
\end{equation*}
\end{proposition}

\subsection{Submanifolds.}
\emph{A submanifold of $\cM$} is a $\mathbb{N}$-graded manifold $\cN=(N,C_\cN)$ together with an 
embedding $(j, j^\sharp)\colon \cN\to \cM$ such that the map $j\colon N\to M$ is closed. 
In addition, we say that it is \emph{wide} if $N=M$. 

\begin{remark}
Observe that  two submanifolds $(\cN, j)$ and $(\cN', j')$ are \emph{equivalent} 
if there exists a diffeomorphism $\Psi\colon \cN\to \cN'$ such that $j=j'\circ\Psi$. We 
will make no distinction between equivalent submanifolds, and we will keep the term 
submanifold to refer to an equivalence class.
\end{remark}

Let $\cM=(M,C_\cM)$ be a $\mathbb{N}$-graded manifold. We say that $\cI\subseteq C_\cM$ is a subsheaf 
of \emph{homogeneous ideals} if for all $U$ open of $M, \ \cI\at{U}$ is an ideal of 
$C_\cM\at{U}$, i.e. $$C_\cM\at{U}\cdot \cI\at{U}\subseteq \cI\at{U},$$ and for all 
$ f\in \cI\at{U}$ its homogeneous components also belong to $\cI\at{U}$.

Given a subsheaf of homogeneous ideals $\cI\subseteq C_\cM$ we can define a subset of $M$ by
\begin{equation*}
Z(\cI)=\bigcap_{f\in \cI\cap C_\cM^0} f^{-1}(\{0\}).
\end{equation*}
We say that $\cJ$, a subsheaf of homogeneous ideals, is \emph{regular} if for each $p\in Z(\cI),$ it
exists $U\subseteq M$ a chart of $\cM$ around $p$
 and local coordinates 
$\{x_i^{\alpha_i}, y_i^{\beta_i}\}$ of $\cM\at{U}$ for which $\cI\at{U}=\langle y_i^{\beta_i}\rangle$. 
In this case, $Z(\cI)$ becomes a closed embedded submanifold of $M$. A direct consequence 
of the local normal form for immersions is the following.

\begin{proposition}\label{sub-ideal}
There is a one to one correspondence between submanifolds of $\cM$ and subsheaves of 
regular homogeneous ideals of $C_\cM$.
\end{proposition}

\begin{proof}
    We give an idea of the proof, for more details in the context of supermanifolds see \cite[Proposition 5.3.8]{fio:sup}. Given a submanifold $j:\cN\to \cM$ we get the pullback map $j^*:C_\cM\to C_\cN$. Define the corresponding subsheaf of homogeneous regular ideals as $$\cI_\cN=\ker(j^*).$$ Clearly $\cI_\cN$ is a subsheaf of homogeneous ideals of $C_\cM$ and the fact that it is regular follows Proposition \ref{prop: SubMfldCoo}. Conversely, given $\cI_\cN$ a subsheaf of homogeneous ideals of $C_\cM$ we define the $n$-manifold $$\cN=(N=Z(\cI_\cN), j^*(C_\cM/\cI_{\cN}))$$ where $j:N\to M$ is the natural embedding. The fact that $\cI_\cN$ is regular implies that $N$ is a closed embedded submanifold and $\cN$ is an $n$-manifold. Finally, the natural map $j^*:C_\cM\to j^*(C_\cM/\cI_{\cN})$ makes $j:\cN\to \cM$ into a submanifold of $\cM$.  
\end{proof}

Let $(\cN,j)$ be a  submanifold of the $NQ$-manifold $(\cM, Q)$ with associated ideal 
$\cI_\cN$. We say that $(\cN,j)$ is a \emph{$NQ$-submanifold} if
\begin{equation}
Q(\cI_\cN)\subseteq \cI_\cN.
\end{equation}

\subsection{Vector bundles.}\label{vector-bundles}

Following \cite[$\S 2.2.2$]{RajThes} we say that a \emph{vector bundle $\pi\colon \cE\to\cM$} is given by two $n$-manifolds $\cE=(E, C_\cE)$  
and $\cM=(M,C_\cM)$ and a submersion $\pi=(p,\pi^\sharp)\colon \cE\to\cM$ such that 
there exists a cover $\{U_\lambda\}_\lambda$ of $M$ with
\begin{equation*}
\cE\at{p^{-1}(U_\lambda)}\cong \cM\at{U_\lambda}\times\bR^{k_0|\cdots|k_n},
\end{equation*}
and such that the transition functions between different opens are degree preserving and fiber-linear.
 The numbers $k_0|\cdots|k_n$ are called the \emph{rank of $\pi\colon \cE\to\cM$}. 
 Therefore, a vector bundle can be defined by specifying trivializations over an open cover
 and gluing them using transition functions, as is common in differential geometry.

All natural operations such as sums, duals, pull-backs or tensor products are also defined for 
graded vector bundles\footnote{Many of these operations leave the world of 
$\mathbb{N}$-graded manifolds but are well defined as $\mathbb{Z}$-graded manifolds, 
for a more detailed treatment see \cite{RajThes}.}. In particular, using the tensor 
product we can define \emph{shifting} by $j\in\bZ$, formally 
$\cE[j]=\cE\otimes\bR_\cM[j]$ where $\bR_\cM[j]$ denotes the trivial bundle over 
$\cM$ with just one linear coordinate in degree $j$, so the new vector bundle has
 fibre coordinates shifted by $j$. 

 A \emph{homogeneous section of degree $j$} is a morphism of $\bZ$-graded manifolds 
 $s\colon \cM\to\cE[j]$ such that $\pi\circ s=\id_\cM$. Sections of degree $j$ are denoted by 
$\Gamma_{j}\cE$ and additionally $\Gamma\cE:=\oplus_{j\in\mathbb{Z}}\Gamma_j\cE$ is 
called the sections of $\cE$. For us, the sections of a vector bundle are always a left $C_\cM$-module. So if $\cE\at{U}$ is trivial, for 
some open $U\subseteq |\cM|$, then
	\begin{align*}
	\Gamma\cE\at{U}=\langle b^i_{\alpha_i} \quad \text{with}\quad |b^i_{\alpha_i}|
	=
	-i \ \text{and} \ 0\leq i\leq n, 1\leq \alpha_i\leq k_i\rangle
	\end{align*}	   
 as a $C_\cM\at{U}$-module. Moreover, as it happens in differential geometry 
 we get the following correspondence.

\begin{proposition}\label{prop:sec-VB}
 Let $\cM=(M, C_\cM)$ be an $n$-manifold. The functor 
 $$\Gamma\colon  \{ \text{Vector bundles over } \cM\}\to \{ \text{Locally finitely generated sheaves of } C_\cM\text{-modules}
 \footnote{If we insist on the total space $\cE$ being an $\bN$-manifold then 
 the local generators of the sheaf must possess non-positive degrees.}\}$$
 is an equivalence of categories.
\end{proposition}
Hence, this result allows us to describe vector bundles by their sheaves of sections. We will make use of this 
correspondence frequently.  

\begin{example}[Pull-backs]\label{ex:pullback}
An standard reference for the technical details on sheaves is $\S$ 2.3 in \cite{scha:sheaf}. Recall that for a sheaf $\cO$ on a manifold 
$M$ and a map $\psi\colon N\to M$, the \emph{inverse image sheaf} is the sheaf associated to 
the presheaf $$(\psi^{-1}\cO)(U):=\lim_{\underset{V\supseteq \psi(U)}{\longrightarrow}}\cO(V).$$ 
Hence for a vector bundle $\cE\to \cM$ and a map $\Psi\colon \cN\to\cM$, the sheaf of sections of 
the pull-back bundle  $\Psi^*\cE$ is given by
$$\Gamma\Psi^*\cE=C_\cN\otimes_{\psi^{-1}C_\cM}\psi^{-1}\Gamma\cE.$$  
This implies that for a given section $s\in\Gamma_j \cE$, we get a section on the pull-back by $\Psi^*s=s\circ \Psi.$
In particular, for the trivial bundle $\bR_\cM$,  we have that   $\Gamma \bR_{\cM}=\langle 1\rangle$ as $C_\cM$-modules, 
so we can identify $C_\cM$ with sections of the trivial bundle and $\Psi^*\bR_\cM=\bR_\cN$. Thus, we define the 
\emph{pull-back of functions} by $$\Psi^*f\in \Gamma \bR_\cN=C_\cN,\quad \text{for any}\quad f\in C_\cM.$$  
Note that if the vector bundle $\cE$ is given in terms of (local) transition functions, the transition functions of $\Psi^*\cE$ are given by the pull-backs of the transition functions of $\cE´$. 
\end{example}

\begin{remark}[Change of trivializations] \label{rm:chacoo}
Let $\pi\colon \cE\to \cM$ be a vector bundle over $\cM=(M, C_\cM)$ and consider $U,V\subseteq M$ two
 trivializing opens with $U\cap V\neq \emptyset$. Let us denote by $$\{ e_i^{\beta_i}\ | \ 0\leq i\leq n, \ 1\leq \beta_i\leq k_i\}
 \quad\text{and}\quad \{ \widehat{e}_i^{\beta_i}\ | \ 0\leq i\leq n, \ 1\leq \beta_i\leq k_i\}$$ the linear coordinates 
 on $\cE\at{p^{-1}(U)}$ and $\cE\at{p^{-1}(V)}$, respectively. Using these coordinates, we define a basis $\{b^i_{\beta_i}\}$
 for $\Gamma\cE\at{U}$ as $C_\cM\at{U}$-module by 
$$ b^i_{\beta_i}\colon \cM\at{U}\to \cE[-i]\at{p^{-1}(U)},\quad  (b^i_{\beta_i})^*(f)=f\ \text{ for } f\in C_\cM|_U\ \text{and} 
\quad (b^i_{\beta_i})^*(e_j^{\beta_j})=\delta_{ij}\delta_{\beta_i\beta_j}.$$
 Therefore, every section $s\in\Gamma \cE\at{U}$ can be written as   
 \begin{align*}
	s=\sum_{i=0}^{n}\sum_{\alpha_i=1}^{k_i} s_i^{\alpha_i} b_{\alpha_i}^i\quad \text{ where }\ s_i^{\alpha_i}=s^* e^{\alpha_i}_i\in C_{\cM}^{|s|+i}.
	\end{align*}		 
 The statement for $\Gamma\cE\at{V}$ is analogous. The transition rule implies that on $\cE\at{p^{-1}(U\cap V)}$ there 
 are functions $(\psi_j^{\beta_j})^i_{\beta_i}\in C^{j-i}_\cM\at{U}$ such that
\begin{align*}
	\widehat{e}_{j}^{\beta_j}=\sum_{i=0}^j\sum_{\beta_i=1}^{k_i} e_i^{\beta_i}(\psi_j^{\beta_j})^i_{\beta_i}. 
\end{align*}
We can organize the functions $(\psi_j^{\beta_j})^i_{\beta_i}$ into an invertible lower triangular block matrix $\Psi$. If 
we denote its inverse by $\Phi$, i.e. $\Phi\Psi=\id$, then the basis of sections on $\cE\at{U\cap V}$ transforms 
using $\Phi^T$, i.e. 
\begin{align*}
	\widehat{b}^j_{\beta_j} =\sum_{i=j}^n\sum_{\beta_i=1}^{k_i}(\phi^j_{\beta_j})^{\beta_i}_i b^{i}_{\beta_i}.
	\end{align*}		 
\end{remark}

\begin{example}[The tangent bundle]\label{tanget-gm}
Let $\cM=(M,C_\cM)$ be an $n$-manifold of dimension $m_0|\cdots|m_n$. Here we show that the tangent bundle of $\cM$, 
$\pi=(p,\pi^\sharp)\colon T\cM=(TM,C_{T\cM})\to\cM$,  is a graded vector bundle of rank $m_0|\cdots|m_n$. 
Let $U\subset M$ be a chart of $\cM$ with local coordinates $\{x_i^{\alpha_i}\}$ where
$0\leq i\leq n, \ 1\leq \alpha_i\leq m_i$. Then $T\cM\at{p^-1(U)}$ has linear fibre coordinates given by 
$$v^{\alpha_i}_i=\D x_i^{\alpha_i}\quad \text{with}\quad |v^{\alpha_i}_i|=|x^{\alpha_i}_i|,\quad \text{for}
\quad 0\leq i\leq n, \ 1\leq \alpha_i\leq m_i.$$
Since we have not introduced  $\D x_i^{\alpha_i}$ yet,  we define them by using transition functions of the coordinate vector fields, 
see Equation \eqref{Eq: TransVF}. 
If $V\subset M$ is another chart of $\cM$ with coordinates $\{\widehat{x}_i^{\alpha_i}\}$ such that on $U\cap V\neq \emptyset$
 we have the change of coordinates 
$$\widehat{v}_j^{\alpha_j}
=\sum_{i=0}^j\sum_{\alpha_i=1}^{m_i}v_i^{\alpha_i}\frac{\partial \hat{x}^{\alpha_j}_j}{\partial x_i^{\alpha^i}}.$$
Finally, notice that we can identify sections of the tangent bundle with vector fields, i.e. $\Gamma_jT\cM=\fX^{1,j}(\cM)$, 
as follows. On the chart $U\subseteq M$, let $b^i_{\alpha_i}\colon \cM\at{U}\to T[-i]\cM\at{p^{-1}(U)}$ be the map defined by 
$$(b^i_{\alpha_i})^*(f)=f \ \text{for}\ f\in C_\cM\at{U}\ \quad \text{and}
\quad (b^i_{\alpha_i})^*(v^{\alpha_j}_j)=\delta_{ij}\delta_{\alpha_i\alpha_j}.$$
Then $b^i_{\alpha_i}$ identifies with the coordinate vector field $\frac{\partial}{\partial x_i^{\alpha_i}}$ on  $C_\cM\at{U}$
via
$$(b^i_{\alpha_i})^*(\sum_{j=0}^{n}\sum_{\alpha_j=1}^{m_j}v_j^{\alpha_j}\frac{\partial f}{\partial x_j^{\alpha_j}})=\frac{\partial f}
{\partial x_i^{\alpha_i}}=\frac{\partial}{\partial x_i^{\alpha_i}}(f)\quad \forall f\in C_{\cM}\at{U}.$$
\end{example}

\begin{example}[Normal bundle to a submanifold]\label{ex:nor}
As usual, we introduce the \emph{normal bundle to a submanifold} $j\colon \cN\to \cM$ as the 
vector bundle over $\cN$ given by $\nu(\cM,\cN)=j^*T\cM/T\cN$.  For later use we give 
a more accurate description of its sheaf of sections.
Let $\cI_\cN$ be the vanishing ideal associated to the submanifold $j\colon \cN\to \cM$ and 
consider the sheaf of vector fields on $\cM$ tangent to $\cN$, 
\begin{equation}\label{eq:vec-N}
    \mathfrak{X}^{1,\bullet}_\cN(\cM)=\{ X\in \mathfrak{X}^{1,\bullet}(\cM) \ | \ X(\cI_\cN)\subseteq \cI_\cN\}.
\end{equation}
As an immediate consequence of Proposition \ref{prop: SubMfldCoo} and Example \ref{ex:pullback}, we get that the 
sheaf of sections of the normal bundle is
$$\Gamma\nu(\cM,\cN)=C_\cN\tensor_{j^{-1}C_\cM} j^{-1}\mathfrak{X}^{1,\bullet}(\cM) \mod 
\big(C_\cN\tensor_{j^{-1}C_\cM} j^{-1}\mathfrak{X}^{1,\bullet}_\cN(\cM)\big).$$
Since $j\colon \cN\to \cM$ is a closed embedding, sometimes we will consider the sheaf over $M$ given by 
$$j_*C_\cN\tensor_{C_\cM} \mathfrak{X}^{1,\bullet}(\cM) \mod \big(j_*C_\cN\tensor_{C_\cM} \mathfrak{X}^{1,\bullet}_\cN(\cM)\big).$$
\end{example}

\begin{lemma}\label{Lem: NormVec}
Let $\cE\to \cM$ a graded vector bundle, and $0_\cM\colon \cM\to \cE$ the submanifold given by the zero section. 
Then, there is a canonical isomorphism of vector bundles $\nu(\cE,\cM)\cong \cE$.  
\end{lemma} 

\begin{proof}
Let $U\subseteq|\cM|$ be an open such that $\cE\at{U}$ is trivial with linear coordinates $\{e_i^{\alpha_i}\}.$ 
 As explained in Remark \ref{rm:chacoo}, they induce the basis of sections $\{b_{\alpha_i}^i\}$. Every section $s\in \Gamma_k\cE\at{U}$
of degree $k$ can locally be expanded in the basis $\{b^i_{\alpha_i}\}$:
	\begin{align*}
	s=\sum_{i=0}^{n} \sum_{\alpha_i=1}^{k_i} s^*(e^{\alpha_i}_i)b^i_{\alpha_i}=\sum_{i=0}^{n} \sum_{\alpha_i=1}^{k_i} s^{\alpha_i}_i b_{\alpha_i}^i
	\end{align*}	    
with $s^{\alpha_i}_i\in C^{k+i}_{\cM |_U}$.  The vertical lift of the section $s$ is the degree $k$ vector field in $\cE\at{U}$ given by
\begin{align*}
s^v:=\sum_{i=0}^{n} \sum_{\alpha_i=1}^{k_i} s^{\alpha_i}_i \frac{\partial}{\partial e_i^{\alpha_i}}\in \fX^{1,k}(\cE\at{U}).
\end{align*}
One can show that the above construction yields a global vector field, since the local basis $b_{\alpha_i}^i$ 
transforms exactly in the same way as the corresponding partial derivatives $\frac{\partial}{\partial e_i^{\alpha_i}}$. 
Therefore we get a well defined map of sheaves of $C_\cM$-modules 
$$\Gamma\cE\to C_\cM\otimes_{0_\cM^{-1}C_\cE}0_\cM^{-1}\fX^{1,\bullet}(\cE), \quad s\to 1\otimes s^v\at{\cM}.$$
Observe that since vertical vector fields are not tangent to the zero section and $\mathfrak{X}^{1,k}_\cM(\cE)$ is 
generated by vector fields of coordinates in $\cM$, this  is the only part that is not in the image. Therefore, we get that 
the quotient map
$$\Gamma \cE\to \Gamma\nu(\cE,\cM),\quad  s\to [1\otimes s^v\at{\cM}]$$
is an isomorphism of sheaves of $C_\cM$-modules.
\end{proof}

Given coordinates adapted to a submanifold $j\colon\cN\hookrightarrow\cM$ as in Proposition 
\ref{prop: SubMfldCoo}, i.e. local coordinates $\{x_i^{\alpha_i},y_i^{\beta_i}\}$ of $\cM$ such that 
$\cJ_\cN=\langle y_i^{\beta_i}\rangle$, one can check that 
	\begin{align*}
	\Big[\tfrac{\partial}{\partial y_i^{\beta_i}}\Big] := 
	\tfrac{\partial}{\partial y_i^{\beta_i}} \mod \big(j_*C_\cN\tensor_{C_\cM} 
	\mathfrak{X}^{1,\bullet}_\cN(\cM)\big)
	\end{align*}
are a local basis for $\Gamma\nu(\cM,\cN)$. If we pick a different set of coordinates 
$\{x_i^{\alpha_i},\tilde{y}_i^{\beta_i}\}$, we get the transition functions 
	\begin{align*}
	\big[\tfrac{\partial}{\partial y_i^{\beta_i}}\big] = 
	j^*(\tfrac{\partial\tilde{y}_k^{\beta_k}}{\partial y_i^{\beta_i}}) 
	\Big[\tfrac{\partial}{\partial \tilde{y}_k^{\beta_k}}\Big].
	\end{align*}

\subsection{Differential forms and Cartan calculus.}

 A \emph{differential $k$-form of degree $j$} on the $n$-manifold $\cM=(M,C_\cM)$, denoted by $\omega\in\Omega^{k,j}(\cM)$, is a map
\begin{equation*}
 \omega\colon \fX^{1,i_1}(\cM)\times\cdots\times\fX^{1,i_k}(\cM)\to C_\cM^{i_1+\cdots+i_k+j}
\end{equation*}
satisfying
 $$\omega(\cdots, X, Y, \cdots)=-(-1)^{|X||Y|}\omega(\cdots,Y,X,\cdots) \quad \text{and}\quad 
 \omega(fX_1,\cdots, X_k)=(-1)^{|f||\omega|}f\omega(X_1,\cdots, X_k).$$

Suppose that 
$\omega\in\Omega^{i,j}(\cM)$ and $\eta\in\Omega^{k,l}(\cM)$. Then, we define a new form 
$\omega\wedge\eta\in\Omega^{i+k,j+l}(\cM)$ by the formula
\begin{equation*}
\omega\wedge\eta(X_1,\cdots, X_{i+k})=\sum_{\sigma\in Sh(i,k)} Ksgn(\sigma)
\omega(X_{\sigma(1)},\cdots, X_{\sigma(i)})\eta(X_{\sigma(i+1)},\cdots,X_{\sigma(i+k)})
\end{equation*}
where $X_1,\cdots, X_{i+k}\in\fX^{1,\bullet}(\cM)$ and $Ksgn(\sigma)$ denotes the 
signature of the permutation multiplied by the Koszul sign. 
$\omega \wedge \eta$ is referred to as the wedge product of $\omega$ and $\eta$. 
Observe that, with our 
sign conventions, we have that $$\omega\wedge\eta=(-1)^{ik+jl}\eta\wedge\omega.$$
The rule for the wedge product implies that $(\Omega^{\bullet,\bullet}(\cM),\wedge)$ 
is generated by $\Omega^{1,\bullet}(\cM)$ as a sheaf of algebras over $C_\cM$.

Analogously to the case of vector fields and tangent vectors, we have the following relation between $1$-forms 
and covectors. The \emph{cotangent space at a point } $p\in M$ is the graded dual of 
the space $T_p\cM$, i.e. $T^*_p\cM=\underline{\Hom}(T_p\cM, \bR)$, which is a non-negatively graded vector 
space\footnote{By $\underline{\Hom}$ we denote the internal $\Hom$ in the category of $\bZ$-graded vector spaces.}.  
Each differential 1-form $\alpha\in \Omega^{1,n}(\cM)$ defines a covector by the formula
	\begin{align*}
	\alpha_p\colon T_p\cM\to \bR[n],\quad  v_p \mapsto \alpha(X)^0(p),
	\end{align*}	 
where $X\in\mathfrak{X}^{1,k}(\cM)$ is a vector field, such that $X_p=v_p$. It is easy to show that this 
map is well-defined. In the same way, for a differential $k$-form $\omega$ one obtains a map 
$\omega_p\colon \bigwedge^kT_p\cM\to \mathbb{R}$.
\begin{remark}
     As it happens for vector fields, differential $1$-forms are not determined by the covectors at a point, i.e.  
      there are non-vanishing differential forms $\alpha$ with  $\alpha_p=0$ for all $p\in M$. 
\end{remark}

Let $U\subseteq M$ be a chart of $\cM$ with coordinates $\{x_i^{\alpha_i}\ | 
\ 0\leq i\leq n, \ 1\leq \alpha_i\leq m_i\}.$ By \eqref{eq:basisvf},
$\{ \frac{\partial}{\partial x_i^{\alpha_i}}\}$ is a basis for $\fX^{1,\bullet}(\cM\at{U})$, thus 
$$\Omega^{1,\bullet}(\cM\at{U})=\langle \D x_i^{\alpha_i}\ | \ 0\leq i\leq n, \ 1\leq \alpha_i\leq m_i\rangle$$ 
defines a basis as $C_{\cM|_U}$-modules as follows  
\begin{equation}\label{eq:def-1form}
\D x_i^{\alpha_i}(\frac{\partial}{\partial x_j^{\alpha_j}})
=(-1)^{ij}
\frac{\partial}{\partial x_j^{\alpha_j}}(x_i^{\alpha_i})=(-1)^{ij}\delta_{ij}\delta_{\alpha_i\alpha_j}.
\end{equation}
Moreover, this identity implies that a basis of $T_p^*\cM$ is 
given by $\{\D x^{\alpha_i}_i\at{p}\}$. 
For a function $f\in C_\cM$, we thus have
	\begin{align*}
	\D f\at{U}= \sum_{i=0}^n\sum_{\alpha_i=1}^{m_i}\D x_i^{\alpha_i} \frac{\partial f}{\partial x_i^{\alpha_i}}.
	\end{align*}
In particular, for a second coordinate system $\{\hat{x}_j^{\alpha_j}\ | \ 0\leq j\leq n, \ 1\leq \alpha_j\leq m_j\}$, we get
	\begin{align*}
	\D \hat{x}^{\alpha_j}_j= 
	\sum_{i=0}^n\sum_{\alpha_i=1}^{m_i}\D x_i^{\alpha_i} \frac{\partial \hat{x}^{\alpha_j}_j}{\partial x_i^{\alpha_i}},
	\end{align*}
which justifies the choices in Example \ref{tanget-gm}.

\begin{remark}\label{Rem. TransForms}
For an $n$-manifold $\cM$, in $\S$\ref{vector-bundles}, we considered the sections of vector bundles as left $C_\cM$-modules, 
but the more canonical choice for differential forms is actually a right $C_\cM$-module which can be seen in the 
transformation rule of the $\D x^{\alpha_i}_i$ above. Nevertheless, we stick to our choice and 
consider them as a left $C_\cM$-module with local bases $\{\D x^{\alpha_i}_i\}$ and the transformation rule 
	 \begin{align*}
	\D \hat{x}^{\alpha_j}_j
	= \sum_{i=0}^n\sum_{\alpha_i=1}^{m_i}(-1)^{i(j-i)}\frac{\partial \hat{x}^{\alpha_j}_j}{\partial x_i^{\alpha_i}} \D x_i^{\alpha_i}. 
	\end{align*}
Note that, up to now, the differential forms are not sections of a vector bundle yet, since $\D x^{\alpha_i}_i$ has positive degree, 
see Proposition \ref{prop:sec-VB}, We will address this point in $\S$\ref{sec: Cotangent}. 
\end{remark}

\begin{example}[Conormal sheaf to a submanifold]\label{ex:conor}
Assume now that $j\colon \cN\hookrightarrow\cM$ is a 
submanifold with vanishing ideal $\cI_\cN$. Define the \emph{conormal sheaf} by
\begin{equation}\label{eq:con-N}
    \Omega^{1,\bullet}_\cN(\cM)=
    \{ \alpha\in\Omega^{1,\bullet}(\cM) \ | \ \alpha(X)\in\cI_\cN \ \text{for all}\ X\in\fX_\cN^{1,\bullet}(\cM)\}.
\end{equation}
Differential $1$-forms are generated by exact $1$-forms, thus using Proposition \ref{prop: SubMfldCoo} together 
with Eq.\eqref{eq:def-1form} one shows that $\Omega^{1,\bullet}_\cN(\cM)=\langle df \ |\ f\in\cI_\cN\rangle$ as $C_\cM$-module.
Analogously to the normal bundle given in Example \ref{ex:nor}, one can show, using again Proposition \ref{prop: SubMfldCoo}, that
$$C_\cN\otimes_{j^{-1}C_\cM}j^{-1}\Omega^{1,\bullet}_\cN(\cM)$$
is locally finitely generated and therefore is the sheaf of sections of a $\bZ$-graded vector bundle over $\cN$.
\end{example}

We recall the basics of the Cartan calculus on graded manifolds. As in classical geometry, 
we define the contraction, Lie derivative, and de Rham differential operators for a graded 
manifold and prove the usual formulas with some extra signs. In order to know the appropriate 
sign, just recall that whenever two symbols are transposed, the sign according to their degree appear.

Let $\cM=(M,C_\cM)$ be a graded manifold and consider $X\in\fX^{1,r}(\cM)$ a vector field. The 
contraction with respect to this vector field is defined by the following rule:
\begin{equation*}
	\begin{array}{ll}
	\iota_X\colon \Omega^{k,l}(\cM)\to\Omega^{k-1,l+r}(\cM)\\
	(\iota_X\omega)(X_2,\cdots, X_{k})=(-1)^{|\omega||X|}\omega (X, X_2,\cdots, X_k).
	\end{array}
\end{equation*}
It is a derivation of the wedge product, i.e.
\begin{equation*}
\iota_X(\omega\wedge \eta)=\iota_X\omega\wedge \eta+(-1)^{(-1)k+|\omega||X|}\omega\wedge \iota_X\eta.
\end{equation*}
Therefore, $\iota_X$ is a derivation of the algebra $(\Omega^{\bullet,\bullet}
(\cM),\wedge)$ of bidegree $(-1, |X|)$. We also define the Lie derivative with respect to 
a vector field as
\begin{equation*}
 \Lie_X f=X(f), \qquad \Lie_X Y=[X,Y],
\end{equation*}
and extend it to forms by the following formula
\begin{equation}\label{lieder}
\Lie_X \iota_{X_1}\cdots \iota_{X_k}\omega=(-1)^{r\sum\limits_{a=1}^k |X_a|} \iota_{X_1}\cdots \iota_{X_k}
\Lie_X\omega+\sum_{j=1}^k (-1)^{r\sum\limits_{a=1}^{j-1}|X_a|} \iota_{X_1}\cdots \iota_{[X,X_j]}
\cdots \iota_{X_k}\omega.
\end{equation}
By definition, $\Lie_X$  is a derivation of the wedge product with bidegree $(0,|X|)$. Finally,
 the de Rham differential can be defined by the usual Cartan formula
\begin{equation*}
\begin{array}{ll}
\iota_{X_0}\cdots \iota_{X_k}\D \omega=&\sum\limits_{i=0}^k (-1)^{i+|X_i|\sum\limits_{a=0}^{i-1}|X_a|}
\Lie_{X_i}\iota_{X_0}\cdots\widehat{\iota_{X_i}}\cdots \iota_{X_k}\omega\\
&+ \sum\limits_{i<j} (-1)^{i+1+|X_i|\sum\limits_{a=i+1}^{j-1}|X_a|}\iota_{X_0}\cdots\widehat{\iota_{X_i}}
\cdots \iota_{[X_i,X_j]}\cdots \iota_{X_k}\omega.
\end{array}
\end{equation*}
This identity implies  that $\D$ is a derivation of the wedge product with bidegree $(1,0)$.
 With all these formulas we are able to prove the usual Cartan  relations.

\begin{proposition}\label{cartan-cal}
Let $\cM$ be a graded manifold. The following formulas hold:
\begin{equation*}
\begin{array}{l}
\left[\D ,\D\right]=2\D^2=0,\quad \Lie_X=\left[\iota_X, \D \right], 
\quad \Lie_{\left[X,Y\right]}=\left[\Lie_X, \Lie_Y\right],\\
\left[\Lie_X, \D \right]=0, \quad \iota_{\left[X,Y\right]}=\left[\Lie_X, \iota_Y\right], 
\quad \left[\iota_X, \iota_Y\right]=0,
\end{array}
\end{equation*}
where the bracket denotes the graded commutator with bidegree, i.e. for example  
$$[\iota_X,\iota_Y]=\iota_X\circ \iota_Y-(-1)^{(-1)(-1)+|X||Y|}\iota_Y\circ \iota_X.$$ 
Moreover, for a map $\Psi:\cM\to \cN$, the pull-back of differential forms 
$\Psi^*:\Omega^{\bullet,\bullet}(\cN)\to \Omega^{\bullet,\bullet}(\cM)$ is the unique algebra 
morphism of bidegree $(0,0)$ that extends the pull-back of functions and commutes with $\D$.
\end{proposition}

\subsection{Graded symplectic manifolds and Lagrangian submanifolds.}\label{sec:sym}
Now we are ready to introduce the main characters of this work.

\begin{definition}
A \emph{symplectic $n$-manifold} $(\cM, \omega)$ is 
an  $n$-manifold $\cM=(M, C_\cM)$ together with $\omega\in\Omega^{2,n}(\cM)$
such that $\D\omega=0$ and $\omega$ is non-degenerate, i.e. the map
\begin{align*}
\omega^\flat\colon \mathfrak{X}^{1,\bullet}(\cM)\to \Omega^{1,\bullet+n}(\cM), \quad  X\to \iota_X\omega,
\end{align*}
is an isomorphism of $C_\cM$-modules.    
\end{definition}

 Since we just consider $\bN$-manifolds, the fact that the degree of the manifold is equal to the 
 degree of the symplectic form is not a restriction at all, see \cite{Roytenberg}.  Since $T_p\cM$ 
 and $T^*_p\cM$ are the stalks at $p\in M$ of the sheaf of vector fields and differential $1$-forms 
 respectively, we can reformulate the non-degeneracy condition of $\omega$ as follows.  
 See also \cite[Theorem 39]{Vit} for a similar statement for $k$-plectic forms on $1$-manifolds.

\begin{proposition}\label{Lem: Symplatp}
Let $\cM=(M,C_\cM)$ be an $n$-manifold. A $2$-form $\omega\in \Omega^{2,n}(\cM)$ is non-degenerate
 if and only if for all $p\in M$ the map 
$$\omega_p^\flat\colon T_p\cM\to (T_p^*\cM)[n],\quad   v_p\mapsto \omega_p(v_p,\cdot ), $$
is an isomorphism. 
\end{proposition}

\begin{proof}
The statement follows from \cite[Proposition 2.2.2]{scha:sheaf} that states that a morphism of sheaves 
is an isomorphism if and only for all points the induced map on the stalks is an isomorphism.
\end{proof}

For every function $f\in C^{k}_\cM$, we denote by $X_f$ the unique degree $k-n$ vector field, such that
\begin{equation}\label{Eq: HamVF}
\omega^\flat(X_f)=\D f.
\end{equation} 
 $X_f$ is called the \emph{Hamiltonian} vector field of $f$. 
Similar to the case $n=0$, Hamiltonian vector fields allow us to define the 
associated Poisson bracket $\{\cdot,\cdot\}\colon C^i_\cM\times C^j_\cM\to C^{i+j-n}_\cM$ 
given by the formula
\begin{equation}
\{f,g\}=\iota_{X_f}\iota_{X_g}\omega.
\end{equation}

A vector field $X\in\fX^{1,k}(\cM)$ is called \emph{symplectic}  if $\Lie_X\omega=0$.

\begin{proposition}\label{bra}
Let $(\cM,\omega)$ be a degree $n$ symplectic manifold. The following equalities hold:
$$\Lie_{X_f}g=\{f,g\},\qquad \Lie_{X_f}\omega=0\quad \text{ and }\quad X_{\{f,g\}}=[X_f,X_g].$$
Moreover, $\{\cdot,\cdot\}$ is a Poisson bracket of degree $-n$.
\end{proposition}

\begin{proof}
We use Proposition \ref{cartan-cal} and formula \eqref{lieder}:
\begin{eqnarray*}
\Lie_{X_f}g&=&(\iota_{X_f}\D +\D \iota_{X_f})g=\iota_{X_f}\D g=\iota_{X_f}\iota_{X_g}\omega=\{f,g\}.\\
\Lie_{X_f}\omega&=& \D \iota_{X_f}\omega+\iota_{X_f}\D\omega=0.\\
\iota_{[X_f,X_g]}\omega&=&[\Lie_{X_f},\iota_{X_g}]\omega=\Lie_{X_f}\iota_{X_g}\omega
=\Lie_{X_f}\D g=\D \Lie_{X_f}g=\D \{f,g\}=\iota_{X_{\{f,g\}}}\omega.
\end{eqnarray*}
Next, we show that the bracket is a Poisson bracket. Since $0=[\iota_X,\iota_Y]=\iota_X\iota_Y-(-1)^{1+|X||Y|}\iota_Y \iota_X$,
 we obtain that $\{\cdot,\cdot\}$ is skew-symmetric.
The graded Leibniz identity:
\begin{equation*}
\begin{array}{rl}
\{f,gh\}=&\iota_{X_f}\iota_{X_{gh}}\omega=\iota_{X_f}\D (gh)=\iota_{X_f}\D g \ h+\iota_{X_f}g\ \D h\\
=&\{f,g\}h+(-1)^{|g||X_f|}g\ \iota_{X_f}\D h=\{f,g\}h+(-1)^{|g|(|f|-n)}g\{f,h\}.
\end{array}
\end{equation*}
The graded Jacobi identity:
\begin{equation*}
\begin{array}{rl}
\{\{f,g\},h\}+(-1)^{(|f|-n)(|g|-n)}\{g,\{f,h\}\}=&\iota_{X_{\{f,g\}}}\iota_{X_h}
\omega+(-1)^{|X_f||X_g|}\iota_{X_g}\iota_{X_{\{f,h\}}}\omega\\
=&\iota_{[X_f,X_g]}\iota_{X_h}\omega+(-1)^{|X_f||X_g|}\iota_{X_g}\iota_{[X_f,X_h]}\omega\\
=&\Lie_{X_f}\iota_{X_g}\iota_{X_h}\omega=\Lie_{X_f}\{g,h\}\\
=&\{f,\{g,h\}\}.
\end{array}
\end{equation*}
\end{proof}

Two special properties of degree $n$ manifolds are stated in the next Proposition, which is due to Roytenberg.

\begin{proposition}[see \cite{Roytenberg}]\label{sym=ham}
Let $(\cM,\omega)$ be a symplectic $n$-manifold. The following 
statements hold:
\begin{enumerate}
	\item If $n\geq 1$, then $\omega$ is exact. Moreover $\omega=\D(\frac{1}{n} 
	i_{\cE_u}\omega)$, where $\cE_u$ is the Euler vector field of $\cM$.
	\item Let $X\in\fX^{1,l}(\cM)$ be a symplectic vector field. If $l+n\neq 0$,	then $X$ is Hamiltonian with $f=\frac{1}{l+n} i_{\cE_u}i_X\omega$.
\end{enumerate}
\end{proposition}

\begin{proof}
Both formulas are straightforward:
\begin{equation*}
n\ \omega=\Lie_{\cE_u}\omega=\D i_{\cE_u}\omega\Rightarrow \omega
=\D (\frac{1}{n}i_{\cE_u}\omega).
\end{equation*}
\begin{equation*}
\D i_{\cE_u}i_X\omega=\Lie_{\cE_u}i_X\omega-i_{\cE_u}\D i_X\omega
=(l+n)i_X\omega\Rightarrow i_X\omega=\D (\frac{1}{n+l}i_{\cE_u}i_X\omega).
\end{equation*}
\end{proof}

\begin{definition}
A \emph{degree $n$ symplectic $NQ$-manifold} is a triple $(\cM,\omega,Q)$ where $(\cM,\omega)$ 
is a symplectic $n$-manifold and $(\cM, Q)$ is an $NQ$-manifold such that 
	\begin{align*}
	\Lie_Q\omega=0.
	\end{align*}
\end{definition}

\begin{corollary}
Let $(\cM,\omega,Q)$ be a degree $n$ symplectic $NQ$-manifold with $n\geq 1$. Then
	\begin{align*}
	Q=X_\theta \ \quad  \text{ with } \quad \theta \in C_\cM^{n+1} 
	\quad  \text{ and } \quad \{\theta,\theta\}=0. 
	\end{align*}	 
\end{corollary}

\begin{proof}
From Proposition \ref{sym=ham}, we know that $Q$ is Hamiltonian with $\theta\in C^{n+1}_\cM$. Moreover, 
by Proposition \ref{bra},  $0=[Q, Q]=X_{\{\theta,\theta\}},$ thus $\D\{\theta,\theta\}=0$ and 
due to the degree $\{\theta,\theta\}=0$. 
\end{proof}

\begin{remark}[Shifted symplectic algebroids]
    Let $(\cM, Q)$ be an $NQ$-manifold of degree $n$. By the Cartan calculus of Proposition 
    \ref{cartan-cal} we obtain that $(\Omega^{\bullet,\bullet}(\cM), \D, \cL_Q)$ becomes a 
    double complex with total differential $D=\cL_Q+(-1)^i\D$.  And, as explained in Remark 
    \ref{tangentcomplex}, we have that $T_p\cM$ and $T_p^*\cM$ are chain complexes. In 
    \cite{pym:shif} a \emph{shifted symplectic algebroid} is defined as a $NQ$-manifold $(\cM, Q)$ 
    endowed with $\omega\in\Omega^{2,n}(\cM)$ such that 
    $$\Lie_Q\omega=0,\quad \D\omega=D(\sum_{i=1}^n \beta_i)\quad \text{for some}\quad \beta_i\in\Omega^{2+i, n-i}(\cM)
    \quad \text{and the map}\quad \omega_p\colon T_p\cM\to (T^*_p\cM)[n]$$
    is a quasi-isomorphism for all $p\in M$. Therefore, the symplectic $NQ$-manifolds of this work 
    are a strict version of the shifted symplectic algebroids introduced in \cite{pym:shif}. We expect 
    that an \emph{up to homotopy version} of all the results obtained in the following sections can be 
    extended to shifted symplectic algebroids.  
\end{remark}

\begin{definition}\label{def:lag}
Let $(\cM, \omega)$ be a symplectic $n$-manifold  and  $j\colon \cN\hookrightarrow \cM$ a  submanifold with 
vanishing ideal $\cI_\cN$. We say that $\cN$ is
\begin{enumerate}
    \item \emph{isotropic}, if $j^*\omega=0$.
    \item \emph{coisotropic},  if 	$\{\cI_\cN,\cI_\cN\}\subseteq \cI_\cN.$
    \item \emph{Lagrangian}, if it is coisotropic with $\Tdim(\cN)=\frac{1}{2}\Tdim(\cM)$. 
\end{enumerate}
Moreover, if $(\cM, \omega, Q)$ is a symplectic $NQ$-manifold, we say that $\cN$  is a 
\emph{Lagrangian $NQ$-submanifold} if $Q$ preserves its vanishing ideal.
\end{definition}
\begin{remark}
Let $(\cM, \omega)$ be a symplectic $n$-manifold and $i\colon \cL\hookrightarrow \cM$  a Lagrangian submanifold. If $\cL$ is wide, i.e. $|\cL|=|\cM|$, by the non-degeneracy of $\omega$, one immediately sees that $\cL$ is an $(n-1)$-manifold.
\end{remark}
Given $p\in N$ we define the \emph{symplectic orthogonal at $p$} as the graded vector space given 
by $$T_p\cN^\omega=\{v\in T_p\cM \ | \ \omega_p(v,u)=0\  \forall u\in T_p\cN\}.$$
Some basic properties of these special submanifolds are stated in the following.

\begin{proposition}\label{Prop: EquLag}
    Let $(\cM, \omega)$ be a symplectic $n$-manifold and  $j\colon \cN\hookrightarrow \cM $ a 
    submanifold with vanishing ideal $\cI_\cN$. Then the following statements hold:
    \begin{enumerate}
        \item If $\cN$ is isotropic, then $T_p\cN\subseteq T_p\cN^\omega$ for all $p\in N$.
        \item If $\cN$ is coisotropic, then $X_f\in \mathfrak{X}^{1,\bullet}_\cN(\cM)$ for all 
        $f\in \cI_\cN$ and  $T_p\cN^\omega\subseteq T_p\cN$ for all $p\in N$.
        \item If $\cN$ is Lagrangian, then $\omega^\flat \colon \mathfrak{X}_\cN^{1,\bullet}(\cM)\to 
        \Omega^{1,\bullet+n}_\cN(\cM)$ is an isomorphism. Moreover $\mathfrak{X}_\cN^{1,\bullet}(\cM)
        =\langle X_f \ | \ f\in \cI_\cN\rangle$ as $C_\cM$-module and  $$X\in \mathfrak{X}_\cN^{1,\bullet}(\cM) \  
        \iff \ \omega(X,Y)\in \cI_\cL \text{ for all } Y\in\mathfrak{X}_\cN^{1,\bullet}(\cM).$$ 
    \end{enumerate}
\end{proposition}

\begin{proof}
    \begin{enumerate}
        \item Since $\cN$ is isotropic,  $0=(j^*\omega)_p=\omega_{j(p)}$ for all $p\in N$, 
        hence $T_p\cN\subseteq T_p\cN^\omega$.
        \item Since $\cN$ is coisotropic, then $\{\cI_\cN,\cI_\cN\}\subseteq \cI_\cN$, thus for 
        $f\in \cI_\cN$ we get that $$X_f(\cI_\cN)=\{f,\cI_\cN\}\subseteq \cI_\cN\quad\text{so}
        \quad X_f\in \fX^{1,\bullet}_\cN(\cM).$$ First observe that for $p\in N$ the isomorphism 
        $\omega^\flat_p$ sends $T_p\cN^{\omega}$ to $ \mathrm{Ann}(T_p \cN)[n]=\langle df_p\ | 
        \ f\in \cI_\cN\rangle$ by Example \ref{ex:conor}. This implies that $T_p\cN^{\omega}=
        \langle X_f|_p\ | \ f\in \cI_\cN\rangle$ and using the above computation we get 
        $X_f|_p\in T_p\cN$ as we want.
        \item Since $\cN$ is Lagrangian, then it is coisotropic, thus the map $\omega^\flat \colon 
        \mathfrak{X}_\cN^{1,\bullet}(\cM)\to \Omega^{1,\bullet+n}_\cN(\cM)$ is well defined by the 
        above item. By \cite[Proposition 2.2.2]{scha:sheaf}, it is enough to show that the induced 
        map on stalks is an isomorphism. If $p\notin j(N)$, then the map on the stalks is 
        $\omega_p^\flat\colon T_p\cM\to (T^*_p\cM)[n]$ which is an isomorphism by Proposition \ref{Lem: Symplatp}. 
        If $p\in j(N)$, then $\omega_p^\flat\colon T_p\cN\to \mathrm{Ann}(T_p \cN)[n].$ Since $\cN$ is coisotropic, 
        we get that this map is surjective and since both spaces have the same dimension, due to the condition 
        $\Tdim(\cN)=\frac{1}{2}\Tdim(\cM)$, we get an isomorphism. The remaining part is an easy exercise. 
    \end{enumerate}
\end{proof}

We end this subsection by commenting on two special issues of graded symplectic manifolds 
which do not appear in usual symplectic geometry. The first one is that the conditions on 
tangent spaces for isotropic and coisotropic submanifolds are not sufficient, as we illustrate 
in the subsequent example. Let $\cM=\bR^{1|4|1}$ be the $2$-manifold with coordinates
$$\{x, q^1, q^2, p^1, p^2, y\} \quad \text{of degree }\quad |x|=0, |q^1|=|q^2|=|p^1|=|p^2|=1, |y|=2$$
together with  the symplectic form given by $\omega=\D y\wedge\D x+ \D p^1\wedge\D q^1+\D p^2\wedge\D q^2.$ Consider 
the submanifold $j\colon \cN=\bR^{1|2|0}\hookrightarrow \cM$ given by the vanishing ideal
$$\cI_\cN=\langle y-p^1q^2,\ p^2,\ q^1\rangle.$$
The submanifold $\cN$ is not isotropic, because $j^*\omega=(q^2\D p^1-p^1 \D q^2)\wedge \D x\neq 0$ 
but for all $p\in N=\bR$ we get $T_p\cN\subseteq T_p\cN^\omega$. Additionally, $\cN$ is also not coisotropic, because 
$$\{ y-p^1q^2, q^1\}=\pm q^2.$$
But $T_p\cN^\omega\subseteq T_p\cN$ for all $p\in M=\bR$.

The second aspect is that not all graded symplectic manifolds admit Lagrangian submanifolds as the 
next example shows. Consider the $2$-manifold $\cM=\bR^{0|1|0}$ with coordinate $\{ e\}$ of degree 
$1$ and symplectic form $\omega= \D e\wedge \D e$, notice that this is indeed a symplectic structure. 
Since $\Tdim{\cM}=1$, $(\cM, \omega)$ cannot possess any Lagrangian submanifolds.  

\subsection{Shifted cotangent bundles.}\label{sec: Cotangent}

One 
main class of examples of graded symplectic manifolds are shifted cotangent bundles. They are crucial for the study of the 
deformation theory of Lagrangian $NQ$-submanifolds since, as we will see later, every symplectic manifold is isomorphic to a cotangent bundle around a Lagrangian submanifold.
For more applications and a systematic study of shifted cotangent 
bundles we refer the reader to \cite{cal:cot, cat:bf, CuecaCoTan}.

 For a degree $n$ manifold $\cM=(M,C_\cM)$ of dimension $m_0|\dots |m_n$,
the \emph{cotangent bundle} $\pi\colon T^*\cM\to \cM$ is defined as the $\bZ$-graded vector bundle 
dual to the tangent bundle introduced in Example \ref{tanget-gm}. Notice that if $n>0$, it is 
not an $\bN$-manifold, since there are coordinates of negative degrees. We can handle
this issue by shifting the cotangent bundle by $k\geq n$, and denote it by $T^*[k]\cM=T^*\cM\otimes\bR_\cM[k]$.
Since vector fields are sections of the tangent bundle, 
 i.e. $\fX^{1,\bullet}(\cM)=\Gamma T\cM$, we could define the \emph{$k$-shifted cotangent bundle} as 
 the vector bundle whose sections are the sheaf $\Gamma T^*[k]\cM=\Omega^{1,\bullet+k}(\cM)$.  

Consider $U\subseteq M$ a chart of $\cM$ with coordinates 
$\{q^{\beta_i}_i \ | \ 0\leq i\leq n, \ 1\leq \beta_i\leq m_i\}$. Then, on $\pi^{-1}(U)\subseteq T^*[k]\cM$, we have 
coordinates 
$$\{q^{\beta_i}_i, p^{i}_{\beta_i}\ | \ \ 0\leq i\leq n, \ 1\leq \beta_i\leq m_i\}
\quad \text{of degree}\quad |p^i_{\beta_i}|=k-i,$$
 where $p_{\beta_i}^i$ are the linear coordinates corresponding to the section 
$\D q^{\beta_{i}}_{i}\in \Gamma_{i-k}(T^*[k]\cM)=\Omega^{1,i}(\cM).$ If $V\subseteq M$ is 
another chart with coordinates $\{\widehat{q}^{\beta_i}_i \ | \ 0\leq i\leq n, \ 1\leq \beta_i\leq m_i\}$
with $U\cap V\neq \emptyset$, 
then the linear cotangent coordinates transform by the rule
\begin{equation}\label{Eq:MomTrans}
p^i_{\alpha_i}=\sum_{j=i}^{n}\sum_{\beta_j=1}^{m_j} (-1)^{i(j-i)}\widehat{p}_{\beta_j}^j
\frac{\partial \widehat{q}^{\beta_j}_j}{\partial q^{\alpha_i}_i},
\end{equation}
see Remark \ref{Rem. TransForms}.  
Note that for $k>n$, the body of the shifted cotangent is just $M$, i.e.  $|T^*[k]\cM|=|\cM|$, 
while 
$|T^*[n]\cM|$ is a vector bundle over $M$ for $n=k$, which is exactly the vector bundle determined by the 
degree $n$ coordinates. This can be proven by using the geometrization of $n$-manifolds
given in Remark \ref{rem:geo}. 

\begin{remark}\label{Rem: formvb}
Let $(\cM, \omega)$ be a symplectic $n$-manifold. With the help of the shifted cotangent bundle 
we reinterpret the map $\omega^\flat\colon \fX^{1,\bullet}(\cM)\to \Omega^{1,\bullet+n}(\cM)$ as a 
graded vector bundle map $\omega^\flat\colon T\cM\to T^*[n]\cM$.   
\end{remark}

Using the coordinate description, one can immediately see that the $1$-form 
	\begin{align*}
	\tau=\sum_{i=0}^n\sum_{\alpha_i=0}^{m_i} p^i_{\alpha_i} \D q^{\alpha_i}_i \in \Omega^{1,k}(T^*[k]\cM) 
	\end{align*}
is globally defined. It is called the Liouville $1$-Form.  

\begin{proposition}\label{prop:cot}
   Let $k\geq n$ and consider the $k$-shifted cotangent bundle $T^*[k]\cM$ of an $n$-manifold $\cM$. 
    \begin{enumerate}
        \item $(T^*[k]\cM,\omega_{\can}=-\D \tau)$ is a degree $k$ symplectic manifold. 
            \item Let $U$ be a chart of $\cM$ with coordinates $\{ q_i^{\alpha_i}\}$ and consider the 
            canonical coordinates $\{ q_i^{\alpha_i}, p^i_{\alpha_i}\}$ on the $k$-shifted cotangent. 
            In this chart, the following expressions hold:
        \begin{itemize}
            \item $\omega_{\can}|_U=-
             \sum_{i=0}^n\sum_{\alpha_i=0}^{m_i}\D (p^i_{\alpha_i} \D q^{\alpha_i}_i)
             =-\sum_{i=0}^n\sum_{\alpha_i=0}^{m_i} \D p^i_{\alpha_i}\wedge \D q^{\alpha_i}_i.$
            \item $\{q_i^{\alpha_i}, p^j_{\alpha_j}\}_{\can}=-\delta_{ij}\delta_{\alpha_i \alpha_j}, 
            \quad \{q_i^{\alpha_i}, q_j^{\alpha_j}\}_{\can}
            =\{p^i_{\alpha_i}, p^j_{\alpha_j}\}_{\can}=0 $
            \item $X_{q_i^{\alpha_i}}=-\frac{\partial}{\partial p^i_{\alpha_i}}, \quad X_{p^i_{\alpha_i}}
            =(-1)^{i(k-i)}\frac{\partial}{\partial q_i^{\alpha_i}}.$
                \end{itemize}     
        \item The map  $J\colon \fX^{1,\bullet+k}(\cM)\to C^\bullet_{T^*[k]\cM}$ defined on coordinates 
        $\{ q_i^{\alpha_i}\}$ by 
        	\begin{align*}
        	J(\frac{\partial}{\partial q_i^{\alpha_i}})=(-1)^{i(k-i)}p_{\alpha_i}^i
        	\end{align*}
        extends to a globally defined left $C_\cM$-module morphism fulfilling   
        	\begin{align*}
        	\{J(X),\pi^* f\}_{\can}=\pi^* X(f) \quad \text{ and } \quad J([X,Y])=\{J(X),J(Y)\}_{\can}. 
        	\end{align*}
        \item $J$ extends to an algebra morphism from the $k$-shifted multivector fields, 
        see Remark \ref{rmk:mult}, to functions
        $$J\colon \bigoplus_i \sym^i(\mathfrak{X}^{1,\bullet}(\cM)[-k]) \to  C^\bullet_{T^*[k]\cM}$$
        and sends the Schouten bracket to the Poisson bracket.  Moreover, it is an isomorphism if $k>n$. 

    \end{enumerate}
\end{proposition}

\begin{proof}
    Items a., b. can be deduced from the coordinate expression for $\omega=-\D \tau$. Point c. is 
    proved using  the transformation rule from Equation \eqref{Eq:MomTrans} and b.
    Finally, d. holds because of c. 
\end{proof}

\begin{remark}
Since the sections of the cotangent can be also considered as a right module (see Remark \ref{Rem. TransForms})
 and the objects $\tau$, $\omega_{\can}$ and $J$ are defined by local coordinated expressions, which in turn come 
 from the choice of treating the sections of $T^*[k]\cM$ as a left module. One can ask if this choice has an 
 influence on the later work. In fact, if one uses the right module convention, one ends up with a symplectomorphic 
 version of the shifted cotangent bundle.   
\end{remark}

At the end of this section we introduce two special examples of Lagrangian submanifolds of $T^*[k]\cM$.  

\begin{example}[Conormal bundles]
 Let $j\colon \cN\hookrightarrow \cM$ a submanifold of the $n$-manifold $\cM$. For $k\geq n$, the \emph{$k$-shifted 
 conormal bundle of $\cN$} is the vector bundle  $N^*[k]\cN\to \cN$ whose sections are     
$$C_\cN\otimes_{j^{-1}C_\cM}j^{-1}\Omega^{1,\bullet+k}_\cN(\cM)$$
as constructed in Example \ref{ex:conor}. This shows that there is an exact sequence of vector bundles over $\cN$
$$0\to N^*[k]\cN\to j^*T^*[k]\cM\to T^*[k]\cN\to 0$$
that can also be taken as the definition of the $k$-shifted conormal bundle of $\cN$. From this exact sequence we 
get that $N^*[k]\cN$ is a submanifold of $T^*[k]\cM$ with $\Tdim N^*[k]\cN=\frac{1}{2}\Tdim T^*[k]\cM$. To see 
that it is coisotropic, we observe that 
$$\cI_{N^*[k]\cN}=\langle f, J(X_f) \ | \ f\in \cI_\cN\rangle.$$
We conclude that conormal bundles are Lagrangian submanifolds.
\end{example}

An application of the conormal bundle is the following.

\begin{proposition}\label{prop:isolag}
Let $(\cM, \omega)$ be a symplectic $n$-manifold and $j\colon \cL\hookrightarrow \cM$ a Lagrangian submanifold. Then 
the symplectic form $\omega$ induces the following isomorphisms of vector bundles over $\cL$:
\begin{equation*}
    \omega^\flat\at{\cL}\colon T\cL\to N^*[n]\cL\quad\text{and}\quad\chi\colon \nu(\cM, \cL)\to T^*[n]\cL.
\end{equation*} 
\end{proposition}

\begin{proof}
    Since $\cL$ is a Lagrangian submanifold, we know by Proposition \ref{Prop: EquLag} that $\omega^\flat\colon  
    \fX^{1,\bullet}_{\cL}(\cM)\to \Omega^{1,\bullet+n}_\cL(\cM)$ is an isomorphism. Therefore, if we pullback 
    this sheaf to $\cL$, we get an isomorphism
    $$\omega^\flat\colon C_\cL\otimes_{j^{-1}C_\cM}j^{-1}\fX^{1,\bullet}_{\cL}(\cM)\to C_\cL\otimes_{j^{-1}C_\cM}j^{-1}
    \Omega^{1,\bullet+n}_\cL(\cM).$$
    By definition, $\Gamma T\cL=C_\cL\otimes_{j^{-1}C_\cM}j^{-1}\fX^{1,\bullet}_{\cL}(\cM)$ and $\Gamma N^*[n]\cL
    =C_\cL\otimes_{j^{-1}C_\cM}j^{-1}\Omega^{1,\bullet+n}_\cL(\cM),$ so we get the first isomorphism. For the second, 
    observe that by Remark \ref{Rem: formvb} the morphism $\omega^\flat\colon T\cM\to T^*[n]\cM$ is an isomorphism. So we 
    get the exact sequence of vector bundle morphisms
    \begin{center}
	\begin{tikzcd}
	0\arrow[r]&T\cL \arrow[r]\arrow[d, "\omega^\flat\at{\cL}"] & j^*T\cM
	\arrow[d, "\omega^\flat"]\arrow[r] &\nu(\cM, \cL)\arrow[d, "\chi"]\arrow[r]&0 \\
	0\arrow[r]&N^*[n]\cL\arrow[r] &  j^*T^*[n]\cM\arrow[r]& T^*[n]\cL\arrow[r]& 0.
	\end{tikzcd}
	\end{center} 
 Therefore, $\chi\colon \nu(\cM, \cL)\to T^*[n]\cL$ is also an isomorphism.
\end{proof}

\begin{example}[Graph of closed $1$-forms]\label{ex:graph}
Let $\cM$ be an $n$-manifold and  $\alpha\in \Omega^{1,k}(\cM)=\Gamma_0T^*[k]\cM$. Then the map 
$\alpha\colon \cM\to T^*[k]\cM$ is an embedding and we get a submanifold of the $k$-shifted 
cotangent bundle, which is denoted by $\graph(\alpha)$. This submanifold has $\Tdim \graph(\alpha)
=\frac{1}{2}\Tdim T^*[k]\cM,$ 
thus we just need to check that the vanishing ideal $\cI_{\graph(\alpha)}$ is closed under the 
Poisson bracket (Definition \ref{def:lag}). For $X\in\fX^{1,\bullet}(\cM)$ we have
 \begin{align*}
	\alpha^* J(X)=\iota_X\alpha\quad \text{and thus}\quad \alpha^*(J(X)-\pi^*\iota_X\alpha)=0. 
    \end{align*}	
Therefore, the vanishing ideal to the submanifold $\graph(\alpha)$ is given by
$$\cI_{\graph(\alpha)}=\langle J(X)-\pi^*\iota_X\alpha\ | \ X\in \fX^{1,\bullet}(\cM)\rangle.$$
By the properties of the Poisson bracket given in Proposition \ref{prop:cot}, we get for $X,Y \in\fX^{1,\bullet}(\cM)$:
\begin{equation*}
    \begin{split}
  \{ J(X)-\pi^*\iota_X\alpha, J(Y)-\pi^*\iota_Y\alpha\}_{\can}=&J([X,Y])-
  \pi^*\Lie_X\iota_Y\alpha+(-1)^{|X||Y|}\pi^*\Lie_Y\iota_X\alpha      \\
  =& J([X,Y])-\pi^*\iota_{[X,Y]}\alpha-\pi^*\iota_X\iota_Y\D \alpha.
    \end{split}
\end{equation*}
Hence, $\graph(\alpha)$ defines a Lagrangian submanifold if and only if $\D\alpha=0$ and in this case 
(and $k\neq 0$) we get 
that $\alpha=\D f$ for some $f\in C^k_\cM$ by Proposition \ref{sym=ham}. 
\end{example}

\section{Weinstein's Lagrangian tubular neighbourhood theorem}

In this section we prove different normal form theorems for graded manifolds. The most relevant 
for us is a graded version of the Lagrangian tubular neighbourhood theorem, which was originally 
stated by Weinstein in \cite{Weinstein} and is sometimes also referred to as Weinstein's Lagrnagian tubular neighbourhood theorem.
Most of the techniques used in the original proof can 
be applied directly to the context of graded symplectic manifolds. Partial versions of this result 
in the context of supermanifolds or algebraic derived stacks can be found in the literature, 
see e.g. \cite{saf:lag, scw:bv}.

\subsection{Tubular neighbourhood theorem.}
One of the key ingredients of the classical Lagrangian tubular neighbourhood theorem is a tubular 
neighbourhood theorem for general submanifolds. Since we are not aware of such a result in graded geometry, 
we use this subsection to state and prove it.

\begin{theorem}\label{Thm: TubNei}
Let $j\colon \cN=(N, C_\cN)\to \cM$ be a submanifold of an $n$-manifold $\cM=(M, C_\cM)$. Then there 
exists $U\subseteq M$ a tubular neighbourhood  of $N$ and a diffeomorphism $\Psi\colon \cM\at{U}\to 
\nu(\cM, \cN)$ such that $j$ is identified with the zero-section.  
\end{theorem} 

The proof will follow an inductive argument, but first we need a tubular neighbourhood theorem in the category of 
vector bundles, which coincides with the category of $1$-manifolds. 

\begin{lemma}\label{lem:tub1}
    Let $j\colon (D\to N)\hookrightarrow (E\to M)$ be a subbundle. Then there exists $U\subseteq M$ a 
    tubular neighbourhood of $N$ and a vector bundle isomorphism $\Psi\colon E\at{U}\to \nu(E, D)$ such 
    that $j$ is identified with the zero-section. 
\end{lemma}

\begin{proof}
  The normal bundle $\nu(E,D)$ has the structure of a double vector bundle 
    \begin{center}
	\begin{tikzcd}
	\nu(E,D) \arrow[r]\arrow[d] & D\ar[d] \\
	\nu(M,N)\arrow[r, "p"] &  N
	\end{tikzcd}
	\end{center}
 with core given by $F=(j^*E)/D\to N$. Therefore, choosing a decomposition of the double vector
  bundle we get a vector bundle isomorphism $\nu(E,D)\cong p^*D\oplus p^*F$.  Let $U\subseteq M$ 
  be a tubular neighbourhood of $N$ with diffeomorphism $\psi\colon \nu(M,N)\to U$. Since $\psi$ is smoothly 
  homotopic to $j\circ p$, we get that $$p^*j^*E\at{U}\cong E\at{U},\quad \text{so there is}\quad 
  \Psi\colon p^*D\oplus p^*F \overset{\sim}{\to} E\at{U}.$$ Thus we get a vector bundle isomorphism 
  $\Psi\colon E\at{U}\to \nu(E,D)$ which identifies $j$ with the zero-section as desired.     
\end{proof}

Since the category of vector bundles is equivalent to the category of $1$-manifolds, we get that the 
above Lemma proves Theorem  \ref{Thm: TubNei} for $n=1$. 

\begin{proof}[Theorem \ref{Thm: TubNei}]
 We are going to prove the theorem by induction on the degree of the manifold using the tower of affine 
 fibrations \eqref{tower}. For $n=1$ the result is proven in the previous Lemma. Assume that the statement is true
  for $(n-1)$-manifolds. Let $j\colon \cN\hookrightarrow \cM$ 
 be a submanifold of the $n$-manifold $\cM$, then we get the following square
 \begin{equation}\label{eqc1}
	\begin{tikzcd}
	\cM  \arrow[r, twoheadrightarrow] & \cM_{n-1} \\
	\cN \arrow[u, hookrightarrow, "j"] \arrow[r, twoheadrightarrow] & \cN_{n-1} \arrow[u, hookrightarrow, "j_{n-1}"]  .
	\end{tikzcd}
 \end{equation}
By the inductive hypothesis, we obtain a tubular neighbourhood $U\subseteq M$ of $N$ 
in $M$ and a diffeomorphism $\Psi_{n-1}\colon \nu(\cM_{n-1}, \cN_{n-1})\to  \cM_{n-1}\at{U}$, 
which intertwines $j_{n-1}$ and  the zero-section.  In order to extend $\Psi_{n-1}$ to $\cM\at{U}$ 
we first observe, that on degree $n$ functions we get the diagram
\begin{center}
	\begin{tikzcd}
	0 \arrow[r]& C^n_{\cM_{n-1}}  \arrow[d, twoheadrightarrow, "j_{n-1}^\sharp"] 
	\arrow[r]& C^n_\cM\arrow[r] \arrow[d, twoheadrightarrow, "j^\sharp"]& 
	\frac{C^n_\cM}{C^n_{\cM_{n-1}}} \arrow[r] \arrow[d, twoheadrightarrow, "j_{n}^\sharp"]& 0  \\
	 0 \arrow[r]&j_*C^n_{\cN_{n-1}} \arrow[r]& j_* C^n_\cN \arrow[r]& \frac{j_*C^n_\cN}{j_*C^n_{\cN_{n-1}}} \arrow[r] &0.
	\end{tikzcd}
\end{center}
By the very definition of graded manifolds, we get that all the involved sheaves are sections of some vector
 bundle over $M$ and $N$, respectively. Therefore we can construct horizontal splittings that commute 
 with $j^\sharp$. Using these splittings the diagram \eqref{eqc1} becomes an inclusion of vector bundles 
 (over some graded manifold). Thus, we can apply Lemma \ref{lem:tub1} to the vector bundle 
 $\cM\at{U}\to \cM_{n-1}\at{U}$ and together with Lemma \ref{Lem: NormVec} we 
 get the desired tubular neighbourhood. 
\end{proof}

For our purposes, we need a special tubular neighbourhood, namely a tubular neighbourhood with
trivial \emph{normal derivative}: given a tubular neighbourhood of a submanifold 
 $j\colon \cN\hookrightarrow \cM$, i.e. a diffeomorphism $\Phi\colon \nu(\cM,\cN)\to \cM\at{U}$ for some open 
 $U\subseteq |\cM|$, we can build a vector bundle isomorphism  
\begin{align*}
\tilde{\Phi}\colon \nu(\nu(\cM,\cN),\cN)\to \nu(\cM,\cN)
\end{align*} 
covering the identity on $\cN$. Precomposing this with the canonical isomorphism $\nu(\cM,\cN)\cong \nu(\nu(\cM,\cN),\cN)$ 
from Lemma \ref{Lem: NormVec}, we get a vector bundle automorphism 
	\begin{align*}
	\nu(\Phi)\colon \nu(\cM,\cN)\to \nu(\cM,\cN),
	\end{align*}
which we call the normal derivative of $\Phi$. 

\begin{proposition}\label{Prop: Norm=id}
Let $j\colon \cN\hookrightarrow \cM$ be a submanifold, then there exists a tubular neighbourhood 
$\Phi\colon \nu(\cM,\cN)\to \cM\at{U}$ for an open $U\subseteq M$, such that $\nu(\Phi)=\id$.
\end{proposition}

\begin{proof}
Let us choose a tubular neighbourhood $\Psi\colon \nu(\cM,\cN)\to \cM\at{U}$, then $\Phi:=\Psi\circ\nu(\Psi)^{-1}$ 
is still a tubular neighbourhood and we have that 
	\begin{align*}
	\nu(\Phi)=\nu(\Psi)\circ\nu(\Psi)^{-1}=\id.
	\end{align*}	 
The latter can be seen by a straightforward computation and by using that for a vector bundle automorphism 
$\Psi\colon \cE\to\cE$ covering the identity  one has 
$(\Psi^{-1})^*\circ s^v\circ\Psi^*=(\Psi(s))^v$ for all sections $s\in \Gamma\cE$. 
\end{proof}

We close this subsection with a tiny observation which will be crucial later on to simplify computations. 

\begin{corollary}\label{Cor: adaptedcoord}
Let $j\colon \cN\hookrightarrow \cM$ and let $\Phi\colon \nu(\cM,\cN)\to \cM\at{U}$ be a tubular 
neighbourhood with $\nu(\Phi)=\id$. Let 
$\{x^{\alpha_i}_i,z^{\beta_j}_j\}$ be coordinates of $\nu(\cM,\cN)$ and let $\{x^{\alpha_i}_i,y^{\beta_j}_j=(\Phi^{-1})^*z^{\beta_j}_j\}$ 
be the corresponding coordinates on $\cM$. The linear coordinates corresponding to the 
local basis $\{[\tfrac{\partial}{\partial y^{\beta_j}_j}]\}$ of $\Gamma\nu(\cM,\cN)$ are given by $\{z^{\beta_j}_j\}$.  
\end{corollary}

\begin{proof}
Let $\{x^{\alpha_i}_i,z^{\beta_j}_j\}$ be coordinates of $\nu(\cM,\cN)$ and denote by 
$\{b_{\beta_j}^j\}$ the corresponding local basis of the sections 
of $\nu(\cM,\cN)$. Note that the vertical lift of $b_{\beta_j}^j$ is given by $\tfrac{\partial}{\partial z^{\beta_j}_j}$. 
Thus, we have that 
	\begin{align*}
	[\tfrac{\partial}{\partial y^{\beta_j}_j}]&
	 = [(\Phi^{-1})^*\circ \tfrac{\partial}{\partial z^{\beta_j}_j}\circ \Phi^*]
	 =\tilde\Phi[ \tfrac{\partial}{\partial z^{\beta_j}_j}]=\tilde\Phi((b_{\beta_j}^j)^v)\\&
	 =\nu(\Phi)(b_{\beta_j}^j)=b_{\beta_j}^j.
	\end{align*}	  
\end{proof}

\begin{remark}
We strongly believe that such a graded version of a tubular neighbourhood theorem 
may be helpful for many other (deformation) problems, such as for Lie subalgebroids inside Lie algebroids, etc, which 
are not directly connected to the realm of graded symplectic geometry.
\end{remark}

\subsection{The Darboux-Weinstein Theorem.}

The second ingredient that we need for the proof of the Lagrangian tubular neighbourhood is a Darboux 
theorem along submanifolds. In the context of  symplectic supermanifolds this result is partially discussed 
in \cite{scw:bv, Schwarz}. It was shown in \cite{vai:morse} that the classical Darboux theorem is a particular 
instance of the equivariant Morse Lemma for graded manifolds. Instead of using this interesting approach, we 
use the classical approach as outlined in \cite{can:sym}. 

\begin{theorem}\label{N-W-Splitting}
Let $\cM=(M,C_\cM)$ be an  $n$-manifold with $n>0$ and $j\colon N\hookrightarrow M$ a submanifold of 
the body. Let $[0,1]\subset (a,b)=I, \ \omega_0,\omega_1\in \Omega^{2,n}(\cM)$ two symplectic structures 
on $\cM$ and let
$\omega_t\in \Omega^{2,n}(\cM)$
be a smooth path of closed $2$-forms of degree $n$ connecting $\omega_0$ and $\omega_1$.
If  for all $p\in N, \ (\omega_t)_p$ is non-degenerate  then there exists $U\subseteq M$ an open 
neighbourhood  of $N$ and a symplectomorphism  $\Phi\colon (\cM|_U,\omega_0)\to (\cM|_U,\omega_1)$.
\end{theorem}

\begin{proof}
Since we assume that $n>0$, Proposition \ref{sym=ham} ensures that there exists a smooth curve of 
degree $n$ $1$-forms $\alpha_t$, such that  $\frac{\D}{\D t}\omega_t=\D\alpha_t$. In particular it is given by
	\begin{align*}
	 \alpha_t=\frac{1}{n}\iota_{\cE_u} (\frac{\D}{\D t}\omega_t)\in\Omega^{1,n}(\cM), 
	\end{align*}
where $\cE_u$ is the Euler vector field of $\cM$ as in Example \ref{Ex: Eulervf}. Since for
 all $p\in M$ we have $\cE_u|_p=0$, we conclude $(\alpha_t)_p=0$.

By hypothesis $(\omega_t)_p$ is non-degenerate for all $p\in N$ and all $t\in I$. Thus there exists $U\subseteq M$
 open neighbourhood of $N$ where 
$(\omega_t)_p$ is non-degenerate for all $p\in U$ and all $t\in I$ (maybe after shrinking $I$). So, by 
Proposition \ref{Lem: Symplatp} we get that $\omega_t$ defines a path of symplectic forms on $\cM\at{U}$. 
Therefore, we can find a unique (time-dependent) vectorfield 
$X_t\in \mathfrak{X}^{1,0}(\cM|_U)$, such that 
	\begin{align*}
	\iota_{X_t}\omega_t=-\alpha_t.
	\end{align*}
Moreover, we have that $X_t|_p=0$ for all $p\in U$, since 
$$X_t(f)= \iota_{X_t}df=\iota_{X_t}\iota_{X_f}\omega= \pm \iota_{X_f}\iota_{X_t}\omega=\pm 
\iota_{X_f}\alpha.$$ 
Remark \ref{Rem: FlowVec0} ensures that the flow of $X_t$ exists for all $t\in I$ and gives a 
smooth curve of diffeomorphisms $\Phi_t\colon \cM\at{U}\to \cM\at{U}$ with 
\begin{align*}
	\frac{d}{dt}\Phi_t^*\omega_t 
	=\Phi_t^* (\mathcal{L}_{X_t}\omega_t +\frac{d}{dt} \omega_t)
	= \Phi_t^* (- d\alpha_t + \frac{d}{dt} \omega_t)=0. 
\end{align*}
For $\Phi_1$ we obtain $\Phi_1^*\omega_1=\omega_0$ 
and the claim is proven.	 
\end{proof}

\begin{example}
To illustrate that Theorem \ref{N-W-Splitting} is non-trivial even if the 
symplectic structures coincide on $M$ we consider the following example: the 
$2$-manifold $\bR^{1|2|1}$ with coordinates $\{x,q,p,y\}$ of degree $|x|=0$, 
$|q|=|p|=1$ and $|y|=2$ and the symplectic structures
	\begin{align*}
	\omega_0=dx\wedge dy+dq\wedge dp \ \text{ and } \ \omega_1
	=dx\wedge d(y+pq)+dq\wedge dp. 
	\end{align*}
Clearly, we have $(\omega_0)_r=(\omega_1)_r$ for all $r\in \mathbb{R}$, thus $\omega_0+t(\omega_1-\omega_0)$ 
is a smooth path connecting $\omega_0$ and $\omega_1$ between 
them. This is a path of symplectic structures and Theorem \ref{N-W-Splitting} ensures now that they are symplectomorphic. However, 
they are obviously not equal.
\end{example}

As an application of Theorem \ref{N-W-Splitting}, we obtain the graded 
version of the classical Darboux Theorem, compare with Theorem 5.3 in \cite{kos:gra} for supermanifolds.

\begin{corollary}[Darboux coordinates]\label{Cor: Darboux}
Let $(\cM,\omega)$ be a  symplectic $n$-manifold of dimension 
$m_0|\dots|m_n$ and $p\in M$. Then there is $U\subseteq M$ a chart of $\cM$ around $p$ 
and a symplectomorphism $\Phi\colon \cM\at{U}\to \cN$ where
\begin{equation*}
    \cN=\left\{\begin{array}{ll}
        \text{for}\quad n=2l+1,& (T^*[n]\bR^{m_0|\cdots|m_l}, \omega_{\can});   \\
        \text{for}\quad n=4l,& 
        (T^*[n]\bR^{m_0|\cdots|m_{2l-1}|\frac{m_{2l}}{2}}, \omega_{\can});\\
        \text{for}\quad n=4l+2, &
        (T^*[n]\bR^{m_0|\cdots|m_{2l}}, \omega_{\can})\times (\bR^{0|\cdots|0 | m_{2l+1}}, 
        \sum_{j=1}^{m_{2\ell+1}} \epsilon_j\D y^j\wedge \D y^j) \quad \text{with}\quad \epsilon_j=\pm 1.
    \end{array}\right.
\end{equation*}
\end{corollary}

\begin{proof}
Let $p\in M$ and let $U\subseteq M$ be a chart of $\cM$ around $p$ with local coordinates  
$\{x^{\alpha_i}_i | 0\leq i\leq n, \ 1\leq \alpha_i\leq m_i\}$. On $\cM\at{U}$ the symplectic 
form can be written as
	\begin{align*}
	\omega=\sum_{j,l=0}^n\sum_{\alpha_j,\alpha_l=1}^{m_j, m_l} \omega^{jl}_{\alpha_j\alpha_l} 
	\D x^{\alpha_j}_j\wedge \D x^{\alpha_l}_l \quad \text{for some}\quad \omega^{jl}_{\alpha_j\alpha_l}\in C_{\cM\at{U}}^{n-j-l}.
	\end{align*}	
 On $\cM\at{U}$ we define now the $2$-form of degree $n$ with constant coefficients 
	\begin{align*}
	\omega_0=\sum_{j,l=0}^n\sum_{\alpha_j,\alpha_l=1}^{m_j, m_l} (\omega^{jl}_
	{\alpha_j\alpha_l})^0(p) \D x^{\alpha_j}_j\wedge \D x^{\alpha_l}_l.
    \end{align*}	 
Therefore, $\omega_{p}=(\omega_0)_{p}$, so they are symplectomorphic in a 
neighbourhood around $p$ by Theorem \ref{N-W-Splitting}. Thus the only missing step is to 
send the matrix $(\omega^{jl}_{\alpha_j\alpha_l})^0(p)$ to its normal form. However this is an easy problem, 
since many of the slots are zero by degree reasons, and we leave the proof to the reader. 
\end{proof}

\begin{remark}
As a trivial consequence from the above corollary, we the deduce: if $(\cM, \omega)$ is a symplectic $n$-manifold 
that does not admit any (local) 
Lagrangian submanifolds, then $n=4l+2$.      
\end{remark}

\subsection{The Weinstein Lagrangian tubular neighbourhood theorem.}
We are now in the position to state and prove the graded version of Weinstein's Lagrangian tubular neighbourhood 
theorem for Lagrangian submanifolds. In the context of supermanifolds the statement already appeared  without a 
proof in \cite{AKSZ, scw:bv}.

\begin{theorem}\label{Thm: dgWeinstein}
Let $(\cM=(M,C_\cM), \omega)$ be a symplectic $n$-manifold and $i\colon \cL=(L,C_\cL)\hookrightarrow\cM$ a 
Lagrangian submanifold. Then there exists an 
open neighbourhood $U$ of $L$ in $M$ and a symplectomorphism $\Psi\colon T^*[n]\cL \to \cM\at{U}$ 
such that the following diagram commutes:
\begin{center}
\begin{tikzcd}
	T^*[n]\cL \arrow[rr, "\Psi"]&&\cM\at{U} \\
	&\cL\ar[ur,hookrightarrow, "i"]\ar[lu,hookrightarrow, "0_{\cL}"']&
\end{tikzcd}
\end{center}
\end{theorem}

\begin{proof}
By the  tubular neighbourhood Theorem \ref{Thm: TubNei} there exists $V\subseteq M$ a tubular 
neighbourhood of $L$  and a diffeomorphism 
	\begin{align*}
	\Phi\colon \nu(\cM,\cL)\to \cM\at{V}
	\end{align*}	 
 which maps the embedding $i\colon \cL\to \cM$ to the zero-section. Since $\cL$ is a Lagrangian submanifold 
 we can apply Proposition \ref{prop:isolag} and get the isomorphism 
 $$\mu=-\chi\colon \nu(\cM,\cL)\to T^*[n]\cL$$
 induced by the symplectic form $\omega^\flat$.
 It remains to compare the two different symplectic structures. Let $p\in L\subseteq M$ and 
 choose an open $W\subseteq M$ around $p$ and coordinates of $\cM$ $\{q^{\alpha_i}_i, y^{\beta_k}_k\}$ 
 adapted to $\cL$. Then  $\{q^{\alpha_i}_i, p^i_{\alpha_i}\}$ are coordinates on $(T^*[n]\cL)\at{W\cap L}$ 
 and $\{q^{\alpha_i}_i, z^{\beta_k}_k\}$ are coordinates on 
 $\nu(\cM,\cL)\at{W\cap L}$, where the $z^{\beta_k}_k$ are the linear coordinates dual to the local basis 
 $[\frac{\partial}{\partial y^{\beta_k}_k}]\in \Gamma\nu(\cM,\cL)\at{W\cap L}$, see $\S$\ref{vector-bundles}. 
 Moreover, we may assume that the tubular neighbourhood fulfils
 $\nu(\Phi)=\id$, see Proposition \ref{Prop: Norm=id}, and hence we may also assume that  
 $\Phi^* y^{\beta_k}_k=z^{\beta_k}_k$, due to Corollary \ref{Cor: adaptedcoord}.
 The morphism $\chi$ has the coordinate expression
 \begin{align*}
 \mu(\big[\tfrac{\partial}{\partial y_j^{\alpha_j}}\big])=
	-\sum_{k=0}^{n-i}\sum_{\beta_k=1}^{m_k} \D q^{\beta_i}_i i^*(\iota_{\tfrac{\partial}{\partial q^{\beta_i}_i}} 
	\iota_{\tfrac{\partial}{\partial y_j^{\alpha_j}}} \omega)
	=-\sum_{k=0}^{n-i}\sum_{\beta_k=1}^{m_k}(-1)^{i(n-i-j)}i^*(\iota_{\tfrac{\partial}{\partial q^{\beta_i}_i}} 
	\iota_{\tfrac{\partial}{\partial y_j^{\alpha_j}}} \omega)\D q^{\beta_i}_i.
 \end{align*}
Since the linear coordinates corresponding to  $[\frac{\partial}{\partial y^{\beta_k}_k}]$ are given 
by $z^{\beta_k}_k$, we get for $\mu$:
 $$\mu^* q_i^{\alpha_i}=q_i^{\alpha_i} 
 \quad\text{ and }\quad 
 \mu^* p^i_{\alpha_i}=\sum_{k=0}^{n-i}\sum_{\beta_k=1}^{m_k}-(-1)^{i(n-i-k)}z^{\beta_k}_k 
	i^*(\iota_{\tfrac{\partial}{\partial q^{\alpha_i}_i}} 
	\iota_{\tfrac{\partial}{\partial y_k^{\beta_k}}} \omega).$$
and finally for the canonical symplectic form, we get
 \begin{align*}
	\mu^*\omega_{\can} &=  - \sum_{i=0}^{n}\sum_{\alpha_i=1}^{m_i}  \D \mu^* p_{\alpha_i}^i 
	\wedge \D \mu^* q^{\alpha_i}_i 
	=  
	 \sum_{i=0}^{n}\sum_{k=0}^{n-i}\sum_{\alpha_i, \beta_k=1}^{m_i, m_k}\D (-1)^{i(n-i-k)}z^{\beta_k}_k i^*(\iota_{\tfrac{\partial}{\partial q^{\alpha_i}_i}}
	 \iota_{\tfrac{\partial}{\partial y_k^{\beta_k}}} \omega)\wedge\D q^{\alpha_i}_i\\&
	=  \sum_{i=0}^{n}\sum_{k=0}^{n-i}\sum_{\alpha_i, \beta_k=1}^{m_i, m_k}\D z^{\beta_k}_k \wedge\D q^{\alpha_i}_i i^*(\iota_{\tfrac{\partial}{\partial q^{\alpha_i}_i}} 
	\iota_{\tfrac{\partial}{\partial y_k^{\beta_k}}} \omega)
	 -z^{\beta_k}_k   \D q^{\alpha_i}_i\wedge \D i^*(\iota_{\tfrac{\partial}{\partial q^{\alpha_i}_i}} 
	 \iota_{\tfrac{\partial}{\partial y_k^{\beta_k}}} \omega).
	\end{align*}	 
 On the other hand, the symplectic form $\omega$ has the coordinate expression
 $$\omega=\sum_{i=0}^{n}\sum_{k=0}^{n-i}\sum_{\alpha_i,\alpha_k=1}^{m_i,m_k} \frac{1}{2}\omega^{ik}_{\alpha_i\alpha_k} \D q^{\alpha_i}_i \wedge\D q_k^{\alpha_k}
 +\sum_{i=0}^{n} \sum_{k=0}^{n-i}\sum_{\alpha_i,\beta_k=1}^{m_i,m_k} \omega^{ik}_{\alpha_i\beta_k} \D q^{\alpha_i}_i \wedge \D y_k^{\beta_k}
 +\sum_{i=0}^{n} \sum_{k=0}^{n-i}\sum_{\beta_i,\beta_k=1}^{m_i,m_k}\frac{1}{2} \omega^{ik}_{\beta_i\beta_k} \D y^{\beta_i}_i \wedge \D y_k^{\beta_k},$$
where the coefficients are given as $\omega^{ik}_{\alpha_i\alpha_k}=(-1)^{(i+k)(n-i-k)} (\iota_{\tfrac{\partial}{\partial q^{\alpha_k}_k}}
\iota_{\tfrac{\partial}{\partial q^{\alpha_i}_i}}\Phi^*\omega)$, etc. 
Since $\Phi^* y^{\beta_k}_k=z^{\beta_k}_k$, we find 
\begin{align*}\Phi^*\omega=\sum_{i=0}^{n}&\sum_{k=0}^{n-i}\sum_{\alpha_i,\alpha_k=1}^{m_i,m_k}\frac{1}{2} \Phi^*\omega^{ik}_{\alpha_i\alpha_k} \D q^{\alpha_i}_i \wedge\D q_k^{\alpha_k}
 +\sum_{i=0}^{n} \sum_{k=0}^{n-i}\sum_{\alpha_i,\beta_k=1}^{m_i,m_k}  \Phi^*\omega^{ik}_{\alpha_i\beta_k} \D q^{\alpha_i}_i \wedge \D z_k^{\beta_k}\\&
 +\sum_{i=0}^{n} \sum_{k=0}^{n-i}\sum_{\beta_i,\beta_k=1}^{m_i,m_k}\frac{1}{2}  \Phi^*\omega^{ik}_{\beta_i\beta_k} \D z^{\beta_i}_i \wedge \D z_k^{\beta_k},
\end{align*}

However, since $\cL$ is a Lagrangian submanifold, the symplectic structure 
$\omega$ at $p'\in L\at{W}$ is given by   
$$(\Phi^*\omega)_{p'}=\sum_{i=0}^{n}\sum_{k=0}^{n-i}\sum_{\alpha_i,\beta_k=1}^{m_i, m_k}  (\omega^{ik}_{\alpha_i\beta_k})^0(p') 
(\D q^{\alpha_i}_i)\at{p'}\wedge (\D y_k^{\beta_k})\at{p'}
+\sum_{i=0}^{n} \sum_{k=0}^{n-i}\sum_{\beta_i,\beta_k=1}^{m_i,m_k}\frac{1}{2}  
(\omega^{ik}_{\beta_i\beta_k})^0(p') (\D z^{\beta_i}_i )\at{p'} \wedge (\D z_k^{\beta_k})\at{p'}.$$
To be more precise, the first term vanishes since $\cL$ is Lagrangian and the rest 
of the equality follows since $(\Phi^*\omega^{ik}_{\alpha_i\beta_k})^0(p')=
(\omega^{ik}_{\alpha_i\beta_k})^0(p')$ and 
$(\Phi^*\omega^{ik}_{\beta_i\beta_k})^0(p')=(\omega^{ik}_{\beta_i\beta_k})^0(p')$ for 
all possible combinations of $i,k,\alpha_i,\beta_k$, etc.
This means in particular that $(\mu^*\omega_{\can})_p+ t(\Phi^*\omega-\mu^*\omega_{\can})_p$ 
is invertible for all $t\in \mathbb{R}$ and $p\in L$, using a classical block matrix argument. 
Therefore we can apply Theorem \ref{N-W-Splitting} and get $U \subseteq M$  a tubular neighbourhood 
of $L$ and a symplectomorphism between $(\Phi^{-1}\circ\mu)^*\omega_{\can}$ 
and $\omega$ on $\cM\at{U}$. 
\end{proof}

As an application we can see that the existence of a wide Lagrangian submanifold already 
determines the global structure of a symplectic $n$-manifold, as it is the case  
for supermanifolds \cite{scw:bv}. 
\begin{corollary}
Let $(\cM,\omega)$ be a symplectic $n$-manifold and let $i\colon \cL\hookrightarrow \cM$ 
be a wide Lagrangian submanifold, then $(\cM,\omega)$ is (non-canonically) symplectomorphic to $(T^*[n]\cL,\omega_{\can})$.
\end{corollary}


\section{Deformations of Lagrangian Q-submanifolds}
This section is the core of our work. We attach to each Lagrangian $NQ$-submanifold 
an $L_\infty$-algebra that controls its (formal) deformation theory.  We discuss in detail 
the gauge equivalences of this $L_\infty$-algebra. Finally, we allow more general 
deformations, namely we discuss the case when the homological 
vector field is deformed at the same time,
which generalizes the case of the simultaneous deformations of a coisotropic submanifold 
and the Poisson bivector field and the deformations of $H$-twisted Poisson structures from \cite{Zambon}.         

\subsection{The $L_\infty$-algebra of a Lagrangian $NQ$-submanifold.}\label{sec:Linf}
The aim of this section is to attach an $L_\infty$-algebra to each
Lagrangian $NQ$-submanifold of a degree $n$ symplectic $NQ$-manifold, 
see Appendix \ref{app:Linfty} for a 
quick reminder on $L_\infty$-algebras. We make heavy use of 
the Weinstein's Lagrangian tubular neighbourhood theorem proved in the previous section which involves non-trivial 
choices. Nevertheless, we will show  that the attached  
$L_\infty$-algebra is unique up to $L_\infty$-isomorphism. 
Let us start with a general observation:   

    \begin{proposition}\label{Thm: weakHomPoiss}
	Let $\cN$ be an $n$-manifold and assume that $(T^*[n]\cN, \omega_{\can}, Q=X_{\theta})$ 
	is a symplectic $NQ$-manifold. Then the maps
	$l^j\colon \bigwedge^j(C_\cN[n-1])\to C_\cN[n+1-j],$ 
		\begin{align*}
		l^j(f_1,\cdots,f_j)
		:= 
		(-1)^{\sum_{i=1}^j(j-i)(|f_i|-n)}0_\cN^*\{\cdots\{\theta,\pi^*f_1\}_{\can},\pi^*f_2\}_{\can},\cdots,\pi^*f_j\}_{\can},
		\end{align*}		 
	define a curved $L_\infty$-structure on $C_\cN[n-1]$ with 
	curvature $l^0=0_\cN^*\theta$. Moreover, the zero-section $0_\cN\colon \cN\hookrightarrow T^*[n]\cN$ is a 
	Lagrangian $NQ$-submanifold if and only if $l^0=0$.    
	\end{proposition}

    \begin{proof}
    We know that $(C_{T^*[n]\cN}[n],\{\cdot,\cdot\}_{\can})$ is a graded Lie algebra and 
    that $C_\cN[n]$ is an abelian subalgebra via  $\pi^*\colon C_{\cN}[n]\to C_{T^*[n]\cN}[n]$. The zero 
    section $0_\cN^*\colon C_{T^*[n]\cN}[n] 
    \to C_{\cN}[n]$ is a projection with the property 
    \begin{equation}\label{eq:zerosec}
        0_\cN^*\{f,g\}_{\can}=0_\cN^*\{\pi^*0_\cN^*f,g\}_{\can}+0_\cN^*\{f,\pi^*0_\cN^*g\}_{\can}    
    \end{equation}
     which can be seen in Darboux coordinates. This turns $(C_{T^*[n]\cN}[n],C_{\cN}[n],0_\cN^*)$ into a V-algebra, 
     see Definition \ref{def:Valgebra}. By Theorem \ref{Them: DerBra} we can use $\theta$ to endow 
     $C_{\cN}[n]$ with an $L_\infty[1]$-algebra 
     structure or, equivalently, $C_\cN[n-1]$ with a curved $L_\infty$-algebra structure given by the 
     above brackets. Using the definition of the 
     vanishing ideal $\cI_{\graph(0)}$ of the zero section from 
     Example \ref{ex:graph}, we can see that $Q=\{\theta,\cdot\}_{\can}$ is tangential 
     to the zero-section if and only if $l^0=0_\cN^*\theta=0$. 
    \end{proof}

	Using the Lagrangian tubular neighbourhood theorem, we attach to each Lagrangian $NQ$-submanifold 
	the structure of an $L_\infty$-algebra on its functions. Moreover we show  that the isomorphism class of 
	this $L_\infty$-algebra is independent of the involved choices.
	
	\begin{theorem}\label{Thm: UniqLinfty}
	Let $(\cM,\omega,Q)$ be a degree $n$ symplectic $NQ$-manifold  with $n\geq 1$ and 
	$j\colon \cL=(L, C_\cL)\hookrightarrow \cM$ a Lagrangian $NQ$-submanifold. 
	For any Lagrangian tubular neighbourhood of $\cL$, 
	$C_\cL[n-1]$ become an $L_\infty$-algebra with brackets given by Proposition \ref{Thm: weakHomPoiss}. 
	Moreover, every two $L_\infty$-algebra 
	structures on $C_\cL[n-1]$, obtained in this way, are (non-canonically) $L_\infty$-isomorphic.    
	\end{theorem}
	
	\begin{proof}
	First we use Theorem \ref{Thm: dgWeinstein} to identify  $\Psi\colon (\cM\at{U},\omega,Q)\overset{\sim}{\to}
	(T^*[n]\cL,\omega_\mathrm{\can}, \{\theta,\cdot\}_{\can})$ where $U\subseteq M$ is a neighbourhood of $L$. It
	follows from Proposition 
	\ref{Thm: weakHomPoiss} that $C_\cL[n-1]$ has the   
	structure of an $L_\infty$-algebra. 
 
    For two Lagrangian tubular neighbourhoods we obtain 
     $$(\cM\at{U},\omega,Q)\cong(T^*[n]\cL,\omega_\can, 
    \{\theta_i,\cdot\}_{\can})\quad \text{for}\quad i=0,1.$$  
	By the construction of the isomorphisms in the proof of Theorem \ref{Thm: dgWeinstein}, which uses Theorem \ref{N-W-Splitting}, 
	we can assume that we find a curve of symplectomorphisms $\Phi_t$, such that 
	$\Phi_0=\id$, $\Phi_t\circ 0_\cN=0_\cN$ and $\Phi_1^*\theta_1=\theta_0$. Then we 
	consider the infinitesimal generator
		\begin{align*}
		\frac{d}{dt}\Phi_t^*=\Phi_t^*\circ X_t
		\end{align*}		
	which fulfils $0_\cL^*\circ X_t=0$. Using the same techniques 
	as in the proof of  Theorem 3.2 from \cite{SchaetzCattaneo} we obtain the 
	desired $L_\infty$-isomorphism. More concretely, it is given by the recursive formula (adapted to our case):
	\begin{align*}
	U&_{n,t}(f_1,\dots,  f_n)
	:=\sum_{\sigma\in S_n}\sign(\sigma) \sum_{k\geq 1} \sum_{\mu_1+\dots+\mu_k=n-1}\frac{1}{nk!\mu_1!\dots\mu_k!}\\&
	0_\cL^*\{\cdots\{\Phi_t^*\pi^*f_{\sigma(1)},U_{\mu_1,t}(\pi^*f_{\sigma(2)},\dots,\pi^*f_{\sigma(\mu_1+1)})\}_{\can},\cdots\}_{\can},
	U_{\mu_k,t}(\pi^*f_{\sigma(\mu_1+\cdots\mu_{k-1}+2)},\dots,\pi^*f_{\sigma(n)})\}_{\can},
	\end{align*}	 
	where $U_{1,t}=0_\cL^*\circ \Phi_t^*\circ \pi^*=\id$. In \cite{SchaetzCattaneo} it is 
	shown that this collection of maps $U_{n,t}\colon \sym C_\cL[n] \to C_\cL[n]$ defines an 
	isomorphism of $L_\infty[1]$-algebra structures on $C_\cL[n]$ induced by $\Phi_t^*\theta_1$. 
	\end{proof}
	
	 \begin{remark}
    For degree $n=1$, i.e. for coisotropic submanifolds inside Poisson manifolds,  the construction of the $L_\infty$-algebra 
    from Theorem \ref{Thm: UniqLinfty} 
    was already observed in \cite{CattaneoFelder}. 
    For $n=2$, i.e. for Dirac structures with support inside Courant algebroids,  there is a partial result in \cite{Gualtieri}. 
    To be more precise, the authors just consider the case when the Dirac structure is supported over the whole base manifold. 
    In contrast to 
    \cite{CattaneoFelder, Gualtieri}, our construction of the $L_\infty$-algebra has a geometric interpretation, 
    therefore the $L_\infty$-isomorphism between different choices are also geometric.
    \end{remark}

    Let us summarize some observations of the structure of the $L_\infty$-algebra constructed on $C_\cL[n-1]$:
    \begin{itemize}
        \item Choosing canonical coordinates $\{q^{\alpha_i}_i,p^i_{\alpha_i}\}$ on $T^*[n]\cM$, one check that for 
        $f_1,\dots,f_k\in C_\cL[n-1]$ the structure maps of the $L_\infty$-algebra are given by 
        \begin{align}\label{Eq: LocLinfty}
	     l^k(f_1,\cdots, f_k)=\pm 
	     0_\cN^* \Big(\frac{\partial^k \theta}{\partial p^{j_1}_{\alpha_{j_1}}\cdots \partial p^{j_k}_{\alpha_{j_k}}}\Big)
	     \frac{\partial f_1}{\partial q_{j_1}^{\alpha_{j_1}}}\cdots \frac{\partial f_k}{\partial q_{j_k}^{\alpha_{j_k}}}.
	     \end{align}
        	
        \item This formula shows that the structure maps $l^j$ from Proposition \ref{Thm: weakHomPoiss} 
         are derivations of $C_{\cL}$ in each slot. 
        \item Since $\cL$ is an $NQ$-submanifold, it has an induced vector field that we denote by $Q\at{\cL}$.
         In this case one easily check the identity $l^1=Q\at{\cL}$. By the previous item, and 
         using standard arguments of the theory of $L_\infty$-algebras, we get that the $Q$-cohomology of $\cL$, i.e. 
         $H^{\bullet}_{Q|_\cL}(\cL)$, inherits a degree $n-1$ Poisson structure, see \cite{cat:poi}. Another consequence of Theorem 
         \ref{Thm: UniqLinfty} is that this Poisson structure is independent of the choice of the Lagrangian 
         tubular neighbourhood. In $\S$\ref{sec:eq-def} we will explain the meaning of $H^{n}_{Q|_{\cL}}(\cL)$.
        \item A priori, there is no obvious reason why there should be just finitely many non-trivial higher brackets. 
        However if we assume that $\cL$ is an $(n-1)$-manifold, then one can see by counting degrees that $l^j=0$ for $j>n+1$. 
    \end{itemize}

Recall from $\S$\ref{sec: Cotangent} that  the body of $T^*[n]\cL$ is a vector bundle over the body of $\cL$. Then the
 $L_\infty$-structure constructed above depends on the vertical $\infty$-jet of $\theta$ in the vector bundle direction 
 at the zero section, which is easily seen by the first point of the previous list. If $\cL$ is actually an 
 $(n-1)$-manifold, the body of  $T^*[n]\cL$ and the body of 
 $\cL$ coincide, thus the vertical $\infty$-jet of $\theta$ is just a polynomial in the variables of positive degrees 
 and hence  $\theta$ can be recovered completely from the $L_\infty$-structure on $C_\cL[n-1]$. This case has been 
 treated in the literature already, see \cite{pri:poi, Raj}.

    \begin{definition}[See \cite{Raj}]\label{Def: HomPoiss}
	A degree $n$ homotopy Poisson manifold is an $n$-manifold $\cN$ with 
	$\theta\in C^{n+2}_{T^*[n+1]\cN}$ such that $\{\theta,\theta\}_\can=0$. It is called 
	strict, if $0_\cN^*\theta=0$
	\end{definition}
    In this case Proposition \ref{prop:cot} ensures that 
    $$C^{k}_{T^*[n+1]\cN}=\bigoplus_{i=0}^k \fX^{i,k-i(n+1)}(\cN),$$ 
	 therefore we can decompose 
	 $$\theta=(H=\pi^0)+(X=\pi^1)+\pi^2+\cdots+\pi^{n+2}\quad \text{with}\quad \pi^i\in \fX^{i, (n+2)-i(n+1)}(\cN).$$  
     The article \cite{Raj} used the Schouten bracket with the multivector fields to define an $L_\infty$-algebra on the 
     functions of a homotopy Poisson manifold $C_\cN[n]$.
	
	\begin{proposition}[Mehta \cite{Raj}]\label{Thm: HomPoiss=HomLieAlg}
	    Let $\cN$ be an $n$-manifold and $\theta\in C^{n+2}_{T^*[n+1]\cN}$ be an strict homotopy Poisson 
	    structure. Then the operations 
	    $$l^1(f)=X(f),\qquad l^j(f_1, \cdots, f_j)
	    =(-1)^{\sum_{i=1}^j(j-i)(|f_i|-n-1)}[\cdots[\pi^{j},f_1], \cdots], f_j],$$
	    for $j=2,\cdots, n+2$ and $f_i\in C_\cN$, define an $L_\infty$-algebra structure on $C_{\cN}[n]$. 
	    Moreover, this $L_\infty$-algebra coincides with the one given in Proposition \ref{Thm: weakHomPoiss}.  
	\end{proposition}
	
	\begin{proof}
    It was proven in Proposition \ref{prop:cot} that the Schouten bracket corresponds
     to the canonical Poisson bracket under the map $J$. Thus the above statement directly 
     follows from Proposition \ref{Thm: weakHomPoiss}.    
	\end{proof}

Combining Theorem \ref{Thm: UniqLinfty} with Proposition \ref{Thm: HomPoiss=HomLieAlg}, we get
	
	\begin{corollary}\label{Prop:Lagsubmanifold}
	${}$
 \begin{itemize}
     \item Let $(\cN, \theta)$ be a degree $n$ homotopy Poisson manifold. The zero-section of $T^*[n+1]\cN$ 
     is a Lagrangian $NQ$-submanifold of $(T^*[n+1]\cN,\{\cdot,\cdot\}_{\can}, Q=\{\theta,\cdot\}_{\can})$ 
	if and only if $\cN$ is a strict degree $n$ homotopy Poisson manifold.
    \item Let $j\colon \cL\hookrightarrow(\cM, \omega, Q)$ be a degree $n$ Lagrangian $NQ$-submanifold with body $L=M$. Then any choice 
    of a Lagrangian tubular neighbourhood gives an strict degree $n-1$ homotopy Poisson structure on $\cL$. Moreover, 
    two different tubular neighbourhoods give $L_\infty$-isomorphic homotopy Poisson structures.
 \end{itemize}
	\end{corollary}	

\begin{remark}\label{rmk:sp}
    According to \cite{pri:poi}, a degree $n$ homotopy Poisson structure $\theta$ on a $NQ$-manifold $(\cN, Q)$, such that $X=Q$,
     is an example of \emph{$n$-shifted Poisson structures}. Additionally, our degree $n$ Lagrangian $NQ$-submanifolds are also examples 
     of \emph{$n$-shifted Lagrangians} as introduced in \cite{pym:shif}.  In derived algebraic geometry there is a well known principle, 
     see e.g. \cite[Corollary 2.41]{saf:poi}, saying that every $n$-shifted Lagrangian induces an $(n-1)$-shifted Poisson structure. 
     Therefore, the second item in Corollary \ref{Prop:Lagsubmanifold} can be seen as a concrete way of realizing this principle for 
     wide Lagrangian $NQ$-submanifolds.    
\end{remark}

\subsection{Deformations of Lagrangian $NQ$-submanifolds.}\label{sec:defq}

In the previous section for each $j\colon \cL\hookrightarrow(\cM, \omega, Q)$ a degree $n$ Lagrangian 
$NQ$-submaniold we attached an $L_\infty$-algebra structure on $C_\cL[n-1]$. Here we explain how 
this $L_\infty$-algebra controls (formal, $1$-parameter) deformations of the Lagrangian 
$NQ$-submanifold inside the symplectic $NQ$-manifold. 

\begin{definition}
   A \emph{$1$-parameter deformation of a Lagrangian $NQ$-submanifold} $j \colon \cL\hookrightarrow (\cM, \omega, Q)$ is a 
fibration $\cF\to I$, where $I$ is an interval containing $0$, together with a smooth map $I\colon \cF\to \cM$ such that: 
\begin{itemize}
    \item $I\at{\cF_t}\colon \cF_t\to \cM$ is an embedding for all $t\in I$,
    \item $\cF_0=\cL$ and $I\at{\cF_0}=j$,
    \item $I\at{\cF_t}:\cF_t\to \cM$ is a Lagrangian $NQ$-submanifold for all $t\in I$.
\end{itemize}
\end{definition}

We call a deformation \emph{small}, if the underlying embedding of bodies 
$|I_t|\colon |\cF_t|\to |\cM|$ is $C^1$-close to $|\cL|$ and the graded coordinates of
$\cL$ stay transversal in some reasonable sense: we want the deformation to be encoded as sections 
of the normal bundle also in the graded world. To be precise:

	\begin{theorem}\label{Thm: SmallDef}
        Let $(\cM,\omega)$ be a symplectic $n$-manifold with $n\geq 1$, let $j\colon \cL\to \cM$ be a Lagrangian submanifold and let  $\Upsilon\colon T^*[n]\cL\to \cM\at{U}$
        be a Lagrangian tubular neighbourhood. The map
        \begin{equation*}
        \begin{array}{c}
            \mathfrak{R}_\Upsilon\colon  \left\{ \begin{array}{c}
                 \text{Smooth curves}  \\
                 f_t\in C^n_\cL, \quad f_0=0
            \end{array} \right\}  \to  \left\{ \begin{array}{cc}
                 \text{Deformations of the}  \\
                 \text{Lagrangian submanifold } \cL
            \end{array} \right\} \quad f_t  \rightsquigarrow  \mathfrak{R}_\Upsilon(f_t):=\Phi_t\circ j\colon \cL\to \cM,\\   
        \end{array} 
        \end{equation*}
        where $\Phi_t$ is the flow of the vector field $-X_{(\Upsilon^{-1}\circ\pi)^*\dot{f}_t}\in \fX^{1,0}(\cM\at{U})$, 
        is injective. We call a deformation \emph{small}, if  it is in the image of $\mathfrak{R}_\Upsilon$ for some Lagrangian tubular neighbourhood $\Upsilon$.
	\end{theorem}

    \begin{proof}	
    First of all we notice that for a curve of functions $f_t\in C_\cL^n$ defined on $I$, the Hamiltonian 
    vector field $X_{\pi^*f}\in\vf^{1,0}(T^*[n]\cL),$ where $\pi:T^*[n]\cL\to \cL$, is given in 
    canonical cotangent coordinates by 
    \begin{align*}
    X_{\pi^*f_t}=-\sum_{i=0}^n\sum_{\alpha_i=1}^{m_i}(-1)^{(n-i)i}\pi^*\Big(\frac{\partial f_t}{\partial q_i^{\alpha_i}}\Big)\frac{\partial} {\partial p_{\alpha_i}^i} 
	\end{align*}     
    So its symbol, see Remark \ref{Rem: FlowVec0}, is a curve of vertical vector fields with fibre-wise constant 
    coefficients on $|T^*[n]\cL|\to |\cL|$ and as such it possess a 
    flow that exists for all $t\in I$.   
    Let us first check that $\mathfrak{R}_\Upsilon$ actually does give a deformation of $\cL$ as a Lagrangian submanifold.  
   Consider the Hamiltonian vector field $-X^{\can}_{\pi^*\dot{f}_t}=-\{ \dot{f}_t, \cdot\}_{\can}\in \fX^{1,0}(T^*[n]\cL)$ and denote its flow by $\Phi_t^{\can}$. 
    Via the map $J$ from Proposition \ref{prop:cot}, we can see 
   $\mathfrak{X}^{1,\bullet-n}(M)\subseteq C_{T^*[n]\cL}^\bullet$  and by Example \ref{ex:graph} we have that 
   $$\mathcal{I}_{\graph(\D f_t)}=\langle J(X)-\pi^*\iota_Xdf_t\ | \ X\in \fX^{1,\bullet}(M)\rangle=\langle J(X)+\{\pi^*f_t, J(X)\}_{\can}\ | \ X\in \fX^{1,\bullet}(M)\rangle$$  
   is the vanishing ideal of the a Lagrangian submanifold $\graph(\D f_t)\hookrightarrow T^*[n]\cL$.
   Note that, since $C_\cL$ is abelian for the canonical Poisson bracket and  $\{f_t,X\}_{\can},\dot{f}_t\in C_\cL,$ we get $\{\dot{f}_t,\{f_t,X\}_{\can}\}_{\can}=0$ 
   and therefore
   \begin{equation*} 
   \frac{d}{dt}(\Phi^{\can}_t)^* (X+\{f_t,X\}_{\can})=(\Phi^{\can}_t)^*(-\{\dot{f}_t,X\}_{\can}-\{\dot{f}_t,\{f_t,X\}_{\can}\}_{\can}+\{\dot{f}_t,X\}_{\can})=0
   \end{equation*} 	 
    This implies  $(\Phi^{\can}_t)^*\mathcal{I}_{\graph(df_t)}=\mathcal{I}_\cL$ and thus $\Phi^{\can}_t(\cL)=\graph(\D f_t)$. Since $\Upsilon$ is a symplectomorphism, 
    $\mathfrak{R}_\Upsilon$ is a one parameter family of Lagrangian submanifolds of $\cM$ as we claimed.  
    
    Let now $f_t,g_t\in C^n_\cL$ be smooth curves with $\mathfrak{R}_\Upsilon(f_t)=\mathfrak{R}_\Upsilon(g_t)$ 
    and let $\Phi_t$ and $\Psi_t$ be the flows of the Hamiltonian vector fields $X_{\pi^*f_t}$ and $X_{\pi^*g_t}$ on $T^*[n]\cL$. 
    We have that 
    	\begin{align*}
    	\Phi_t\circ 0_\cL =\Psi_t\circ 0_\cL \quad \iff \quad \Psi_t^{-1}\circ \Phi_t\circ 0_\cL=0_\cL. 
    	\end{align*} 
    Note that $ \Psi_t^{-1}\circ \Phi_t$ is the flow of $X_{\Phi_t^*(\pi^*f_t-\pi^*g_t)}$, which can be shown by deriving. Moreover, we conclude that 
    $0_\cL^*X_{\Phi_t^*(\pi^*f_t-\pi^*g_t)}=0$. To simplify this expression we note that 
    	\begin{align*}
    	\frac{\D}{\D t} \Phi_t^*\pi^*h= \Phi^*_t\{\pi^*f_t,\pi^*h\}=0
    	\end{align*}
    for all $h\in C_\cL$, since $C_\cL$ is an abelian subalgebra. This implies that $\Phi_t^*\pi^*=\pi^*$ 
    for all $t$ and thus $\Phi_t^*(\pi^*f_t-\pi^*g_t)=\pi^*f_t-\pi^*g_t$. 
    Finally, we choose canonical coordinates $\{q^{\alpha_i}_i,p^{i}_{\alpha_i}\}$ and compute 
    	\begin{align*}
    	X_{\pi^*(f_t-g_t)}=-\sum_{i=0}^n \sum_{\alpha_i=1}^{m_i} \pi^*\big(\frac{\partial (f_t-g_t)}
    	{\partial q^{\alpha_i}_i}\big)\frac{\partial}{\partial p^i_{\alpha_i}}
    	\end{align*}
    using Proposition \ref{prop:cot}. So we get $0_\cL^*X_{\pi^*(f_t-g_t)}=0$ if and only if $\D (f_t-g_t)=0$ 
    and by degree reasons $f_t=g_t$.
    \end{proof}	
   
    \begin{remark}
    Note that the Lagrangian tubular neighbourhood theorem only gives us a diffeomorphism of a neighbourhood $U$ of the zero 
    section of $T^*[n]\cL$ to a neighbourhood of  $\cL$ in $\cM$. So to be precise, we can only consider curves 
    of functions  $f_t\in C_\cL^n$, such that the image of $\D f_t\colon \cL\to T^*[n]\cL$ is contained in $U$.
    Nevertheless, since $f_0=0$, we can always achieve that by restricting the time interval, at least if $|\cL|$ is compact.  
    \end{remark}

    \begin{remark}
    Even though we specified a tubular neighbourhood in Theorem \ref{Thm: SmallDef}, we can even get a slightly 
    stronger statement, if $|\cL|$ is compact: let $\Psi\colon \cL\times (a,b)\to \cM$ be a deformation of $\cL$ in $\cM$ 
    as a Lagrangian submanifold. For the choice of a Lagrangian tubular neighbourhood 
    $\Phi\colon T^*[n]\cL\to \cM\at{U}$, we may assume that $\Psi$ is a deformation of the zero-section 
    for a possibly smaller interval $I'\subseteq (a,b)$. Since $\pi\circ \Psi(\argument,0)=\id_\cL$, 
    we know that $\chi_t:=\pi\circ\Psi(\argument,t)\colon \cL\to \cL$ is a smooth curve of diffeomorphisms, 
    after a possible shrinking of $I'$. This means in particular that 
    $\hat{\Psi}(\argument,t):=\Psi(\argument,t)\circ\chi^{-1}_t\colon \cL\to T^{*}[n]\cL$ is a map 
    fulfilling $\pi\circ \hat{\Psi}(\argument,t)=0$ and hence a smooth family of sections. By Example \ref{ex:graph}, 
    we even get that $\hat{\Psi}$ is induced by the graph of the differential of a smooth curve of functions $f_t\in C_\cL$.           
    \end{remark}
    
    Since we only consider small deformations, we concentrate on shifted cotangent bundles in the following lines.  
	
	\begin{corollary}\label{Prop: DefLagsilly}
	Let $f_t\in C_\cL^n$ be a smooth curve of functions such that $f_0=0$.
	The graph $\graph (\D f_t)\subseteq T^*[n]\cL$ is a Lagrangian $NQ$-submanifold at $t$, if and only if 
		\begin{align*}
		0_\cL^*\Phi_t^*\theta=0
		\end{align*}
	where $\Phi_t$ is the flow of $X_{\dot{f}_t}=-\{\dot{f_t},\cdot\}_{\can}\in \mathfrak{X}^{1,0}(T^*[n]\cL)$
	 at $t$. Thus, if $0_\cL^*\Phi_t^*\theta=0$ for all $t$, we have that 
	$t\mapsto \graph(df_t)$ is a deformation of Lagrangian $NQ$-submanifolds of the zero-section. 
	\end{corollary} 
	\begin{proof}
	   A Hamiltonian vector field $X_\theta$ is tangent to a Lagrangian submanifold if and 
	    only if the restriction of $\theta$ vanishes along the Lagrangian submanifold by Proposition \ref{Prop: EquLag}. Therefore 
	    $\graph(\D f_t)=\Phi_t\circ 0_\cL\colon \cL\to T^*[n]\cL$ defines a Lagrangian $NQ$-submanifold 
	    for $t\in I$ if and only if $$0_\cL^*\Phi_t^*\theta=0.$$
	\end{proof}

  In order to relate the $L_\infty$-algebra constructed in $\S$ \ref{sec:Linf} with smooth curves 
  parametrizing deformations of Lagrangian $NQ$-submanifolds, let us recall 
  the notion of Maurer-Cartan elements in a given $L_\infty$-algebra, see also Definiton \ref{Def: MC}.

\begin{definition}
Let $(\liealg{g},\{l^k\})$ be a curved $L_\infty$-algebra. An element
$\pi\in \liealg{g}^1$ is called a \emph{Maurer-Cartan} element if 
	\begin{equation}\label{eq:mc}
	\sum_{k=0}^\infty \frac{1}{k!} l^k(\overbrace{\pi,\dots, \pi}^{k\text{-times}})=0.    
	\end{equation}
We denote the set of Maurer-Cartan elements by $\mathrm{MC}(\mathfrak{g})$.
\end{definition}  
	
	Before we approach the general deformation problem, let us study the illustrative case of strict homotopy Poisson manifolds, 
	see Definition \ref{Def: HomPoiss}. Their advantage is that the $L_\infty$-algebra structure on its functions from Proposition 
	\ref{Thm: HomPoiss=HomLieAlg} has only finitely many non-vanishing higher brackets.  
	So the framework is as follows: let $\cN$ be an $(n-1)$-manifold and 
	$\theta\in C^{n+1}_{T^*[n]\cN}$, such that $\{\theta,\theta\}_{\can}=0$ and $0_\cN^*\theta=0.$ Then 
	$(T^*[n]\cN, \omega_{\can}, X_\theta)$ is a degree $n$ symplectic $NQ$-manifold and $0_\cN\colon \cN\to T^*[n]\cN$ a 
	Lagrangian $NQ$-submanifold by Corollary \ref{Prop:Lagsubmanifold}. 
    
    \begin{proposition}\label{Prop: Tangent to Lagrangian}
	Let $(\cN, \theta)$ be a strict homotopy Poisson manifold of degree $(n-1)$ and let $f_t\in C^n_{\cN}$ be a smooth curve with $f_0=0$. 
	Then $f_t$ defines a deformation of the Lagrangian $NQ$-submanifold via  $\graph(\D f_t)\colon \cN\to T^*[n]\cN$ 
	for the homological vector field $Q=X_\theta$ 
	if and only if  $f_t$ is a 
	Maurer-Cartan element in the $L_\infty$-algebra from Proposition \ref{Thm: HomPoiss=HomLieAlg} for all $t\in I$. 
	\end{proposition} 	
	
	\begin{proof}
	  Since $\cN$ has degree $(n-1)$, Proposition \ref{prop:cot} ensures that 
	$$C^{k}_{T^*[n]\cN}=\bigoplus_{i=0}^k \fX^{i,k-in}(\cN).$$
 Moreover, since $(\cN,\theta)$ is a strict homotopy Poisson manifold,
 $\theta$ decomposes as 
 $$\theta=\pi^1+\pi^2+\dots+\pi^{n+1}\quad \text{with}\quad  \pi^i\in \fX^{i,n+1-in}(\cN).$$
 The  Hamiltonian vector field of a smooth curve of functions $f_t\in C_\cL^n$ can be seen as an operator 
  $$X_{f_t}=\{f_t,\cdot\}_{\can}=[f_t, \cdot]\colon \fX^{i,k-in}(\cN)\to \fX^{i-1,k-(i-1)n}(\cN).$$  
They are, by degree reasons, locally nilpotent. The exponential is therefore well-defined by the usual formula
 $$\E^{-X_{f_t}}(g):=\sum_{k=0}^{\infty}\frac{(-1)^k}{k!} [f_t,[f_t,[f_t,\cdots,[f_t,g]\cdots]\quad \text{for}\quad g\in C^\bullet_{T^*[n]\cN}.$$
  Since $f_t$ and $\dot{f}_t$ Poisson-commute, $\E^{-X_{f_t}}$ fulfils 
    \begin{align*}
		\frac{d}{dt}\E^{-X_{f_t}}g =  \E^{-X_{f_t}}(-\{\dot{f}_t,g\}_{\can}).
	\end{align*}
 This implies that	$\E^{-X_{f_t}}$ is the pull-back with the flow of $-\{\dot{f},\cdot\}_{\can}$. By Corollary 
 \ref{Prop: DefLagsilly} we get that the curve $f_t$ defines a deformation of the Lagrangian 
 $NQ$-submanifold $0_\cN\colon \cN\hookrightarrow (T^*[n]\cN, \omega_{\can}, \{\theta,\cdot\}_{\can})$ if and only if 
		\begin{align*}
		0=0_\cN^*(\E^{-X_{f_t}}(\theta))
		= \sum_{k=1}^{n+1}\frac{(-1)^k}{k!} [f_t,[f_t,[f_t,\cdots,[f_t,\pi^k]\cdots]
		=\sum_{k=1}^{n+1}\frac{1}{k!}l^k(f_t,\dots,f_t)
		\end{align*}
	for all $t\in I$. So $f_t$ must be a Maurer-Cartan element for all $t\in I$.
	\end{proof}    

    Maurer-Cartan elements on the $L_\infty$-algebra defined by a homotopy Poisson manifold already 
    appeared in \cite{she:hom}, surprisingly there is no connection with deformation theory in that work.  
    An immediate consequence of the above Proposition together with Corollary \ref{Prop:Lagsubmanifold} is the following.

	\begin{proposition}
	Let $j\colon \cL\hookrightarrow (\cM, \omega, Q)$ be a degree $n$ Lagrangian $NQ$-submanifold with body $L=M$.
	 Small deformations of $\cL$ as Lagrangian $NQ$-submanifold are controlled by Maurer-Cartan elements 
	in $( C_\cL[n-1],l^j)$ for the choice of some tubular neighbourhood $\Psi$.
	\end{proposition}	
	
	\begin{proof}
	    Using the Lagrangian tubular neighbourhood Theorem \ref{Thm: dgWeinstein}, we get a symplectomorphism 
	    $\Psi\colon (\cM, \omega)\to (T^*[n]\cL,\omega_{\can})$. Then $\Psi_*Q=X_\theta$ for some 
	    $\theta\in C_{T^*[n]\cL}^{n+1}$ and since $L=M$,  $(\cL, \theta)$ is a degree $(n-1)$ strict 
	    homotopy Poisson structure by Corollary \ref{Prop:Lagsubmanifold}. So we can apply the previous 
	    Proposition and get the result. 
	\end{proof}

 In the general case the role of the $L_\infty$-algebra is bit more subtle, since we may have convergence 
 issues of the Maurer-Cartan equation. 
A way to circumvent this issue is given by the following statement.  
 
   \begin{theorem}\label{Thm:formDefProb}
   Let $j\colon \cL \to (\cM,\omega,Q)$ be a degree $n$ Lagrangian $NQ$-submanifold, let 
   $\Psi\colon T^*[n]\cL \to \cM\at{U}$ be a tubular neighbourhood and let 
   $f_t\in C_\cL^n$ be a curve of smooth functions, such that $f_0=0$ and $\mathfrak{R}_\Phi(f_t)$ 
   deforms $\cL$ as a Lagrangian $NQ$-submanifold. Then the formal Taylor series of $f_t$,      
   defined by 
   \begin{align*}
    \hat{f_t}:=\sum_{k=1}^\infty  \frac{1}{k!} (\tfrac{d^k f_t}{dt^k}\at{t=0}) \nu^k  \in \nu C^1_\cL[n-1][[\nu]]
   \end{align*}
   is a Maurer-Cartan element in the $L_\infty$-algebra $(C_\cL[n-1][[\nu]],\{l^i\}_{i\in \mathbb{N}})$, 
   where all the structure maps are extended $\nu$-linearly.  
   \end{theorem} 
   
   \begin{proof}
   Let $f_t\in C^n_{\cL}$ be a smooth 
   curve such that $f_0=0$. By Corollary \ref{Prop: DefLagsilly} we know that $f_t$ defines a 
   deformation of the Lagrangian $NQ$-submanifold via  $\graph(\D f_t) \colon \cL\to T^*[n]\cL$  
   if and only if $0_\cL^*\Phi_t^*\theta=0$,
    where $\Phi_t$ is 
    the flow of $-\{\dot{f_t},\argument\}_{\can}\in\fX(T^*[n]\cL)$,
    which means that in particular, we have that 
   	\begin{align*}
   	0= \frac{d^k}{dt^k}\at[\Big]{t=0}0_\cL^*\Phi_t^*\theta \ \text{ for all }\  k\geq 0.
   	\end{align*}
   	This is equivalent to $\hat{f_t}:=\sum_{k=0}^\infty \frac{\nu^k}{k!} \frac{d^k}{dt^k}\at{t=0}f$, 
   	where $\nu$ is a formal parameter, 
   	fulfilling the Maurer-Cartan equation: one can show by induction that 
   		\begin{align*}
   		\frac{d^k}{dt^k} \Phi_t^*\theta 
   		=\Phi^*_t \Big(\sum_{l=1}^k\frac{1}{l!} 
   		\sum_{i_1+\dots+i_l=k;\ i_j>0} \frac{(-1)^lk!}{i_1!\cdots i_l!} 
   		\{f_t^{(i_1)},\{\dots,\{f_t^{(i_l)},\theta\}_{\can}\cdots\}_{\can}\}_{\can}\Big)
   		\end{align*}
   	where $f^{(i)}_t$ indicates the $i$-th $t$-derivative of $f_t$. With this, one computes
   	 for $k\in\mathbb{N}$
   		\begin{align*}
   	    0= \frac{d^k}{dt^k}\at[\Big]{t=0}0_\cL^*\Phi_t^*\theta &= 
        \sum_{l=1}^k\frac{1}{l!} 
        \sum_{i_1+\dots+i_l=k;\ i_j>0} \frac{(-1)^lk!}{i_1!\cdots i_\ell!} 
        0_\cL^*\{f_0^{(i_1)},\{\dots,\{f_0^{(i_\ell)},\theta\}_{\can}\cdots\}_{\can}\}_{\can}\\&
        =\sum_{l=1}^k\frac{1}{l!} \sum_{i_1+\dots+i_l=k;\ i_j>0} 
        \frac{k!}{i_1!\cdots i_\ell!}l^\ell (f_0^{(i_1)},\cdots,f_0^{(i_\ell)}), 	
   	\end{align*}	
   	which is exactly order $k$ (in $\nu$) for the formal Maurer-Cartan equation of $\hat{f}$, i.e. of
   		\begin{align*}
   		\sum_{k=1}^\infty\frac{1}{k!} l^k (\hat{f_t},\cdots,\hat{f_t})
   		\end{align*}
   	for the $\nu$-linear extended brackets $l^{k}$.
	\end{proof}      

 \begin{remark}
 Theorem \ref{Thm:formDefProb} shows that the formal deformation problem of a Lagrangian $NQ$-submanifold is 
 controlled by the $L_\infty$-algebra constructed in Theorem \ref{Thm: weakHomPoiss} for the 
 choice of a tubular neighbourhood. 
Theorem \ref{Thm: UniqLinfty} states additionally that this choice is irrelevant. 
 \end{remark}
 
 \subsection{Gauge equivalence and trivial deformations.}\label{sec:eq-def}
 
 Having discussed the deformation theory for Lagrangian $NQ$-submanifolds, 
 a natural question to ask is when do we consider deformations to be equivalent? 
 Recall that in symplectic geometry one usually uses Hamiltonian isotopies, see e.g. \cite{oh:coi}, which is not suitable for 
 our purposes, since by Theorem \ref{Thm: SmallDef} we get that every small deformation of a Lagrangian submanifold is
 realized by the flow of a Hamiltonian vector field. Therefore if we use Hamiltonian isotopies every two small deformations would 
 be equivalent. However, in the previous section we have found that small (formal) deformations correspond to Maurer-Cartan 
 elements of an $L_\infty$-algebra. The set of Maurer-Cartan elements comes naturally equipped with an equivalence relation given by the gauge action, see Appendix 
 \ref{app:Linfty}. In this section we aim to describe gauge equivalent elements in detail and give them a geometric meaning. 
 
Recall that two Maurer-Cartan elements $f_0,f_1$ in an $L_\infty$-algebra $(\mathfrak{g}, \{ l^i\})$ are called 
\emph{gauge equivalent}, if there is a curve $f_t$ of Maurer-Cartan elements starting at $f_0$ and ending at $f_1$ together with 
a curve of degree $0$ elements $\lambda_t$, such that\footnote{Note that at the moment we forget about the convergence issues of this equation.} 
	\begin{equation}\label{eq:eqiv}
	    \frac{\D}{\D t}f_t=\sum_{k=0}^\infty \frac{1}{k!}l^{1+k}(f_t,\cdots,f_t,\lambda_t),
	\end{equation} 
see also Definition \ref{prop: EquMC}.

The Maurer-Cartan equation \eqref{eq:mc} and the equation \eqref{eq:eqiv} are generally non-linear. Their linearizations at $t=0$ are
\begin{equation*}
    l^1(\dot{f}_0)=0\quad\text{and}\quad \dot{f}_0=l^1(\lambda_0). 
\end{equation*}
This means that, if a curve $f_t$ of smooth functions with $f_0=0$ fulfils the 
Maurer-Cartan Equation \ref{eq:mc}, then its derivative at $t=0$ is closed with respect 
to the differential $l^1$. If it fulfils additionally Equation \eqref{eq:eqiv}, $\dot{f}_0$ is even exact.  

Solutions to these equations are called \emph{infinitesimal deformations} and \emph{infinitesimal gauge equivalence}. 
Therefore, as an immediate consequence of Theorem \ref{Thm:formDefProb}, we get

\begin{proposition}
Let $j\colon \cL\hookrightarrow(\cM, \omega, Q)$ be a Lagrangian $NQ$-submanifold. Infinitesimal 
deformations of $\cL$ as a Lagrangian $NQ$-submanifold modulo infinitesimal gauge equivalences are parameterized by $H^{n}_{Q|_\cL}(\cL).$
\end{proposition}

Let $(\cM, \omega, Q)$ be a degree $n$ symplectic $NQ$-manifold and $h_t\in C^{n-1}_\cM$ a smooth curve of functions. 
Then $\{Q(h_t), \cdot\}\in \fX^{1, 0}(\cM)$ is a degree zero Hamiltonian vector field that commutes with $Q$ because 
$[Q, X_{Q(h_t)}]=X_{Q^2(h_t)}=0$. Thus if $j\colon \cL\hookrightarrow (\cM, \omega, Q)$ is a Lagrangian $NQ$-submanifold 
we have that $\Psi_t\circ j\colon \cL\hookrightarrow (\cM, \omega, Q)$ defines a deformation of $\cL$ as a Lagrangian 
$NQ$-submanifold, where $\Psi_t$ is the flow of $\{Q(h_t), \cdot\}$. Motivated by this observation, we propose the following definition.

\begin{definition}\label{Def: TrivDef}
Let $j\colon \cL\to (\cM, \omega, Q)$ be a degree $n$ Lagrangian $NQ$-submanifold. Let 
$f_t\in C_{\cL}^n$ be a smooth curve such that $\mathfrak{R}_\Phi (f_t)$
 is a small deformation of $\cL$ as a Lagrangian $NQ$-submanifold. We say that $f_t$ 
induces a \emph{trivial deformation}, if there 
exists a curve of functions $h_t\in C_{\cM}^{n-1}$, such that 
	\begin{align*}
	\mathfrak{R}_{\Phi} (f_t)=\Psi_t\circ j\colon \cL\to \cM,
	\end{align*}
where $\Psi_t$ is the flow of $\{Q(h_t),\argument\}$.
\end{definition}

Before we discuss the relation of trivial deformations and gauge equivalences in our $L_\infty$-algebra, 
we need a technical result, relating flows of different kind of Hamiltonian vector fields:   

\begin{proposition}\label{Prop: equivalenceDef}
Let $\cN$ be an $n$-manifold and assume that $(T^*[n]\cN, \omega_{\can}, X_\theta)$ is a $NQ$-manifold for some 
$\theta\in C^{n+1}_{T^*[n]\cN}$. Let   $f_t\in C_{\cN}^n$ and $g_t\in C_{T^*[n]\cN}^n$ smooth curves and denote 
by $\Phi_t$  the flow of $-\{\pi^*\dot{f}_t,\cdot\}_{\can}$ and by
$\Psi_t$ the flow of $X_{g_t}=\{g_t,\cdot\}_{\can}$.  Then 
\begin{align*}
	\Phi_t \circ 0_\cN=\Psi_t\circ 0_\cN \quad \text{if and only if}\quad \dot{f_t}=-0_\cN^*\Phi_t^*g_t.
	\end{align*}
\end{proposition}

\begin{proof}
We compute 
	\begin{align*}
	\frac{\D}{\D t} \Psi_t^*\circ (\Phi_t^{-1})^*&
	=\Psi_t^*\circ \{g_t,\cdot\}_{\can}\circ(\Phi_t^{-1})^* + 
	\Psi_t^*\circ\{\pi^*\dot{f}_t,\cdot\}_{\can}\circ(\Phi_t^{-1})^*\\&
	=\Psi_t^*\circ(\Phi_t^{-1})^*\circ \{\Phi_t^*g_t+\Phi_t^*\pi^*\dot{f}_t,\cdot\}_{\can}\\&
	=\Psi_t^*\circ(\Phi_t^{-1})^*\circ \{\Phi_t^*g_t+\pi^*\dot{f}_t,\cdot\}_{\can},
	\end{align*}
where in the last step we used that $\Phi_t$ is the Hamiltonian flow of the pull-back of a function on 
$\cN$, so it preserves the projection $\pi$. The above equation tells us that  
$\Phi_t^{-1}\circ \Psi_t$ is the flow of the vector field 
$\{\Phi_t^*g_t+\pi^*\dot{f}_t,\cdot\}_{\can}$.
Therefore this Hamiltonian vector field  preserves the zero section, if and only if 
	\begin{align*}
	\dot{f}=-0_\cN^*\Phi_t^*g_t,
	\end{align*}
which can be seen, for example, in the canonical cotangent coordinates.  
\end{proof}

It is again illustrative to consider first the case of strict homotopy Poisson manifolds, i.e. wide Lagrangian
submanifolds.  

\begin{theorem}\label{Lem:HomPoissEqMC}
Let $(\cN,\theta)$ be a strict degree $n-1$ homotopy Poisson manifold. 
Let $f\in C_{\cN}^n$ be a Maurer-Cartan element in the
$L_\infty$-algebra given by Proposition \ref{Thm: HomPoiss=HomLieAlg}. 
Then there exists a smooth curve $f_t$ of degree $n$ functions on $\cL$ 
with $f_0=0$ and $f_1=f$, such that $f_t$ induces a trivial deformation 
of the Lagrangian $NQ$-submanifold 
$0_\cN\colon \cN\hookrightarrow (T^{*}[n]\cN, \omega_{\can}, X_\theta)$ if and only if 
$f$ is gauge equivalent to $0$ as Maurer-Cartan elements. 
\end{theorem}

\begin{proof}
Let us assume first that $f$ induces a trivial deformation of 
$0_\cN\colon \cN\hookrightarrow T^{*}[n]\cN$. By the definition of trivial deformations and 
Proposition \ref{Prop: equivalenceDef}, there exist smooth curves 
$h_t\in C_{T^*[n]\cN}^{n-1}$ and  $f_t\in C_{\cN}^{n}$, such that
	\begin{equation*}
	f_0=0, \quad 	f_1=f \quad\text{and}\quad 	\dot{f}_t
	=0_\cN^*\Phi_t^*\{h_t, \theta\}_{\can} =0_\cN^*\{\pi^*0_\cN^*\Phi_t^*h_t, \Phi_t^*\theta\}_{\can},   
	\end{equation*}
 where the last equality is obtained by equation \eqref{eq:zerosec}. Let us denote 
 $\lambda_t:=0_\cN^*\Phi_t^*h_t\in C^{n-1}_{\cN}$ and as in  the proof of Proposition 
 \ref{Prop: Tangent to Lagrangian}, we have that the flow $\Phi_t$ is given by
$\Phi_t^*=\E^{-\{f_t,\argument\}}$ and hence
	\begin{align*}
	\frac{\D}{\D t}f_t=\dot{f}_t=0_\cN^*\{\pi^*\lambda_t, \Phi_t^*\theta\}_{\can}=
	\sum_{k=0}^{n+1}\frac{(-1)^k}{k!} [\lambda_t, [f_t,[\dots,[f_t,\theta]\dots]]]
	= \sum_{k=0}^{n+1}\frac{1}{k!}l^{k+1}(f_t,\dots,f_t,\lambda_t).
	\end{align*}
Thus $f$ is gauge equivalent to $0$ as Maurer-Cartan elements, 
since $f_1=f$. 

The other direction is easier: let us 
assume we are given smooth curves $f_t\in C_{\cN}^n$ and $\lambda_t\in C_{\cN}^{n-1}$ fulfilling  
\begin{align*}
	\frac{\D}{\D t}f_t=\sum_{k=0}^\infty \frac{1}{k!}l^{1+k}(f_t,\cdots,f_t,\lambda_t).
	\end{align*}	  
We define the curve  $g_t=\{\lambda_t,\theta\}_{\can}\in C^n_{T^*[n]\cN}$ and compute
	\begin{align*}
	\frac{\D}{\D t}f_t &=\sum_{k=0}^\infty \frac{1}{k!}l^{1+k}(f_t,\cdots,f_t,\lambda_t)
	=\sum_{k=0}^{n+1}\frac{(-1)^k}{k!} [\lambda_t, [f_t,[\dots,[f_t,\theta]\dots]]]=\\&
	=\sum_{k=0}^{n+1}\frac{(-1)^k}{k!} 0_\cN^* [f_t,[\dots,[f_t,[\lambda_t,\theta]\dots]]]
	=\sum_{k=0}^{n+1}\frac{(-1)^k}{k!} 0_\cN^* [f_t,[\dots,[f_t,g_t]\dots]]=\\&
	= 0_\cN^* \Phi_t^*g_t,
	\end{align*}
where we used that $f_t$ and $\lambda_t$ commute because $C^\bullet_{\cN}$ is abelian for the Schouten bracket. 
Using Proposition \ref{Prop: equivalenceDef},  we conclude that $f_t$ induces a trivial deformation 
of the zero-section. 
\end{proof}

Note that in the general case, i.e. an $n$-shifted cotangent bundle of a degree $n$ manifold, 
we have to be a bit more careful with the gauge equivalence relation of Maurer-Cartan elements (\ref{prop: EquMC}), since the 
$L_\infty$-algebra from  Theorem \ref{Thm: UniqLinfty} controls just the formal deformation problem.
Nevertheless, we have a similar statement to Theorem \ref{Lem:HomPoissEqMC} for the formal 
deformation problem in the spirit of 
Theorem \ref{Thm:formDefProb}: 

\begin{theorem}
Let $i\colon \cL\hookrightarrow (\cM,\omega,X_\theta)$ be a Lagrangian submanifold, let $\Phi\colon T^*[n]\cL\to \cM\at{U}$
be a Lagrangian tubular neighbourhood and let $f_t\in C_{\cL}^{n}$ be a smooth curve of functions with $f_0=0$,
such that $\mathfrak{R}_\Phi(f_t)$ induces a trivial deformation. Then the formal Taylor expansion of $f_t$ 
is a Maurer-Cartan element in $(C_{\cL}[[\nu]],\{l^j\}_{j\in \mathbb{N}})$ gauge equivalent to $0$.  
\end{theorem}

\begin{proof}
Since $f_t$ induces a trivial deformation, we can find $\lambda_t\in C_{\cL}^{n-1}$, such that
	\begin{align*}
	\dot{f}_t =0_\cN^*\{\pi^*\lambda_t, \Phi_t^*\theta\}_{\can}
	\end{align*}	  
by the proof of Theorem \ref{Lem:HomPoissEqMC}. If we Taylor expand both sides of the equation, we get 
	\begin{align*}
	\frac{\D}{\D \nu} \hat{f_t} =\widehat{\dot{f}_t}=  \sum_{k=0}^\infty l^{k+1}(\hat{f_t},\dots,\hat{f_t},\hat{\lambda_t}).
	\end{align*}
\end{proof}

\subsection{Simultaneous Deformations including the $Q$-structure.}\label{sec:sim}
Before we discuss some examples of our deformation theory, we want to briefly discuss the case when we do not 
only deform a Lagrangian $NQ$-submanifold, but also vary the ambient homological vector field. 
Nevertheless, we keep the symplectic form fixed, for deformations of the symplectic form on 
$2$-manifolds see \cite{bof:sym}. For convenience, we only treat the case of 
$(T^*[n]\cN,\omega_{\can},\theta)$ 
where $\cN$ is an $n$-manifold and $\theta\in C_{T^*[n]\cN}^{n+1}$, 
such that $\{\theta,\theta\}_{\can}=0$.
From the discussion in the previous sections it should be clear that, at 
least for small deformations, this is already the most general case. Moreover, 
the proofs are more or less an adaption of the previous sections and are kept short. 

\begin{proposition}\label{Thm: flow QLag}
Let $\cN$ be an $n$-manifold and $\theta_t\in C_{T^*[n]\cN}^{n+1}$  a smooth curve such
 that $\{\theta_t,\theta_t\}_\can=0$. Let $f_t\in C_{\cN}^{n}$ be a smooth curve of functions, then
$\graph(\D f_t)\hookrightarrow (T^*[n]\cN,\omega_{\can},\theta_t)$  is a Lagrangian $NQ$-submanifold 
for all $t$  if and only if $0_\cN^*\Phi_t^*\theta_t=0$ for all $t$, where $\Phi_t$ is the flow 
of $-\{\pi^*\dot{f}_t,\argument\}_\can$. 
\end{proposition}

\begin{proof}
Follows automatically from  Corollary \ref{Prop: DefLagsilly}.  
\end{proof}

The idea is to relate the deformation theory of the Lagrangian $NQ$-submanifold $0_\cN:\cN\hookrightarrow 
(T^*[n]\cN,\omega_\can,\theta_t)$, to the 
$L_\infty$-algebra from Theorem \ref{Thm: DerBrackII}.

\begin{corollary}\label{Cor: Lag+QDef}
Let $\cN$ be an $n$-manifold and let $\theta\in C_{T^*[n]\cN}^{n+1}$ be a function such that 
$\{\theta,\theta\}_\can=0$ and $0_\cN^*\theta=0$. 
Then there is an $L_\infty$-algebra structure on $\mathcal{L}:=C_{\cN}[n-1]
\oplus C_{T^*[n]\cN}[n] $, where the structure maps $L_\theta^i \colon \bigwedge^i 
\mathcal{L}\to \mathcal{L}$ are given by 
\begin{equation*}
    \begin{array}{rcl}
          L^1_{\theta}(g+f)&=&0_\cN^*(f+\{\theta,g\}_{\can})-\{\theta,f\}_{\can}, \\
          L^2_{\theta}(f_1, f_2)&=&-\{f_1,f_2\}_{\can},\\
          L^{k+1}_{\theta}(f, g_1,\dots,g_k)&=&(-1)^{k(|f|-n+1)+ 
	\sum_{i=1}^{k}(k-i)(|g_i|-n)}0_\cN^*\{\cdots\{f,g_1\}_{\can},g_2\}_{\can}\cdots, g_k\}_{\can},\\
        L^k_{\theta}(g_1,\dots, g_k)&=&l^k(g_1,\dots, g_k),
    \end{array}
\end{equation*}
for $g, g_1,\cdots, g_k\in C_\cN[n-1]$ and $f, f_1, f_2\in C_{T^*[n]\cN}[n].$
\end{corollary} 

\begin{proof}
As already mentioned in the proof of Proposition \ref{Thm: weakHomPoiss}, the pair 
$C_{\cN}[n-1]$ and  $ C_{T^*[n]\cN}[n] $ forms a $V$-algebra, so we can apply  
Theorem \ref{Thm: DerBrackII} to obtain a curved $L_\infty[1]$-algebra structure. 
The curvature is given by $0_\cN^*\theta_t=0$. The structure maps above are just the maps from Theorem 
\ref{Thm: DerBrackII} identified via the d\'ecalage isomorphism. 
\end{proof}

Let us now proceed as in Proposition \ref{Prop: Tangent to Lagrangian} and treat the case of 
homotopy Poisson manifolds where there are no convergence issues. 

\begin{theorem}\label{Cor: MCshit}
	Let $(\cN, \theta)$ be a strict degree $n-1$ homotopy Poisson manifold and let 
	$f_t\in C_\cN^{n}$ be given such that $f_0=0$ as well as $\theta_t\in C_{T^*[n]\cN}^{n+1}$ with 
	$\theta_0=0$.  
	Then $Q_t$ is a homological vector field and $$\graph(\D f_t)\hookrightarrow (T^*[n]\cN,\omega_{\can}, 
	Q_t=X_{\theta+\theta_t})$$ is a Lagrangian $NQ$-submanifold
	at time $t$ if and only if $f_t+\theta_t$ is a Maurer-Cartan element in the 
	$L_\infty$-algebra from Corollary \ref{Cor: Lag+QDef} at time $t$. 
\end{theorem}

\begin{proof}
On the one hand, notice that $Q_t$ is a homological vector fields if and only if 
$$\frac{1}{2}\{\theta+\theta_t,\theta+\theta_t\}_{\can}=\{\theta,\theta_t\}_{\can}+\frac{1}{2}\{\theta_t,\theta_t\}_{\can}=0.$$
On the other hand, the flow of $-\{\pi^*\dot{f}_t,\argument\}_{\can}$ is given by $\E^{-\{\pi^*f_t,\cdot\}_{\can}}$, 
as in the proof of Proposition \ref{Prop: Tangent to Lagrangian}. Therefore, using Proposition \ref{Thm: flow QLag}, 
we get that the equation for $\graph(\D f_t)$ being a Lagrangian $NQ$-submanifold inside $(T^*[n]\cN,\omega_{\can}, Q_t)$ is
\begin{align*}
	0=0_\cN^*\E^{-\{f_t,\cdot\}}(\theta+\theta_t)&
	=\sum_{k=0}^{n+1}\frac{(-1)^k}{k!} [f_t,[f_t,[f_t,\cdots,[f_t,\pi^k+\pi^k_t]\cdots]\\&
	=
	\sum_{k=0}^{n+1}\frac{1}{k!}\Big( L^{k}_\theta(f_t,\dots, f_t)
	+L^{k+1}_\theta(\theta_t, f_t,\dots, f_t)\Big),
	\end{align*}
where $\theta=\pi^1+\pi^2+\dots \in \bigoplus_{i=1}^{n+1} \fX^{i,k-in}(\cN)$ and 
$\theta_t=\pi^1_t+\pi^2_t+\dots \in \bigoplus_{i=1}^{n+1} \fX^{i,k-in}(\cN)$.

The Maurer-Cartan equation for the element $f_t+\theta_t$ on the $L_\infty$-algebra of  Corollary \ref{Cor: Lag+QDef} 
reads as follows:
\begin{align*}
    0=&\sum_{k=1}^\infty \frac{1}{k!}L_\theta^k(f_t+\theta_t, \cdots,f_t+\theta_t)\\
    =& 0_\cN^*(\theta_t+\{\theta,f_t\}_{\can})-\{\theta,\theta_t\}_{\can} +\frac{1}{2}( l^2(f_t, f_t)+
    2\ 0^*_{\cN}\{\theta_t, f_t\}_{\can} -\{\theta_t,\theta_t\}_{\can})\\
    &+\sum^{\infty}_{k=3}\frac{1}{k!}(L_\theta^k(f_t,\dots, f_t)+kL_\theta^{k}(\theta_t, f_t,\dots, f_t)).
\end{align*}
According to the splitting $\mathcal{L}=C_{\cN}[n-1]
\oplus C_{T^*[n]\cN}[n]$, we get the two independent equations 
\begin{align*}
    & \{\theta,\theta_t\}_{\can}+\frac{1}{2}\{\theta_t,\theta_t\}_{\can}=0 \quad \text{ and }\\&
    \sum_{k=0}^{n+1}\frac{1}{k!}\Big( L^{k}_\theta(f_t,\dots, f_t)
	+L^{k+1}_\theta(\theta_t, f_t,\dots, f_t)\Big)=0.
\end{align*}
These equations are exactly the ones presented above, therefore the result is proven.
\end{proof}

\begin{remark}\label{rmk:simul}
One can study simultaneous deformations for arbitrary Lagrangian $NQ$-submanifolds 
inside symplectic $NQ$-manifolds. The $L_\infty$-algebra governing this problem will be the same as in the previous theorem. However this will just control the formal deformation 
problem as happens in Theorem \ref{Thm:formDefProb} and convergence issues will appear.
\end{remark}

\section{Examples}\label{sec:ex}
This last section is mainly devoted to give examples of Lagrangian $NQ$-submanifolds and their deformations and relate our
results with previous works. 

\begin{example}[The case $n=1$]
 It is well known, see \cite{Roytenberg, Vor:bia}, that non-graded Poisson manifolds are in one to one 
 correspondence with degree $1$ symplectic $NQ$-manifolds via 
 $$(M, \pi)\rightsquigarrow (T^*[1]M, \omega_{\can}, Q=\{\pi, \cdot\}_{\can}).$$ 
 Using the conormal bundle of Example \ref{ex:conor} 
 we have that coisotropic submanifolds of $(M, \pi)$ are in one to one correspondence 
 with the Lagrangian $NQ$-submanifolds of $(T^*[1]M, \omega_{\can}, Q=\{\pi, \cdot\}_{\can})$, i.e.  
 $$j:C\hookrightarrow (M,\pi)\rightsquigarrow j:N^*[1]C\hookrightarrow (T^*[1]M, \omega_{\can}, \{\pi, \cdot\}_{\can}),$$
which can be seen in \cite{cat:coi}. 
Therefore, studying deformations of degree $1$ Lagrangian $NQ$-submanifolds 
is the same as studying deformations of coisotropic submanifolds inside Poisson manifolds. In this case,  
the $L_\infty$-algebra of Theorem \ref{Thm: UniqLinfty} controlling deformations is the same as the one 
constructed in \cite{CattaneoFelder}. The special case of deformations of coisotropic submanifolds was already discussed in
\cite{oh:coi}.  For the study of simultaneous 
deformations of coisotropic submanifolds of Poisson manifolds one needs the more general version of Theorem \ref{Cor: MCshit} explained in 
Remark \ref{rmk:simul}. This problem was already addressed in \cite{Zambon}.    
\end{example}

\begin{example}[The case $n=2$]\label{ex:n2}
    It was established in \cite{Roytenberg, Severa}   that degree $2$ symplectic $NQ$-manifolds 
    correspond to \emph{Courant algebroids}, see \cite{wei:cou} for definition. Under this 
    correspondence we have that Lagrangian $NQ$-submanifolds become \emph{Dirac structures with support}, 
    see \cite{Support} for a definition. In particular, wide Lagrangian $NQ$-submanifolds are the usual 
    Dirac structures. 

    The $L_\infty$-algebra controlling the deformation problem of a Dirac structure inside a 
    Courant algebroid already appeared in \cite{wei:cou}. However the full deformation problem 
    was just studied recently in \cite{Gualtieri}. Their main Theorem 3.7 is a particular case 
    of our Theorem \ref{Thm: UniqLinfty}, however the methods used are quite different. While 
    we explore the graded geometry side, they use Courant algebroid techniques. As far as we 
    know the deformation theory for Dirac structures with support was not worked out yet. Therefore 
    Theorem \ref{Thm:formDefProb} is a genuine new contribution to the theory of Courant algebroids. 
    The simultaneous deformation of a \emph{generalized complex structure}, a particular case of complex
     Dirac structure, inside an exact Courant algebroid appeared in \cite{Zambon}. Their result can be seen as a 
     particular case of Theorem \ref{Cor: MCshit}.
\end{example}

\begin{example}[Higher Courant algebroids] It was observed by the first author that the degree $n$ 
symplectic $NQ$-manifolds $(T^*[n]T[1]M, \omega_{\can}, Q)$ codify the higher Courant algebroid structures 
on the vector bundle $TM\oplus \wedge^{n-1}T^*M$, see \cite{CuecaCoTan}. Under this correspondence the 
Lagrangian $NQ$-submanifolds become the higher Dirac structures, examples of them include graph of 
$n$-forms and graph of twisted Nambu-structures, see \cite{CuecaCoTan}. Therefore our work presents 
a natural framework to deal with the deformation theory of these objects. More generally one can 
replace $(T[1]M, Q_{dr})$ by any Lie algebroid $(A[1], Q)$ as in \cite{CuecaCoTan}.
\end{example}

\begin{example}[The classical master equation]
    Let $(\cM, \omega)$ be a symplectic $n$-manifold and denote the associated Poisson bivector 
    field by $\pi=\omega^{-1}\in \fX^{2,-n}(\cM)$. Since $\pi$ is Poisson we get that 
    $(T^*[n+1]\cM,\omega_{\can}, X_\pi=\{\pi, \cdot\}_{\can})$ is a degree $n+1$ symplectic 
    $NQ$-manifold with the zero-section as a Lagrangian $NQ$-submanifold
    $$0_\cM:\cM\hookrightarrow(T^*[n+1]\cM,\omega_{\can}, X_\pi=\{\pi, \cdot\}_{\can}).$$
    Using Proposition \ref{Prop: Tangent to Lagrangian}, small deformations of the 
    zero-section as a Lagrangian $NQ$-submanifold are in one to one correspondence with functions 
    $$f\in C^{n+1}_\cM\quad \text{satisfying}\quad \{f,f\}=0,$$
    this equation is known as the \emph{classical master equation}, see e.g. \cite{AKSZ}.
\end{example}
The next three examples are variations of the previous one:

\begin{example}[Q-deformations as Lagrangian deformations]
    Let $(\cM, \omega, Q=\{\alpha, \cdot\})$ be a degree $n$ symplectic $NQ$-manifold 
    with $\alpha\in C^{n+1}_\cM$ satisfying the classical master equation 
    $\{\alpha, \alpha\}=0$. As before denote $\pi=\omega^{-1}\in \fX^{2,-n}(\cM)$ and 
    consider the degree $n+1$ symplectic $NQ$-manifold $(T^*[n+1]\cM, \omega_{\can}, X_{Q+\pi})$ 
    with the zero-section as a Lagrangian $NQ$-submanifold 
    $$0_\cM:\cM\hookrightarrow(T^*[n+1]\cM,\omega_{\can}, X_{Q+\pi}=\{Q+\pi, \cdot\}_{\can}).$$
    Again, by Proposition \ref{Prop: Tangent to Lagrangian}, small deformations of the zero-section 
    as a Lagrangian $NQ$-submanifold are in one to one correspondence with functions 
    $$f\in C^{n+1}_\cM\quad \text{satisfying}\quad 0=Q(f)+\frac{1}{2}\{f,f\}=\frac{1}{2}\{\alpha+f, \alpha+f\}.$$
    Thus the function $f$ provides a deformation of the homological vector field $\widehat{Q}=\{\alpha+f,\cdot\}=Q+X_f.$ 
    
    An example of this situation is given by \emph{exact Courant algebroids}, see \cite{Sev:let}, 
    as a deformation of the standard Courant algebroid by using a $3$-form. In this case we see 
    that two such deformations are equivalent if and only if the $3$-forms are cohomologous. 
    Thus we recover the classification of exact Courant algebroids as obtained in \cite{Sev:let}
     by  Lagrangian deformation theory.
\end{example}

\begin{example}[Twisted Courant algebroids as deformations of Courant algebroids]\label{ex:tca}
This is a generalization of Example 4.2 in \cite{ike:can}. They where interested in boundary 
conditions for AKSZ sigma models. Here we reinterpret their computations as a deformation 
of a Lagrangian $NQ$-submanifold inside of a degree $3$ symplectic $NQ$-manifold. 

As explained in Example \ref{ex:n2}, Courant algebroids are in correspondence with degree $2$ symplectic 
$NQ$-manifolds $(\cM, \omega, Q=X_\alpha)$ where $\alpha\in C^3_{\cM}$. Pick a degree $2$ 
symplectic manifold  $(\cM, \omega)$. Since $\cM$ is symplectic we have that $T[1]\cM$ is a symplectic $3$-manifold with 
symplectic form given by the tangent lift $\omega^T$. Moreover, the map 
$$\omega^\flat:(T[1]\cM, \omega^T)\to(T^*[3]\cM, \omega_{\can})$$ 
is a symplectomorphism. By construction of $T[1]\cM$ we have that $\Omega^4(|\cM|)\subset C^4_{T[1]\cM}$ . 
Thus given $H\in\Omega^4(|\cM|)$ we can consider the vector field 
$$Q=Q_{dr}+\{H, \cdot\}_{T}\in \fX^{1,1}(T[1]\cM)\quad\text{and}\quad Q^2=0\Leftrightarrow dH=0.$$ 
In this case, $0_{\cM}:\cM\hookrightarrow(T[1]\cM, \omega^T, Q)$ is a Lagrangian $NQ$-submanifold. 
Using $\omega^\flat$ as a tubular neighbourhood, we choose $\alpha\in C^3_{\cM}$ and consider the deformation 
$$\E^{-\{\pi^*\alpha, \cdot\}_{T}}\circ 0_\cM:\cM\hookrightarrow T[1]\cM,$$
where $\pi:T[1]\cM\to \cM$ is the projection. Example 4.2 in \cite{ike:can} shows that this is a 
deformation of the zero section as a Lagrangian $NQ$-submanifold if and only if $(\cM, \omega, \alpha)$ 
is a \emph{twisted Courant algebroid} as introduced in \cite{thom:cou}. It is easy to see that in this 
example the $L_\infty$-algebra whose Maurer-Cartan elements are deformations is an honest $L_4$-algebra. 
For the coordinate expression of the above computation see also \cite{ike:can}.
\end{example}

\begin{example}[The Casimir as an obstruction for deformation] 
This example corresponds to $\S$7.2.1 in \cite{eck:clif}. They study the 
Poisson structure on the Weil algebra and 
the cubic Dirac operator. We reinterpret their computations as an obstruction to deform a Lagrangian 
$NQ$-submanifold inside of a degree $3$ symplectic $NQ$-manifold.

Let $(\mathfrak{g}, [\cdot, \cdot], B(\cdot, \cdot))$ be a quadratic Lie algebra, i.e. a Courant 
algebroid over a point. Then $(\mathfrak{g}[1], \omega, X_\phi)$ is a degree $2$ symplectic 
$NQ$-manifold with $\phi\in C^3_{\mathfrak{g}[1]}$. Using the pairing we identify 
$B^\flat:\mathfrak{g}\overset{\sim}{\to}\mathfrak{g}^*.$ The Weil algebra of $\mathfrak{g}$ is 
given by the $2$-manifold $\cN=T[1]\mathfrak{g}[1]=\mathfrak{g}[1]\oplus \mathfrak{g}[2]$ and 
carries a degree $-2$ Poisson structure $\pi\in \fX^{2,-2}(\cN)$ defined by
$$\{\xi_1, \xi_2\}_{\pi}=2B(\xi_1, \xi_2), \quad \{\widehat{\xi}_1, \widehat{\xi}_2\}_{\pi}=\widehat{[\xi_1, \xi_2]}, \quad \{\xi_1, \widehat{\xi}_2\}_{\pi}=0,$$
where $\xi_i\in \mathfrak{g}[1]$ and $\widehat{\xi}_i\in\mathfrak{g}[2].$ So $(\cN, \pi)$ is 
a strict degree $2$ homotopy Poisson manifold and  $(T^*[3]\cN, \omega_{\can}, X_{\pi})$ is a 
degree $3$ symplectic $NQ$-manifold with the Lagrangian $NQ$-submanifold 
$$0_\cN:\cN\to (T^*[3]\cN, \omega_{\can}, X_{\pi}).$$  
Since $X_\pi\at{\cN}=0$ we have that any element $f\in C^3_\cN$ gives an infinitesimal 
deformation of the zero-section as introduced in $\S$\ref{sec:eq-def}.  
If we have a smooth curve of functions $f_t\in C_{\cN}^{3}$, which deforms
$\cN$, then 
	\begin{align*}
	\tfrac{1}{2}\{f_t,f_t\}_\pi=0
	\end{align*}	 
and $f_0=0$.
If we derive this equation twice along $t$ and evaluate at $t=0$, we get that $\{\dot{f}_0,\dot{f}_0\}_\pi=0$. So, if we 
are searching for a deformation $f_t$ such that $\dot{f}_0$ is prescribed,
then $\{\dot{f}_0,\dot{f}_0\}_\pi=0$.  

In this particular 
example, to see if an infinitesimal deformation $D$ has the chance to extend  to a full deformation  we compute $\{D, D\}_{\pi}\in C^4_{\cN}$. The assignment $f\rightsquigarrow [\{f, f\}]$ is  called the \emph{Kuranishi map} \cite[$\S 4$]{man:def}.

Consider the element 
$$D=\sum_{i} \widehat{e}_i e^i+ \phi\in C^3_\cN,$$
where $\{e^i\}$ is a basis of $\mathfrak{g}[1], \{\widehat{e}_i\}$ basis of $\mathfrak{g}[2]$ and
 $\phi\in C_{\mathfrak{g}[1]}^3$ as in the beginning. Then Proposition 7.5 in \cite{eck:clif} shows that
$$\{D, D\}_\pi=2\sum_{i}\widehat{e}_i\widehat{e}^i \in C^4_{\cN}$$ that is the \emph{quadratic Casimir element.} 
Thus we conclude that the quadratic Casimir element gives an obstruction to deforming the zero-section of 
$ (T^*[3]\cN, \omega_{\can}, X_{\pi})$  as a Lagrangian $NQ$-submanifold in the direction of the element $D$. 
\end{example}

\begin{appendices}
\section{Reminder on $L_\infty$-algebras, Maurer-Cartan elements and derived brackets}\label{app:Linfty}
There are a lot of definitions of $L_\infty$-algebras in the literature: some prefer to 
use the anti-symmetric view-point, other the symmetric one.  
Since we are using both, let us discuss all of them briefly. This appendix is far from
 being complete, it is meant to fix notations and introduce 
the sign conventions we are using. For an introduction of (curved) $L_\infty$-algebras, 
we refer the reader to \cite{KraftSchnitzer}.

\begin{definition}
  A \emph{curved $L_\infty[1]$-algebra} is a $\bZ$-graded vector space $L$ 
	endowed with a degree one codifferential $NQ$ on the 
	symmetric coalgebra $(\sym(L),\Delta_\mathrm{sh})$.  $Q(1)\in L^1$ is called the curvature, 
	if it vanishes $(L,Q)$ is called an $L_\infty[1]$-algebra. 
  An $L_\infty[1]$-morphism between two curved $L_\infty[1]$-algebras $F\colon (L,Q) 
  \rightarrow (L',Q')$ is a morphism of graded coalgebras
  \begin{equation}
    F \colon 
		\sym(L)
		\longrightarrow
		\sym(L')
  \end{equation}
	such that $F(1)=1$ and $F \circ Q = Q' \circ F$.
\end{definition}
Alternatively, one represent the codifferential $NQ$ as a sequence of maps 
	\begin{align*}
	Q_k\colon \sym^k(L)\to L[1],
	\end{align*}
known as the Taylor coefficients, 
given by the appropriate projections in the direct sum, that exists since the coalgebra $\sym(L)$ 
is cofree in some category. The condition $Q^2=0$ reflects also in 
a quadratic relation between the $Q_k$.  Moreover, we can represent an 
$L_\infty$-morphism $F\colon (L,Q) 
  \rightarrow (L',Q')$ also by its Taylor coefficients
	\begin{align*}
	F_k\colon \sym^k L  \rightarrow L'
	\end{align*}
fulfilling certain compatibilities with the Taylor coefficients of $NQ$ and $Q'$. 

An $L_\infty$-algebra is just an $L_\infty[1]$-algebra on a shifted vector space. 
To be more precise: 

\begin{definition}
A \emph{curved $L_\infty$-algebra $(\liealg{g},Q)$} is a $\bZ$-graded vector space $\liealg{g}$ 
endowed with a degree one codifferential $NQ$ on the 
	symmetric coalgebra $(\sym(\liealg{g}[1]),\Delta_\mathrm{sh})$. The \emph{curvature} is $Q(1)\in\liealg{g}^2$  and if it vanishes, $(\liealg{g},Q)$ 
	is called an \emph{$L_\infty$-algebra}. A morphism between two $L_\infty$-algebras is defined to be a 
	morphism of the corresponding $L_\infty[1]$-algebras.  
\end{definition}

For a graded vector space $V$, there is the so-called  d\'{e}calage isomorphism 
$\text{d\'{e}c}\colon \sym^k V[1]\to (\bigwedge^kV)[k]$
given by 
	\begin{align*}
	\text{d\'{e}c}(x_1\vee\dots\vee x_k):= 
	(-1)^{\sum_{i=1}^k(k-i)(|x_i|-1)}x_1\wedge\dots\wedge x_k
	\end{align*}
for all $x_i\in V[1]^{|x_i|-1}=V^{|x_i|}$. Using this, one immediately sees that 
we can also define an $L_\infty$-algebra structure on $\liealg{g}$  as a sequence of maps
 $\ell^k\colon \bigwedge^k\liealg{g}\to \liealg{g}[2-k]$
	for $k\geq 0$, called the Taylor coefficients, such that 
	\begin{align*}
	\sum_{k=1}^n  \sum_{\sigma\in Sh(k,n-k)} 
		\chi(\sigma) 
		(-1)^{k(n-k)}\ell^{n-k+1}(\ell^k(x_{\sigma(1)},\cdots, x_{\sigma(k)}),
		x_{\sigma(k+1)}, \dots ,	x_{\sigma(n)})=0
	\end{align*}
for all $x_1,\dots,x_n \in \liealg{g}$, where $\chi(\sigma)=\sign(\sigma)\epsilon(\sigma)$.

\begin{definition}\label{Def: MC}
Let $(\liealg{g},\{l^k\}_{k\in\mathbb{N}_0})$ be an $L_\infty$-algebra. An element 
$\pi\in \liealg{g}^1$ is called a \emph{Maurer-Cartan} element if 
	\begin{align*}
	\sum_{k=0}^\infty \frac{1}{k!} l^k(\overbrace{\pi,\dots, \pi}^{k\text{-times}})=0,
	\end{align*}
where we assume the sum converges in some sense. The set of Maurer-Cartan elements of a 
curved  $L_\infty$-algebra $\mathfrak{g}$ is denoted by $\mathrm{MC}(\mathfrak{g})$.
\end{definition}

\begin{lemma}
Let $(\liealg{g},\{l^k\}_{k\in\mathbb{N}_0})$ be a curved $L_\infty$-algebra 
and let $\pi\in\liealg{g}^1$. The maps 
	\begin{align*}
	l_\pi^k(x_1,\dots, x_k)
	:=\sum_{i=0}^\infty \frac{1}{i!} l^{k+i}
	(\overbrace{\pi,\dots\pi }^{i\text{-times}},x_1,\dots, x_k)
	\end{align*}
are Taylor coefficients of a curved $L_\infty$-algebra, $\liealg{g}^\pi$, with curvature 
$\sum_{k=0}^\infty \frac{1}{k!} l^k(\overbrace{\pi,\dots,\pi}^{k\text{-times}})$.
In particular, if $\pi\in\mathrm{MC}(\liealg{g})$ then $\liealg{g}^\pi$ is an $L_\infty$-algebra.
\end{lemma}

On the set of Maurer-Cartan elements there is a notion of equivalence:

\begin{definition}\label{prop: EquMC}
Let  $(\liealg{g},\{l^k\}_{k\in\mathbb{N}_0})$ be a curved $L_\infty$-algebra and let 
$\pi_0,\pi_1\in \mathrm{MC}(\liealg{g})$. The elements $\pi_0$ and $\pi_1$ are called \emph{gauge equivalent} if there 
are smooth curves $\pi\colon [0,1]\to \mathrm{MC}(\liealg{g})$ and 
$\lambda\colon [0,1]\to \liealg{g}^0$, such that 
	\begin{align*}
	\begin{cases}
	\pi(0)=\pi_0 \text{ and } \pi(1)=\pi_1 \\
	\frac{\D}{\D t} \pi(t)= l_{\pi(t)}^1(\lambda(t))
	\end{cases}.
	\end{align*}
\end{definition}

\begin{remark}
The term "smooth" in Definition \ref{prop: EquMC} can be interpreted in many ways, but it is
 necessary to specify since we need a derivative in the required equations. 
In the purely algebraic setting one usually takes polynomials, but in the setting of the paper,
we have a good notion of smoothness. Moreover, if one has a version of the Picard-Lindelöf theorem one can consider the initial value problem 
  \begin{align*}
	\begin{cases}
	\pi(0)=\pi \\
	\frac{\D}{\D t} \pi(t)= l_{\pi(t)}^1(\lambda(t))
	\end{cases}
	\end{align*}
for all $\pi\in\mathrm{MC}(\liealg{g})$ and all 
smooth paths $\lambda\colon [0,1]\to \liealg{g}^0$ and evaluate the result 
at $t=1$. One can show that this yields a Maurer-Cartan element. The gauge action of $\lambda$ on $\pi$ is now defined as 
	\begin{align*}
	\lambda\acts \pi=\pi(1).
	\end{align*}
For more details see \cite{KraftSchnitzer}.
\end{remark}

In the following we recall the construction of derived brackets from Voronov \cite{Voronov}.

\begin{definition}\label{def:Valgebra}
    A \emph{$V$-algebra} is a triple  $(\liealg{g},\liealg{a},p)$, where $(\liealg{g}, [\argument,\argument])$ is a graded Lie algebra, $\iota\colon\liealg{a}\hookrightarrow \liealg{g}$ is an abelian subalgebra an $p\colon \liealg{g}\to \liealg{a}$ is a projection such that
    $$p\circ \iota=\id \quad \text{and}\quad p[x,y]=p[\iota(p(x)),y]+p[x,\iota(p(y))].$$
\end{definition}

\begin{theorem}[\cite{Voronov}, Theorem 1]\label{Them: DerBra}
Let $(\liealg{g},\liealg{a},p)$ be a $V$-algebra and let $\Delta\in\liealg{g}^1$ be given such 
that $[\Delta,\Delta]=0$. The maps $Q^i_\Delta\colon \sym^i\liealg{a}\to \liealg{a}$ given by
	\begin{align*}
	Q^i_\Delta(a_1,\dots, a_i)= p[\cdots[\Delta,\iota(a_1)],\iota(a_2)],\cdots],\iota(a_i)]
	\end{align*}
are the Taylor coefficients of a curved $L_\infty[1]$-algebra structure on $\liealg{a}$ with 
curvature $p(\Delta)$.	  
\end{theorem}
Since there is an equivalence of $L_\infty$-algebra structures on a vector space $V$ and 
$L_\infty[1]$-algebra structures on $V[1]$ via the d\'{e}calage-isomorphism, we immediately 
get 

\begin{corollary}
Let $(\liealg{g},\liealg{a},p)$ be a $V$-algebra and let $\Delta\in\liealg{g}^1$ be given 
such that $[\Delta,\Delta]=0$. The maps $l^i_\Delta\colon \bigwedge^i \liealg{a}[-1]\to 
\liealg{a}[1-k]$ given by 
   \begin{align*}
	l^i_\Delta(a_1, \dots , a_i)=(-1)^{\sum_{j=1}^i(i-j)|a_j|}
	 p[\cdots[\Delta,\iota(a_1)],\iota(a_2)],\cdots],\iota(a_i)]
	\end{align*}
are the Taylor coefficients of a curved $L_\infty$-algebra structure on $\liealg{a}[-1]$ 
with curvature $p(\Delta)$, where $a_i\in \liealg{a}^{|a_i|}=\liealg{a}[-1]^{|a_i|+1}$.
\end{corollary}

\begin{theorem} [\cite{VoronovII}, Theorem 2]\label{Thm: DerBrackII}
Let $(\liealg{g},\liealg{a},p)$ be a $V$-algebra and $\Delta\in\liealg{g}^1$ such that $[\Delta,\Delta]=0$. Denote by $(\liealg{a},\{Q_\Delta^i\}_{i\in\mathbb{N}})$ the 
corresponding curved $L_\infty[1]$-algebra of Theorem \ref{Them: DerBra}. Then the direct sum 
$\liealg{h}:=\liealg{a}\oplus\liealg{g}[1] $ carries a curved $L_\infty[1]$-algebra structure 
given by the Taylor coefficients  $R_\Delta^{i}\colon \sym^i\liealg{h}\to \liealg{h}$ which 
are defined by 
	\begin{itemize}
	\item $R^1_\Delta(a+x)=p([\Delta,a]+x)-[\Delta,x]$ for $a+x\in\liealg{a}\oplus\liealg{g}[1]$,
	\item $R^2_\Delta(x, y)=-(-1)^{|x|}[x,y]$ for $x,y\in \liealg{g}[1]$ 
	(to be more precise $x\in \liealg{g}[1]^{|x|}=\liealg{g}^{|x|+1}$),
	\item $R^{k+1}_\Delta(x, a_1,\dots, a_k)=p[\cdots[x,a_1],a_2,\cdots],a_k]$ 
	for $x\in\liealg{g}[1]$ and $a_1,\dots,a_k\in \liealg{a}$ and 
	\item $R^k_\Delta(a_1,\dots, a_k)=Q_\Delta^k(a_1,\dots,a_k)$ for  $a_1,\dots,a_k\in \liealg{a}$.
	\end{itemize}		 
\end{theorem}

\begin{remark}
One can obtain the $L_\infty$-algebra structure in a different way: one constructs 
it for $\Delta=0$ with $[\Delta,\Delta]=0$ and then twists the structure by 
$\Delta\in \liealg{h}^1$. It is clear now that we started with an $L_\infty$-algebra 
without curvature and get one with curvature $p(\Delta)+\frac{1}{2}[\Delta,\Delta]=p(\Delta)$. 
Using this, one proves that for two gauge equivalent $\Delta_0,\Delta_1\in \mathrm{MC}(\liealg{h})$ 
there is an $L_\infty$-isomorphism between 
$(\liealg{h},\{R^i_{\Delta_0}\}_{i\in \mathbb{N}})$ and $(\liealg{h},
\{R^i_{\Delta_1}\}_{i\in \mathbb{N}})$, which is proven in \cite{KraftSchnitzer}. One can 
additionally show that this constructed $L_\infty$-morphism restricts to $\liealg{a}$ and 
coincides with the one constructed in \cite{SchaetzCattaneo}.
 
\end{remark}
\end{appendices}
\bibliographystyle{nchairx}

\end{document}